\documentclass[11pt]{article}

\usepackage[margin=1in]{geometry}
\usepackage[utf8]{inputenc}
\usepackage{amsmath,amssymb,amsfonts,amsthm} 
\usepackage{array}
\usepackage{booktabs}
\usepackage{float}
\usepackage{enumerate}
\usepackage{mathtools}
\usepackage{mathrsfs}
\usepackage{xcolor}

\usepackage[hidelinks]{hyperref}
\newtheorem{theorem}{Theorem}[section]
\newtheorem{proposition}[theorem]{Proposition}
\newtheorem{lemma}[theorem]{Lemma}
\newtheorem{corollary}[theorem]{Corollary}

\theoremstyle{definition}
\newtheorem{definition}[theorem]{Definition}
\newtheorem{example}[theorem]{Example}

\newcommand{\examplecontinuationname}{}
\newtheorem*{examplecontinuationblock}{\examplecontinuationname}
\newenvironment{examplecontinuation}[2]{%
  \renewcommand{\examplecontinuationname}{Example~\ref*{#1}}%
  \begin{examplecontinuationblock}[#2]%
  \leavevmode\phantomsection\ignorespaces
}{%
  \end{examplecontinuationblock}%
}

\theoremstyle{remark}
\newtheorem{remark}[theorem]{Remark}
\newif\ifextendedproofs
\extendedproofsfalse

\newcommand{\inner}[2]{\langle#1,#2\rangle}      
\newcommand{\R}{\mathbb{R}}        

\newcommand{\norm}[1]{\|#1\|}                    
\newcommand{\T}{\top\hspace{-1pt}}               
     
\newcommand{\D}{\operatorname{Diag}}             

\newsavebox{\firstgeometrytablebox}
\newlength{\firstgeometrytablewidth}

\usepackage{tikz}
\usetikzlibrary{arrows.meta, positioning}

\definecolor{myred}{RGB}{255,50,50}         
\definecolor{myblack}{RGB}{0,0,0}           

\definecolor{revisiongreen}{RGB}{0,120,60}

\title{A Smooth Equality Reformulation of Nonlinear
Exponential Cone Optimization}

\author{
Huimin Li\thanks{Graduate School of Informatics, Kyoto University,
Yoshida-Honmachi, Sakyo-ku, Kyoto \mbox{606--8501}, Japan
(\texttt{hm.li@amp.i.kyoto-u.ac.jp}, \texttt{ellen@i.kyoto-u.ac.jp}).}
\and
Ellen H. Fukuda\footnotemark[1]
}

\date{}

\begin{document}

\maketitle

\begin{abstract}

\noindent Exponential cone optimization has received considerable attention in the modern literature, but studies involving nonlinear objective and constraint functions are still limited. This work can be seen as one of the first steps toward understanding nonlinear exponential cone optimization (NECO) problems. Here, we ask the question ``Can smooth equality constraints replace nonlinear exponential cone constraints?''. 
From a practical point of view, even squared slack reformulations in nonlinear programming introduce difficulties. However, such reformulations can be useful for studying optimality conditions. Motivated by this observation, in this work, we use a smooth surjective parameterization of the exponential cone to reformulate NECO as an equality-constrained nonlinear program (RNLP). We then analyze Karush--Kuhn--Tucker (KKT)
and second-order optimality conditions for NECO through RNLP.
Although feasibility is represented exactly, the KKT
conditions for RNLP do not in general guarantee those for NECO.
We characterize this correspondence through multiplier restrictions that are equivalent to a second-order condition involving only slack variables.
We also establish an exact curvature identity showing that RNLP reproduces the sigma term of NECO without additional assumptions.
Finally, we prove the relationships between second-order necessary and sufficient conditions for the two problems. We identify when RNLP can be used to rule out local minimizers of NECO or establish strict local optimality for NECO without explicitly evaluating the sigma term.

\end{abstract}

\medskip
\noindent\textbf{Keywords.}
Exponential cone, equality reformulation, nonlinear conic optimization,
second-order optimality conditions, sigma term, slack variables.

\medskip
\noindent\textbf{MSC codes.}
90C30, 90C25, 49K10.

\section{Introduction}

Let \(\mathcal K_{\exp}\subseteq\mathbb R^3\) denote
the exponential cone defined by
\[
\mathcal K_{\exp}
\coloneqq
\{y\in\mathbb R^3\colon y_3>0,\
y_2\geq y_3\exp(y_1/y_3)\}
\cup
\{(y_1,y_2,0)\colon y_1\leq0,\ y_2\geq0\}.
\]
For a positive integer \(p\), let
\(\mathcal K_{\exp}^p\) denote the Cartesian product
of \(p\) cones \(\mathcal K_{\exp}\).
In this paper, we consider the following nonlinear
exponential cone optimization (NECO) problem:
\begingroup
\renewcommand{\theHequation}{NECO}
\begin{equation}
\begin{array}{ll}
\underset{x\in\mathbb R^n}{\operatorname{minimize}}
& f(x)\\[1mm]
\operatorname{subject\ to}
& g(x)\in\mathcal K_{\exp}^p,
\end{array}
\tag{NECO}
\label{prob:neco}
\end{equation}
\endgroup
where \(f\colon\mathbb R^n\to\mathbb R\) and
\(g\colon\mathbb R^n\to\mathbb R^{3p}\) are twice
continuously differentiable.
We investigate how a smooth equality-constrained reformulation
can be used to analyze optimality conditions for NECO. We note that equality constraints can be also included here, but we omit them for simplicity.

The exponential cone gives compact conic representations
of models involving exponential, logarithmic, entropy,
and relative-entropy functions; see, for example,
\cite[Chapter~5]{MOSEKCookbook2025} and
\cite{ChandrasekaranShah2017}.
It is supported by modern conic optimization solvers, and primal-dual interior-point methods have been developed for
NECO when $f$ and $g$ are linear~\cite{DahlAndersen2022}.
These developments motivate an optimality theory for
nonlinear problems with exponential-cone constraints.
Direct second-order conditions for such problems have
recently been derived from the tangent, normal, and
second-order tangent sets of
\(\mathcal K_{\exp}\)~\cite{Xu2026Thesis}.
These geometric formulas depend on the location of
the point within the cone, particularly on whether it lies on
a curved or flat part of its boundary.

The KKT conditions for NECO require the multiplier
to belong to a normal cone.
On the other hand, second-order optimality conditions involve
the so-called sigma term, which characterizes the curvature of the exponential cone through its second-order tangent sets.
Evaluating these conditions requires geometric information
beyond derivatives of the objective and constraint functions.
We ask whether a smooth equality reformulation can express
these conditions through derivatives of its Lagrangian,
and which optimality properties are preserved in doing so.

An approach with a long history in nonlinear programming
(NLP) is to replace inequality constraints by smooth
equality constraints with additional squared-slack variables
\cite{Spedicato1982,Tapia1974,Tapia1980}.
Armand and Orban~\cite{ArmandOrban2012} examined numerical
difficulties associated with this transformation.
Jongen and Stein~\cite{JongenStein2003} show that
squared-slack reformulations can admit KKT points that
do not correspond to KKT points of the original problem
and can complicate the feasible-set geometry.
Their results illustrate why an exact equality reformulation
need not make a problem easier to solve.
Fukuda and Fukushima~\cite{FukudaFukushima2017NLP}
studied the correspondence between KKT points and
second-order optimality conditions of the original and
reformulated problems.
For an equality-constrained nonlinear program,
first-order conditions consist of feasibility and
stationarity equations.
Its standard second-order optimality conditions use
the Hessian of the Lagrangian on the linearized feasible
space, with no separate sigma term.
A correspondence between the two formulations can therefore
provide an alternative way to derive or verify conic
optimality conditions.
Whether local minimizers and optimality conditions are
preserved depends on the choice of parameterization.
Levin, Kileel, and Boumal~\cite{LevinKileelBoumal2025}
provide general criteria for these correspondences.
Related work also examines how parameterizations affect
convergence. Ouyang et al.~\cite{OuyangLiuPongWang2025}
study the Kurdyka--{\L}ojasiewicz (KL) exponent under
Hadamard parametrization of \(\ell_1\)-regularized problems
and use it to derive local convergence rates for a gradient
method under suitable assumptions.
Ouyang~\cite{Ouyang2025Square} studies the corresponding
exponent for square transformations of objectives comprising
a smooth term and a polyhedral term.

Beyond NLP, squared-slack
reformulations have been studied for nonlinear second-order
cone programming (NSOCP)~\cite{FukudaFukushima2016NSOCP}
and nonlinear semidefinite programming (NSDP)~\cite{LourencoFukudaFukushima2017NSDP}.
More generally, every symmetric cone is a cone of squares
in a Euclidean Jordan algebra, which provides a unified
equality-reformulation framework for nonlinear
symmetric-cone problems
\cite{LourencoFukudaFukushima2018SymmetricCone}.
These studies show that second-order conditions for the
equality reformulation can recover conic multiplier
conditions that its first-order conditions alone
do not guarantee.
They also relate the quadratic terms involving slack
variables to the curvature terms in conic second-order
conditions.
In NLP, Fukuda and Fukushima
\cite{FukudaFukushima2017NLP} observe that the second-order
sufficient condition is stronger than necessary for
recovering multiplier signs.
This distinction motivates studying multiplier recovery
separately from second-order sufficiency.
Ding and Wright further established second-order
correspondences for scalar squared-variable reformulations
\cite{DingWright2025}, as well as second-order and
local-minimizer correspondences for matrix-factor
reformulations of nonlinear semidefinite constraints
\cite{DingWright2025NSDP}.
Their necessary-condition correspondences concern weak
second-order conditions for the original problems,
which test a subspace rather than the full critical cone.

The exponential cone, however, does not admit a square
representation in a Euclidean Jordan algebra.
More generally, it admits no exact global polynomial
slack representation with finitely many unconstrained
variables
\cite{BlekhermanParriloThomas2013,LindstromLourencoPong2023}.
An exact smooth slack representation therefore requires
a nonpolynomial mapping.
To obtain such a representation, we propose the following
smooth surjection
\[
\Phi(a,b,c,r)
\coloneqq
(a^2r-b^2,a^2\exp(r)+c^2,a^2)
=
a^2(r,\exp(r),1)+(-b^2,c^2,0).
\]
We denote its blockwise extension by \(\Phi_p\).
The term \(a^2(r,\exp(r),1)\) is a nonnegative multiple
of a curved-boundary point of the exponential cone, and \((-b^2,c^2,0)\)
belongs to its flat face.
Introducing \(z\in\mathbb R^{4p}\), we consider the
equality-constrained reformulation
\begingroup
\renewcommand{\theHequation}{RNLP}
\begin{equation}
\tag{RNLP}
\label{prob:RNLP}
\begin{array}{ll}
\underset{x\in\mathbb R^n,\ z\in\mathbb R^{4p}}
{\operatorname{minimize}}
& f(x)\\[0.2em]
\operatorname{subject\ to}
& g(x)-\Phi_p(z)=0.
\end{array}
\end{equation}
\endgroup

For a feasible point \(x\), any slack vector \(z\)
satisfying \(g(x)=\Phi_p(z)\) is called a \emph{preimage}
of \(g(x)\).
The same cone point may have more than one preimage by the surjectivity of \(\Phi_p\).
Moreover, the surjectivity of \(\Phi_p\) preserves feasibility exactly,
but does not by itself ensure that optimality conditions
transfer between the two formulations.
In particular, the KKT conditions for RNLP do not,
in general, imply the KKT conditions for NECO with
the same multiplier. 
We then study here how local minimizers, KKT points, and
second-order optimality conditions of NECO and RNLP are related. These results do not necessarily establish faster
evaluation of optimality conditions or convergence
of a numerical method.
Figure~\ref{fig:results-overview} summarizes the results of the paper, with the main relationships between NECO and its equality reformulation.

\begin{figure}[!ht]
\centering
\ifdefined\tikzpicture
  \let\overviewstart\relax
  \let\overviewfinish\relax
\else
  \documentclass[11pt]{article}
  \usepackage[margin=1in]{geometry}
  \usepackage{amsmath,amssymb,xcolor,tikz}
  \usetikzlibrary{arrows.meta}
  \usepackage[active,tightpage]{preview}
  \PreviewEnvironment{tikzpicture}
  \setlength\PreviewBorder{4pt}
  \def\overviewstart{\begin{document}}
  \def\overviewfinish{\end{document}}
\fi
\overviewstart
\begingroup
\definecolor{overviewNECO}{RGB}{148,53,43}
\definecolor{overviewRNLP}{RGB}{0,92,135}
\begin{tikzpicture}[
  x=\linewidth,y=1mm,
  font=\small,
  statement/.style={
    draw=black!40, line width=.45pt, rounded corners=4pt,
    align=center, inner xsep=4pt, inner ysep=6pt,
    text width=.33\linewidth, minimum height=10mm
  },
  neco/.style={statement,draw=overviewNECO!65,fill=overviewNECO!3},
  rnlp/.style={statement,draw=overviewRNLP!65,fill=overviewRNLP!3},
  note/.style={font=\footnotesize,align=center,inner sep=0pt},
  implication/.style={-{Latex[length=1.8mm,width=1.4mm]},line width=.65pt},
  equivalence/.style={Latex-Latex,line width=.65pt}
]

\node[note] at (.5,-3)
  {Fix \(g(\bar x)=\Phi_p(\bar z)\);\\ use the same multiplier \(\lambda\).};

\node[neco,minimum height=32mm,label distance=3mm,
  label={[font=\Large\bfseries,text=overviewNECO,fill=white,
    inner xsep=4pt,inner ysep=5pt]above:NECO}]
  (nkkt) at (.18,-23)
  {\textbf{KKT point} \((\bar x,\lambda)\)};
\node[rnlp,minimum height=32mm,label distance=3mm,
  label={[font=\Large\bfseries,text=overviewRNLP,fill=white,
    inner xsep=4pt,inner ysep=5pt]above:RNLP}]
  (rkkt) at (.82,-23)
  {KKT point \((\bar x,\bar z,\lambda)\)\\[1pt]
   {\footnotesize
    \( + \) second-order slack inequalities\\
    \textbf{at every preimage}}
   \par\smallskip
   {\footnotesize
    Minimum for uniform recovery\\
    from the slack equation:\\[2pt]
    \(\boldsymbol{0,\ 1,\ +\infty}\)
    \quad(\(\Phi\), one block)}};
\draw[equivalence] (nkkt.east)--(rkkt.west);

\draw[black!25] (0,-44)--(1,-44);

\node[neco,minimum height=14mm] at (.18,-55) {\textbf{Conic sigma term}};
\node[rnlp,minimum height=14mm] at (.82,-55)
  {\textbf{Supremum} of the slack\\quadratic term};
\node[font=\Large] at (.5,-55) {\(=\)};
\node[note,text width=.28\linewidth] at (.5,-60)
  {KKT multiplier of NECO};
\draw[black!25] (0,-66)--(1,-66);

\node[neco,minimum height=14mm] (exclude) at (.18,-77)
  {\(\bar x\) is \textbf{not a local minimizer}};
\node[rnlp,minimum height=14mm] (fails) at (.82,-77)
  {Second-order necessary\\condition \textbf{fails}
   at \((\bar x,\bar z,\lambda)\)};
\draw[implication] ([yshift=-3mm]fails.west)--
  node[note,above=2pt,text width=.27\linewidth]
    {{KKT point of NECO;\\nondegeneracy at \(\bar x\)}}
  ([yshift=-3mm]exclude.east);

\node[neco,minimum height=14mm] (certify) at (.18,-95)
  {\(\bar x\) is a \textbf{strict local minimizer}};
\node[rnlp,minimum height=14mm] (sosc) at (.82,-95)
  {Second-order sufficient\\condition \textbf{holds}
   at \((\bar x,\bar z,\lambda)\)};
\draw[implication] (sosc.west)--(certify.east);

\end{tikzpicture}
\endgroup
\overviewfinish
\caption{Overview of the main relationships between NECO and RNLP.}
\label{fig:results-overview}
\end{figure}
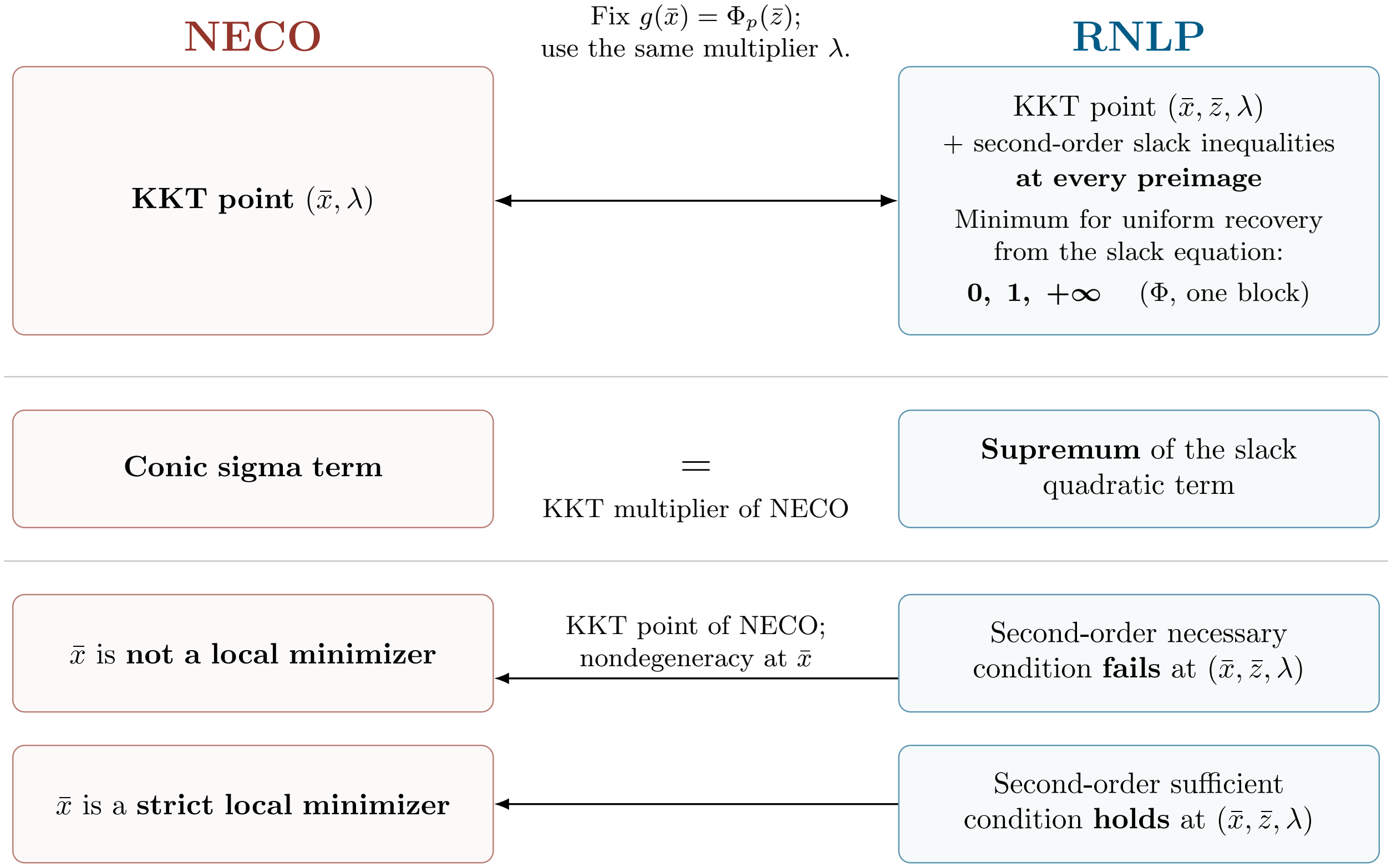

Our main contributions are the following.

\begin{enumerate}
\item
\emph{KKT transfer and second-order multiplier recovery.}
Every KKT point of NECO gives a KKT point of RNLP
at every preimage with the same multiplier.
We characterize when a KKT point of RNLP yields a
KKT point of NECO with the same multiplier through
explicit blockwise conditions, called normal compatibility.
These multiplier restrictions hold if and only if
the second-order inequalities involving only the slack
variables hold at every preimage of the same cone point.
The equivalence follows from a characterization of
the tangent cone through derivatives of \(\Phi\).

For the proposed mapping and a single cone block,
we further determine how many selected preimages
are needed to recover these restrictions uniformly over
multipliers satisfying the KKT equation for the slack
variables in RNLP.
The selection may depend on the cone point, but not
on the multiplier.
No test is needed in the interior.
On the curved boundary and the relative interior
of the flat face, one preimage suffices, as in
scalar squared-slack and unrestricted square
matrix-factor formulations
\cite[Theorem~3.3]{DingWright2025},
\cite[Theorem~3.1]{DingWright2025NSDP}.
On the remaining boundary, however, no finite selection
suffices.
The minimum number of preimages for uniform recovery
is therefore \(0\), \(1\), or \(+\infty\).

\item
\emph{Exact curvature consistency beyond symmetric cones.}
We show that, for any KKT multiplier of NECO,
the RNLP reproduces exactly the sigma term
in the second-order conditions for NECO at every preimage
on directions represented by the reformulation.
The squared-slack NSOCP and NSDP analyses
develop their curvature comparisons within proofs
assuming the second-order sufficient condition for
the reformulation, which implies strict complementarity
\cite[Proposition~3.3]{FukudaFukushima2016NSOCP},
\cite[Proposition~10]{LourencoFukudaFukushima2017NSDP}.
For the nonsymmetric exponential cone, we establish
the curvature identity without nondegeneracy, strict
complementarity, or a second-order optimality assumption.
The matrix-factor lemmas in
\cite[Lemmas~2.1--2.2]{DingWright2025NSDP}
also imply a corresponding identity for the semidefinite cone.
Our result establishes this identity for the exponential cone
without relying on the symmetric-cone square structure.

\item
\emph{Second-order optimality-condition transfer and its limits.}
We determine when second-order conditions for RNLP can
rule out a point as a local minimizer of NECO or establish
that it is a strict local minimizer of NECO.
Under nondegeneracy for NECO, the second-order necessary
condition for NECO (SONC) transfers to RNLP at every preimage
with the unique KKT multiplier of NECO, even when
the linear independence constraint qualification
(LICQ) for RNLP fails.
The resulting necessary condition uses the Lagrangian
Hessian for RNLP without separately evaluating the sigma term.
Conversely, the RNLP SONC implies the NECO SONC
when the same multiplier gives a KKT point of NECO
satisfying strict complementarity and no cone block
lies on the nonexposed ray away from the origin.

The standard second-order sufficient condition
for RNLP (SOSC) holds if and only if every cone block
lies on the curved boundary, the same multiplier gives
a KKT point of NECO satisfying strict complementarity,
and the direct fixed-multiplier second-order expression
is positive on every nonzero critical direction.
Whenever the RNLP SOSC holds, it certifies a strict
local minimizer of NECO without an additional
constraint qualification.
\end{enumerate}

We also examine how the slack representation affects
local minimizers and constraint qualifications,
providing the conditions needed to interpret candidates
obtained from RNLP in the original problem.
In the regularity analysis, we identify an intrinsic
restriction on the derivatives of \(C^2\) slack mappings
at nonzero points of the nonexposed ray.
This restriction applies to every such mapping
defined on an open set and taking values in the
exponential cone.
It shows why changing the slack mapping cannot, in general,
make nondegeneracy for NECO sufficient to ensure LICQ
for RNLP.

This paper is organized as follows.
Section~\ref{sec:preliminaries} introduces the notation,
cone geometry, and KKT and regularity conditions for NECO.
Section~\ref{sec:equality-reformulation} presents the equality
reformulation and compares local minimizers.
Section~\ref{sec:kkt} examines KKT transfer and constraint
qualifications.
Section~\ref{sec:soc} studies multiplier recovery through
second-order inequalities.
Section~\ref{sec:curvature-optimality-comparisons}
establishes curvature consistency and compares second-order
optimality conditions.
Finally, Section~\ref{sec:conclusion} discusses the results and
future research directions.


\section{Preliminaries}
\label{sec:preliminaries}
We first introduce the notation used in this paper, as well as the geometry of the exponential cone.
Throughout this paper, we use \(\inner{\cdot}{\cdot}\) and \( \norm{\cdot} \) for the Euclidean inner product and norm, respectively.
We write \(\mathbb N\coloneqq\{1,2,\ldots\}\) and
\(\mathbb R_+\coloneqq[0,+\infty)\).
For a set \(S\subseteq\mathbb R^m\),
\(\operatorname{int}S\), \(\operatorname{ri}S\), and
\(\operatorname{span}S\) denote its interior, relative interior,
and linear span, respectively. Its orthogonal complement is
\[
S^\perp\coloneqq
\{v\in\mathbb R^m\colon\langle v,s\rangle=0\ \forall s\in S\}.
\]
For a vector \(v\), we also write \(v^\perp\) for
\(\{v\}^\perp\).
For sets \(A,B\subseteq\mathbb R^m\), their sum is
\(A+B\coloneqq\{a+b\colon a\in A,\ b\in B\}\).
For a closed convex cone \(K\subseteq\mathbb R^m\),
its dual and polar are
\(K^* \coloneqq \{\xi \colon \langle\xi,y\rangle\geq0\ \forall y\in K\}\)
and \(K^\circ\mathrel{\coloneqq}-K^*\), respectively.
The tangent cone of \(K\) at \(y\in K\) is given by
\[
T_{K}(y)
\coloneqq
\{d \colon \exists\,t_k\downarrow0,\ d^k\to d,\
y+t_kd^k\in K\ \forall k\}.
\]
The normal cone of $K$ at $y$ is
\(N_{K}(y) \coloneqq T_{K}(y)^\circ\). Also, we denote
\(\operatorname{lin}T_{K}(y)
\coloneqq T_{K}(y)\cap(-T_{K}(y))\).
The tangent-cone formula for convex sets
\cite[Theorem~6.9]{RockafellarWets1998} gives
\begin{equation}
\label{eq:tangent-cone-convex-cone}
T_K(y)
=
\operatorname{cl}\bigl(
K+\operatorname{span}\{y\}
\bigr),
\end{equation}
where \(\operatorname{cl}\) denotes closure.
A convex subset \(F\subseteq K\) is a
\emph{face} if \(u,v\in K\) and \(u+v\in F\)
imply \(u,v\in F\).
It is \emph{exposed} if
\(F=K\cap\{\xi\}^{\perp}\) for some \(\xi\in K^*\).

We will also use the following standard polar-cone identity; see,
for example, \cite[Chapter~6]{RockafellarWets1998}.

\begin{lemma}
\label{lem:ri-polar-lineality}
Let \(K\subseteq\mathbb R^m\) be a closed convex cone,
and \(\xi\in\operatorname{ri}K^\circ\). Then
\( K\cap\{\xi\}^{\perp}=\operatorname{lin} K\).
\end{lemma}

The gradient and the Hessian of a function \(f \colon \mathbb{R}^n \to \mathbb{R}\) at \(x \in \mathbb{R}^n\) are denoted by \(\nabla f(x)\) and \(\nabla^2 f(x)\), respectively. For a function \(h \colon \mathbb{R}^n \times \mathbb{R}^m \to \mathbb{R}\), its partial gradient with respect to \(x\) at \((x, y)\) is denoted by \(\nabla_x h(x, y)\), and its second-order partial derivatives are denoted by \(\nabla^2_{xx} h(x, y)\).
For \(g\colon \mathbb R^n \to \mathbb R^{3p}\),
its derivative \(Dg(x)\) is identified with its Jacobian matrix, and its transpose is
denoted by \(Dg(x)^\T\).
For a matrix \(A\), \(\operatorname{Im}A\) and \(\ker A\)
denote its image and kernel, respectively.
For a twice continuously differentiable mapping
\(h=(h_1,\ldots,h_q)\colon U\to\mathbb R^q\), where
\(U\subseteq\mathbb R^m\) is open, we write
\[
D^2h(z)[w,w]\coloneqq
\bigl(w^\T\nabla^2h_j(z)w\bigr)_{j=1}^q
\]
for \(z\in U\) and \(w\in\mathbb R^m\).

For \(q\in\mathbb N\) and \(v\in\mathbb R^{qp}\), we write
\(v=(v^{(1)},\ldots,v^{(p)})\), where
\(v^{(i)}\in\mathbb R^q\), and use subscripts to denote coordinates
within each block. We write
\(g=(g^{(1)},\ldots,g^{(p)})\), where
\(g^{(i)}\colon\mathbb R^n\to\mathbb R^3\); thus,
\(g_j^{(i)}(x)\) denotes the \(j\)-th coordinate of
\(g^{(i)}(x)\). When needed, for \(z\in\mathbb R^{4p}\), we write
\(z^{(i)}=(a_i,b_i,c_i,r_i)\).

\subsection{Geometry of the exponential cone}
\label{subsec:goe-exp-cone}

Following \cite{LindstromLourencoPong2023}, with coordinates adapted
to our convention, define
\[
\begin{aligned}
F_{-\infty}
&\coloneqq
\left\{
y\in\mathbb R^3\colon
y_3=0,\ y_1\leq0,\ y_2\geq0
\right\},\\
F_\infty
&\coloneqq
\left\{
y\in\mathbb R^3\colon
y_3=0,\ y_1\leq0,\ y_2=0
\right\},\\
F_{\rm ne}
&\coloneqq
\left\{
y\in\mathbb R^3\colon
y_3=0,\ y_1=0,\ y_2\geq0
\right\},\\
F_{\rm curv}
&\coloneqq
\left\{
y\in\mathbb R^3\colon
y_3>0,\
y_2=y_3\exp\!\left(y_1 / y_3\right)
\right\}.
\end{aligned}
\]
The sets \(F_\infty\) and \(F_{\rm ne}\) are the two extreme rays of
\(F_{-\infty}\), and \(F_{\rm ne}\) is the unique nonexposed face.
The set \(F_{\rm curv}\) is the curved boundary and is not a face, as an abuse of notation.

\begin{figure}[H]
\centering
\resizebox{0.96\linewidth}{!}{%
\input{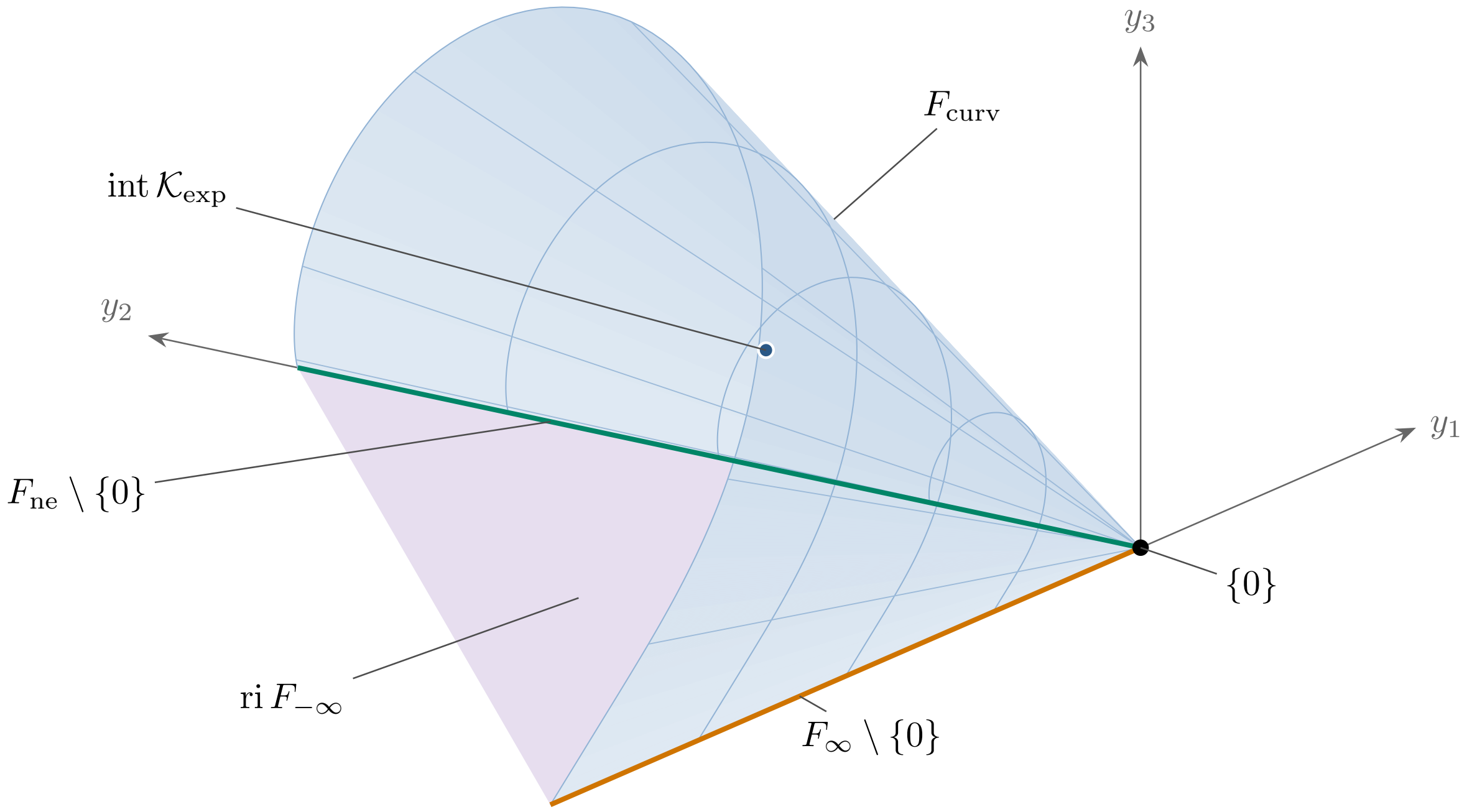}%
}
\caption{Geometry of the exponential cone: a schematic
illustration of its six disjoint subsets.
}
\label{fig:exp-cone-regions}
\end{figure}

The exponential cone then admits the disjoint decomposition
\[
\mathcal K_{\exp}
=
\operatorname{int}\mathcal K_{\exp}
\mathbin{\cup}
F_{\rm curv}
\mathbin{\cup}
\operatorname{ri}F_{-\infty}
\mathbin{\cup}
\bigl(F_{\rm ne}\setminus\{0\}\bigr)
\mathbin{\cup}
\bigl(F_\infty\setminus\{0\}\bigr)
\mathbin{\cup}
\{0\}.
\]
These six subsets distinguish the cases considered throughout
the paper. We record their coordinate descriptions for later use:
\[
\begin{aligned}
\operatorname{int}\mathcal K_{\exp}
&=
\left\{
y\in\mathbb R^3\colon
y_3>0,\
y_2>y_3\exp\!\left(y_1 / y_3\right)
\right\},\\
F_{\rm curv}
&=
\left\{
y\in\mathbb R^3\colon
y_3>0,\
y_2=y_3\exp\!\left(y_1 / y_3\right)
\right\},\\
\operatorname{ri}F_{-\infty}
&=
\left\{
y\in\mathbb R^3\colon
y_3=0,\ y_1<0,\ y_2>0
\right\},\\
F_{\rm ne}\setminus\{0\}
&=
\left\{
y\in\mathbb R^3\colon
y_3=0,\ y_1=0,\ y_2>0
\right\},\\
F_\infty\setminus\{0\}
&=
\left\{
y\in\mathbb R^3\colon
y_3=0,\ y_1<0,\ y_2=0
\right\},\\
\{0\}
&=
\left\{
y\in\mathbb R^3\colon
y_1=y_2=y_3=0
\right\}.
\end{aligned}
\]
Figure~\ref{fig:exp-cone-regions} illustrates this decomposition.

We now list the tangent cone, tangent lineality space, and normal cone
of \(\mathcal K_{\exp}\), which will be used to compare feasible directions, constraint
qualifications, and multiplier conditions in NECO and RNLP.

\begin{proposition}
\label{prop:tangent-normal-lineality-exp-cone}
Let \(y\in\mathcal K_{\exp}\), and, for \(y\in F_{\rm curv}\), set \(\nu(y)
\coloneqq
\bigl(y_2,-y_3,y_2(1-y_1/y_3)\bigr)\).
Then the tangent cone, tangent lineality space, and normal cone are
given in Table~\ref{tab:first-order-geometry}.
\end{proposition}

\begin{table}[H]
\caption{First-order geometry of the exponential cone.}
\label{tab:first-order-geometry}
\centering
\small
\sbox{\firstgeometrytablebox}{%
\begin{tabular}{@{}llll@{}}
\toprule
\(y\)
&
\(T_{\mathcal K_{\exp}}(y)\)
&
\(\operatorname{lin}T_{\mathcal K_{\exp}}(y)\)
&
\(N_{\mathcal K_{\exp}}(y)\)
\\
\midrule
\(\operatorname{int}\mathcal K_{\exp}\)
&
\(\mathbb R^3\)
&
\(\mathbb R^3\)
&
\(\{0\}\)
\\
\(F_{\rm curv}\)
&
\(\{d\colon\nu(y)^\T d\leq0\}\)
&
\(\{d\colon\nu(y)^\T d=0\}\)
&
\(\mathbb R_+\nu(y)\)
\\
\(\operatorname{ri}F_{-\infty}\) or
\(F_{\rm ne}\setminus\{0\}\)
&
\(\{d\colon d_3\geq0\}\)
&
\(\{d\colon d_3=0\}\)
&
\(\{(0,0,\lambda_3) \colon\lambda_3\leq0\}\)
\\
\(F_\infty\setminus\{0\}\)
&
\(\{d\colon d_2,d_3\geq0\}\)
&
\(\{d\colon d_2=d_3=0\}\)
&
\(\{(0,\lambda_2,\lambda_3)\colon
\lambda_2,\lambda_3\leq0\}\)
\\
\(0\)
&
\(\mathcal K_{\exp}\)
&
\(\{0\}\)
&
\(\mathcal K_{\exp}^{\circ}\)
\\
\bottomrule
\end{tabular}%
}
\global\firstgeometrytablewidth=\wd\firstgeometrytablebox
\usebox{\firstgeometrytablebox}
\end{table}

\begin{proof}
The tangent-cone formulas follow from
\cite[Theorem~5.1]{Xu2026Thesis}.
The tangent lineality spaces are obtained from its definition.
Since \(\mathcal K_{\exp}\) is closed and convex, the normal cones
are given by \(N_{\mathcal K_{\exp}}(y)
=
T_{\mathcal K_{\exp}}(y)^\circ\).
\end{proof}

For \(y=(y^{(1)},\ldots,y^{(p)})\in\mathcal K_{\exp}^p\), with
\(y^{(i)}\in\mathbb R^3\), the tangent cone, tangent lineality space, and normal cone of \(\mathcal  K_{\exp}^p\) are the Cartesian products of their one-block
counterparts.

\subsection{KKT conditions and regularity conditions}
\label{subsec:kkt-neco}

Returning to~\eqref{prob:neco}, we denote its feasible set by
\(\mathcal F_{\rm NECO}\mathrel{\coloneqq}\{x\colon g(x)\in\mathcal K_{\exp}^p\}\).
The Lagrangian function of \eqref{prob:neco} is
\[
L(x,\lambda)
\coloneqq
f(x)+\langle\lambda,g(x)\rangle,
\qquad
\lambda\in\mathbb R^{3p}.
\]
Thus, a point \(\bar x \in \mathcal F_{\rm NECO}\) is said to satisfy the KKT conditions if
there exists \(\bar\lambda\in\mathbb R^{3p}\) such that
\begin{equation} \label{eq:kkt_neco_stationary}
    \nabla f(\bar x)+Dg(\bar x)^\T\bar\lambda=0 ,
\end{equation}
\begin{equation}
    \label{eq:kkt_neco_multiplier_blockwise}
    \bar\lambda\in
N_{\mathcal K_{\exp}^p}\bigl(g(\bar x)\bigr).
\end{equation}
Note that~\eqref{eq:kkt_neco_stationary} corresponds to $\nabla_x L(\bar x,\bar\lambda) = 0$. A vector \(\bar\lambda\) satisfying
\eqref{eq:kkt_neco_stationary} and
\eqref{eq:kkt_neco_multiplier_blockwise} is called a
\emph{KKT multiplier} of \eqref{prob:neco} at \(\bar x\).
We say that a pair \((\bar x,\bar\lambda)\in
\mathcal F_{\rm NECO}\times\mathbb R^{3p}\) is a KKT point of
\eqref{prob:neco} if \(\bar\lambda\) is a KKT multiplier at
\(\bar x\). We denote the set of all KKT multipliers at a feasible
point \(\bar x\) by
\[
\Lambda(\bar x)
\coloneqq
\left\{
\lambda\in
N_{\mathcal K_{\exp}^p}\bigl(g(\bar x)\bigr)
 \colon \nabla f(\bar x)+Dg(\bar x)^\T\lambda=0
\right\}.
\]
In the following, we present some well-known definitions specialized to~\eqref{prob:neco}. For the general theory of conic optimization, we refer to~\cite{BonnansShapiro2000}.

\begin{definition}
\label{def:strict-complementarity-neco}
Let \((\bar x,\bar\lambda)\) be a KKT point of
\eqref{prob:neco}. We say that \emph{strict complementarity} for~\eqref{prob:neco}
holds at \((\bar x,\bar\lambda)\) if
\(\bar\lambda\in\operatorname{ri}
N_{\mathcal K_{\exp}^p}(g(\bar x))\).
\end{definition}
Because the normal cone is blockwise separable, strict
complementarity is equivalent to requiring each block multiplier to
lie in the relative interior of the corresponding block normal cone.
We also recall two regularity conditions for
\eqref{prob:neco}, which guarantee the boundedness and uniqueness of KKT multipliers, respectively.

\begin{definition}\label{def:rcq_neco}
For \(\bar x\in\mathcal F_{\rm NECO}\), Robinson's constraint
qualification holds for \eqref{prob:neco} at \(\bar x\) if
\begin{equation}
\label{eq:neco-robinson-cq}
\operatorname{Im}Dg(\bar x)
+
T_{\mathcal K_{\exp}^p}\bigl(g(\bar x)\bigr)
=
\mathbb R^{3p}.
\end{equation}
Equivalently,
\(0\in\operatorname{int}\bigl(g(\bar x)+Dg(\bar x)\mathbb R^n
-\mathcal K_{\exp}^p\bigr)\).
\end{definition}

\begin{definition} \label{def:nondegenaracy_neco}
    For \(\bar x\in\mathcal F_{\rm NECO}\), the problem
\eqref{prob:neco} is nondegenerate at \(\bar x\) if
\begin{equation}
\label{eq:neco-nondegeneracy}
\operatorname{Im}Dg(\bar x)
+
\operatorname{lin}T_{\mathcal K_{\exp}^p}
\bigl(g(\bar x)\bigr)
=
\mathbb R^{3p}.
\end{equation}
\end{definition}

Since
\(\operatorname{lin}T_{\mathcal K_{\exp}^p}(g(\bar x))
\subseteq T_{\mathcal K_{\exp}^p}(g(\bar x))\),
nondegeneracy~\eqref{eq:neco-nondegeneracy} implies
Robinson's constraint qualification~\eqref{eq:neco-robinson-cq}.
Nondegeneracy also ensures that \(\Lambda(\bar x)\)
is a singleton whenever it is nonempty.


\section{Exact equality representation and minimizers}
\label{sec:equality-reformulation}

This section establishes the exact equality representation and
characterizes its preimages.
We then determine when
a local minimizer of one formulation corresponds to a
local minimizer of the other, and how the slack variables affect
whether these minimizers are strict.

\subsection{Parameterization and preimages}
\label{subsec:slack-map}

Recall the mapping \(\Phi\colon\mathbb R^4\to\mathbb R^3\)
defined by
\begin{equation*}
\Phi(a,b,c,r)
=
\left(
a^2r-b^2,\,
a^2\exp(r)+c^2,\,
a^2
\right).
\end{equation*}

The following result shows that the mapping \(\Phi\) parameterizes the exponential cone exactly.

\begin{proposition}
\label{prop:image-Phi-exp-cone}
One has \(\Phi(\mathbb R^4)=\mathcal K_{\exp}\).
\end{proposition}

\begin{proof}
We first show that
\(\Phi(\mathbb R^4)\subseteq\mathcal K_{\exp}\).
Let \(y\mathrel{\coloneqq}\Phi(a,b,c,r)\). If \(a\neq0\), then \(y_3=a^2>0\) and
\[
\frac{y_1}{y_3}
=r-\frac{b^2}{a^2}
\leq r.
\]
Hence
\(y_2=a^2\exp(r)+c^2
\geq y_3\exp(y_1/y_3)\), so
\(y\in\mathcal K_{\exp}\).
If \(a=0\), then \(y=(-b^2,c^2,0)\), and again
\(y\in\mathcal K_{\exp}\).
Conversely, let \(y\in\mathcal K_{\exp}\).
If \(y_3>0\), choose
\(a=\sqrt{y_3}\), \(b=0\),
\(r=y_1/y_3\), and
\(c=\sqrt{y_2-y_3\exp(y_1/y_3)}\).
The last quantity is well-defined because
\(y\in\mathcal K_{\exp}\), and these choices give
\(\Phi(a,b,c,r)=y\).
If \(y_3=0\), then \(y_1\leq0\) and \(y_2\geq0\).
Taking \(a=0\), \(b=\sqrt{-y_1}\),
\(c=\sqrt{y_2}\), and any \(r\in\mathbb R\) again gives
\(\Phi(a,b,c,r)=y\).
\end{proof}

Note that the parameterization can be written as
\begin{equation}
\label{eq:Phi-decomposition}
\Phi(a,b,c,r)
=a^2(r,\exp(r),1)+(-b^2,c^2,0).
\end{equation}
For every \(r\in\mathbb R\), the vector
\((r,\exp(r),1)\) lies in \(F_{\rm curv}\).
Thus the first summand belongs to \(F_{\rm curv}\) when
\(a\neq0\) and is zero otherwise, while the second belongs
to \(F_{-\infty}\).

To describe the slack vectors representing a given cone point
\(y\in\mathcal K_{\exp}\), define
\[
\Phi^{-1}(y)\coloneqq
\{z\in\mathbb R^4\colon\Phi(z)=y\}.
\]
The decomposition shows, in particular, that
\(\Phi(0,b,c,r)=(-b^2,c^2,0)\) for every \(r\in\mathbb R\).
Hence varying \(r\) while keeping \(b,c\) fixed gives
infinitely many preimages of the same point.
Although \(D\Phi(0,b,c,r)\) is also independent of \(r\),
the second derivative
\(\frac{\partial^2\Phi}{\partial a^2}(0,b,c,r)
=2(r,\exp(r),1)\) depends on \(r\).
Thus different preimages can yield different second-order
inequalities for the same multiplier, a distinction examined
in Subsection~\ref{subsec:second-order-recovery-normality}
when recovering the multiplier restrictions required by
the KKT conditions for NECO.

For a general point \(y\in\mathcal K_{\exp}\), all its
preimages can be described explicitly.
Writing \(z=(a,b,c,r)\), the condition
\(z\in\Phi^{-1}(y)\) is equivalent to
\begin{equation}
\label{eq:general-preimage-relations}
a^2=y_3,\qquad
b^2=y_3r-y_1,\qquad
c^2=y_2-y_3\exp(r).
\end{equation}
These relations describe the preimages on each of the six
subsets of the exponential cone, as we can see in the following lemma.

\begin{lemma}
\label{lem:preimage-structure-exp-cone}
For every \(y\in\mathcal K_{\exp}\), the preimage
\(z = (a, b, c, r) \in \Phi^{-1}(y)\) is described in
Table~\ref{tab:preimage-structure}.
\end{lemma}

\begin{table}[htbp]
\caption{Preimages of points in the exponential cone.}
\label{tab:preimage-structure}
\centering
\small
\renewcommand{\arraystretch}{1.3}
\begin{tabular*}{\firstgeometrytablewidth}
  {@{\extracolsep{\fill}}lcccc@{}}
\toprule
Location of \(y\)
& \(a^2\)
& \(b^2\)
& \(c^2\)
& Admissible \(r\)\\
\midrule

\(\operatorname{int}\mathcal K_{\exp}\)
& \(y_3\)
& \(y_3r-y_1\)
& \(y_2-y_3\exp(r)\)
& \(\left[y_1/y_3,\,
   \log\!\left(y_2/y_3\right)\right]\)
\\

\(F_{\rm curv}\)
& \(y_3\)
& \(0\)
& \(0\)
& \(\{y_1/y_3\}\)
\\

\(\operatorname{ri}F_{-\infty}\)
& \(0\)
& \(-y_1\)
& \(y_2\)
& \(\mathbb R\)
\\

\(F_{\rm ne}\setminus\{0\}\)
& \(0\)
& \(0\)
& \(y_2\)
& \(\mathbb R\)
\\

\(F_\infty\setminus\{0\}\)
& \(0\)
& \(-y_1\)
& \(0\)
& \(\mathbb R\)
\\

\(\{0\}\)
& \(0\)
& \(0\)
& \(0\)
& \(\mathbb R\)
\\

\bottomrule
\end{tabular*}
\end{table}

\begin{proof}
By~\eqref{eq:general-preimage-relations}, it suffices to distinguish
the cases \(y_3>0\) and \(y_3=0\).

Suppose first that \(y_3>0\). Then \(a^2=y_3>0\), and the
conditions \(b^2\geq0\) and \(c^2\geq0\) are equivalent to
\[
\frac{y_1}{y_3}
\leq r\leq
\log\!\left(\frac{y_2}{y_3}\right).
\]
If \(y\in\operatorname{int}\mathcal K_{\exp}\), then
\(y_2>y_3\exp(y_1/y_3)\), so the two endpoints are distinct.
Moreover, \(b\) and \(c\) cannot both vanish, since this would give
\(r=y_1/y_3=\log(y_2/y_3)\). If \(y\in F_{\rm curv}\), the two
endpoints coincide, and hence \(r=y_1/y_3\) and \(b=c=0\).

Suppose now that \(y_3=0\). Then
\eqref{eq:general-preimage-relations} gives
\(a=0\), \(b^2=-y_1\), and \(c^2=y_2\), while
\(r\in\mathbb R\) is arbitrary. The remaining four rows of
Table~\ref{tab:preimage-structure} follow directly from the
definitions of the subsets of \(\mathcal K_{\exp}\).
\end{proof}

Note that \(\Phi^{-1}(y)\) consists of exactly two points if
\(y\in F_{\rm curv}\), whereas it is infinite for every
\(y\in\mathcal K_{\exp}\setminus F_{\rm curv}\).

For the product cone, define the block map
\(\Phi_p\colon\mathbb R^{4p}\to\mathbb R^{3p}\) as follows. For
\(z=(z^{(1)},\ldots,z^{(p)})\in\mathbb R^{4p}\), set
\(\Phi_p(z)\mathrel{\coloneqq}(\Phi(z^{(1)}),\ldots,\Phi(z^{(p)}))\).
For \(y=
\bigl(y^{(1)},\ldots,y^{(p)}\bigr)
\in\mathcal K_{\exp}^p\),
define
\[
\Phi_p^{-1}(y)
\coloneqq
\left\{
z\in\mathbb R^{4p}
 \colon 
\Phi_p(z)=y
\right\}.
\]
By block separability,
\(\Phi_p^{-1}(y)=\prod_{i=1}^p\Phi^{-1}(y^{(i)})\).

By Proposition~\ref{prop:image-Phi-exp-cone}, the reformulation
\eqref{prob:RNLP} introduced above represents the feasible set of
\eqref{prob:neco} exactly.
We denote the feasible set of \eqref{prob:RNLP} by
\[
\mathcal F_{\rm RNLP}
\coloneqq
\left\{
(x,z)\in\mathbb R^n\times\mathbb R^{4p}
 \colon 
g(x)=\Phi_p(z)
\right\}.
\]
Consequently, the feasible set of \eqref{prob:neco} is exactly the projection of
the feasible set of \eqref{prob:RNLP} onto the \(x\)-space:
\begin{equation}
\label{eq:feasible-set-projection}
\mathcal F_{\rm NECO}
=
\operatorname{proj}_x\bigl(\mathcal F_{\rm RNLP}\bigr).
\end{equation}
Indeed, by surjectivity of \(\Phi_p\),
\(g(x)\in\mathcal K_{\exp}^p\) if and only if there exists
\(z\in\mathbb R^{4p}\) such that \(g(x)=\Phi_p(z)\).

We introduce two running examples whose solutions can
be found directly, allowing the conclusions obtained from
the reformulation to be checked against the original problems.
The minimizers lie on the curved boundary in the first
example and on the nonexposed ray in the second.
We revisit these problems throughout the paper to illustrate
how RNLP can be used to analyze NECO and how the conclusions
depend on the cone geometry.

\begin{example}[Minimizers on the curved boundary]
\label{ex:running-nonconvex}
For \(n=2\), \(p=1\), and \(x=(u,v)\), consider
\begin{equation}
\tag{NECO--Ex.~\ref*{ex:running-nonconvex}}
\label{prob:running-nonconvex-neco}
\begin{aligned}
\min_{u,v\in\mathbb R}\quad
& f(x)\coloneqq v-\exp(u)+(u^2-1)^2\\
\text{subject to}\quad
& g(u,v)\coloneqq(u,v,1)\in\mathcal K_{\exp}.
\end{aligned}
\end{equation}
Feasibility is equivalent to \(v\geq\exp(u)\).
Since the objective is strictly increasing in \(v\),
every local minimizer satisfies \(v=\exp(u)\).
Thus, its cone point lies on \(F_{\rm curv}\).
The corresponding RNLP is
\begin{equation}
\tag{RNLP--Ex.~\ref*{ex:running-nonconvex}}
\label{prob:running-nonconvex-rnlp}
\begin{aligned}
\min_{u,v,a,b,c,r\in\mathbb R}\quad
& f(x)\\
\text{subject to}\quad
& u=a^2r-b^2,\\
& v=a^2\exp(r)+c^2,\\
& a^2=1.
\end{aligned}
\end{equation}
We will use RNLP to identify local minimizers of NECO
and examine how its optimality conditions distinguish
them from other candidates.
\end{example}

\begin{example}[A minimizer on the nonexposed ray]
\label{ex:running-boundary}
For \(n=p=1\), consider
\begin{equation}
\tag{NECO--Ex.~\ref*{ex:running-boundary}}
\label{prob:running-boundary-neco}
\begin{aligned}
\min_{x\in\mathbb R}\quad
& f(x)\coloneqq x^2\\
\text{subject to}\quad
& g(x)\coloneqq(0,1,x)\in\mathcal K_{\exp}.
\end{aligned}
\end{equation}
The feasible set is \([0,1]\), and the unique minimizer is
\(\bar x=0\).
The feasible cone points satisfy
\[
\begin{aligned}
g(0)&=(0,1,0)\in F_{\rm ne}\setminus\{0\},\\
g(x)&\in\operatorname{int}\mathcal K_{\exp}
\qquad(0<x<1),\\
g(1)&=(0,1,1)\in F_{\rm curv}.
\end{aligned}
\]
The corresponding RNLP is
\begin{equation}
\tag{RNLP--Ex.~\ref*{ex:running-boundary}}
\label{prob:running-boundary-rnlp}
\begin{aligned}
\min_{x,a,b,c,r\in\mathbb R}\quad
& f(x)\\
\text{subject to}\quad
& a^2r-b^2=0,\\
& a^2\exp(r)+c^2=1,\\
& x=a^2.
\end{aligned}
\end{equation}
Here, we will examine how RNLP recovers the minimizer
and why its multiplier and second-order conditions require
additional care at the nonexposed ray.
\end{example}

\subsection{Local and strict local minimizers}
\label{subsec:inner-semicontinuity}
\label{subsec:local-minimizers}

The feasible-set projection identity \eqref{eq:feasible-set-projection} guarantees that
every feasible point \(x\) of \eqref{prob:neco} admits at least one slack
\(z\in\mathbb R^{4p}\) satisfying \(g(x)=\Phi_p(z).\)
To determine whether a local minimizer of RNLP yields a
local minimizer of NECO, we also need
to compare nearby feasible points at a fixed preimage
\(\bar z\). Suppose that \(g(\bar x)=\Phi_p(\bar z).\)
If \(\{x^k\}\subseteq\mathcal F_{\rm NECO}\) satisfies \(x^k\to\bar x\), then
\eqref{eq:feasible-set-projection} guarantees that, for each \(k\), there exists at
least one \(z^k\in\mathbb R^{4p}\) such that \(g(x^k)=\Phi_p(z^k).\)
However, it does not ensure that the chosen preimages satisfy
\(z^k\to\bar z\).  
Without such a nearby choice, we cannot use the fact that
\((\bar x,\bar z)\) is a local minimizer of
\eqref{prob:RNLP} to compare its objective value with
those at the feasible points \(x^k\) approaching \(\bar x\).
This motivates the following continuity property of
\(\Phi^{-1}\).

We first state the one-block version.

\begin{definition}
\label{def:cone-level-inner-semicontinuity}
Let \(\bar y\in\mathcal K_{\exp}\), and \(\bar z\in \Phi^{-1}(\bar y)\).
We say that \(\Phi^{-1}\) is inner semicontinuous at \((\bar y,\bar z)\) relative to
\(\mathcal K_{\exp}\), if, for every sequence \(\{y^k\}\subseteq\mathcal K_{\exp}\) such that \(y^k\to\bar y\),
there exists a sequence \(\{z^k\}\subseteq\mathbb R^4\) such that
\(z^k\to\bar z\) and \(\Phi(z^k)=y^k\) for all sufficiently
large \(k\).
\end{definition}

The next result completely characterizes inner semicontinuity of
\(\Phi^{-1}\) relative to \(\mathcal K_{\exp}\) at an arbitrary preimage.

\begin{proposition}
\label{prop:cone-level-inner-semicontinuity}
Let \(\bar y\in\mathcal K_{\exp}\) and
\(\bar z\in\Phi^{-1}(\bar y)\). Then the following statements are
equivalent:
\begin{enumerate}[(i)]
\item
\(\bar y\in
\operatorname{int}\mathcal K_{\exp}
\cup F_{\rm curv}
\cup\operatorname{ri}F_{-\infty}\).
\item
The mapping \(\Phi^{-1}\) is inner semicontinuous at
\((\bar y,\bar z)\) relative to \(\mathcal K_{\exp}\).
\end{enumerate}
\end{proposition}

\begin{proof}
Write \(\bar z=(\bar a,\bar b,\bar c,\bar r)\). 

\smallskip
\noindent\emph{(i) \(\Rightarrow\) (ii).}
Let \(y^k\in\mathcal K_{\exp}\) satisfy \(y^k\to\bar y\). 
If $\bar y_3 > 0$, i.e., \(\bar y \in \operatorname{int}\mathcal K_{\exp} \cup F_{\rm curv}\), then 
\(y_3^k>0\) for all
sufficiently large \(k\). For such $k$, 
set
\(\ell_k\coloneqq y_1^k/y_3^k\),
\(u_k\coloneqq\log(y_2^k/y_3^k)\), and
\(r_k\coloneqq\min\{u_k,\max\{\ell_k,\bar r\}\}\).
Since \(\bar z\in\Phi^{-1}(\bar y)\), Lemma~\ref{lem:preimage-structure-exp-cone} gives
\(\bar y_1/\bar y_3\leq\bar r\leq
\log(\bar y_2/\bar y_3)\). Moreover,
\(\ell_k\to\bar y_1/\bar y_3\) and
\(u_k\to\log(\bar y_2/\bar y_3)\), and hence
\(r_k\to\bar r\).
Since \(\ell_k\leq r_k\leq u_k\), the quantities
\(y_3^kr_k-y_1^k\) and
\(y_2^k-y_3^k\exp(r_k)\) are nonnegative. Choose
\(a_k,b_k,c_k\) such that
\(a_k^2=y_3^k\),
\(b_k^2=y_3^kr_k-y_1^k\), and
\(c_k^2=y_2^k-y_3^k\exp(r_k)\), with their signs chosen to agree
with the corresponding nonzero components of \(\bar z\).
Then \eqref{eq:general-preimage-relations} gives
\(z^k\coloneqq(a_k,b_k,c_k,r_k)\in\Phi^{-1}(y^k)\).
Furthermore,
\(a_k^2\to\bar a^2\),
\(b_k^2\to\bar b^2\), and
\(c_k^2\to\bar c^2\), so the choice of signs, together with
\(r_k\to\bar r\), gives \(z^k\to\bar z\).

Suppose now that
\(\bar y\in\operatorname{ri}F_{-\infty}\).
Then \(\bar y_1<0\), \(\bar y_2>0\), and \(\bar y_3=0\).
Set \(r_k\coloneqq\bar r\). Since
\(y_3^k\bar r-y_1^k\to-\bar y_1>0\) and
\(y_2^k-y_3^k\exp(\bar r)\to\bar y_2>0\), for all sufficiently
large \(k\) we may choose \(a_k,b_k,c_k\) from
\eqref{eq:general-preimage-relations}, with the signs of
\(b_k,c_k\) chosen to agree with those of \(\bar b,\bar c\).
Then \(z^k\coloneqq(a_k,b_k,c_k,r_k)\in\Phi^{-1}(y^k)\) and
\(z^k\to\bar z\). Thus \(\Phi^{-1}\) is inner semicontinuous at
\((\bar y,\bar z)\) in all three subsets in~{\rm (i)}.

\smallskip
\noindent\emph{(ii) \(\Rightarrow\) (i).}
Suppose first that
\(\bar y=(0,\bar y_2,0)\), where \(\bar y_2\geq0\), which covers
\(F_{\rm ne}\setminus\{0\}\) and the origin. Then
\(\bar z=(0,0,\bar c,\bar r)\), where $\bar r$ is an arbitrary scalar. Let \(t_k\downarrow0\) and set
\(y^k\coloneqq((\bar r+1)t_k,
\bar y_2+t_k\exp(\bar r+1),t_k)\).
Then \(y^k\in\mathcal K_{\exp}\) and \(y^k\to\bar y\).
For any
\(z^k=(a_k,b_k,c_k,r_k)\in\Phi^{-1}(y^k)\),
\eqref{eq:general-preimage-relations} gives
\(r_k=\bar r+1+b_k^2/t_k\geq\bar r+1\).
Hence no such sequence \(z^k\) can converge to \(\bar z\).

Finally, suppose that
\(\bar y=(\bar y_1,0,0)\in F_\infty\setminus\{0\}\).
Let \(t_k\downarrow0\) and set
\(y^k\coloneqq(\bar y_1,t_k\exp(\bar y_1/t_k),t_k)\).
Then \(y^k\in F_{\rm curv}\) and \(y^k\to\bar y\).
By Lemma~\ref{lem:preimage-structure-exp-cone}, every
\(z^k\in\Phi^{-1}(y^k)\) has
\(r_k=\bar y_1/t_k\to-\infty\), and therefore no such sequence
can converge to \(\bar z\).
Thus \(\Phi^{-1}\) fails to be inner semicontinuous at
\((\bar y,\bar z)\) on each of the remaining three subsets.
\end{proof}

Since \(\Phi_p\) is block-separable, the one-block analysis
extends to the product cone. For local-minimizer transfer, however,
only cone points generated by feasible points of NECO are relevant.
For \(x\in\mathcal F_{\rm NECO}\), define
\[
(\Phi_p^{-1}\circ g)(x)
\coloneqq
\Phi_p^{-1}(g(x))
=
\{z\in\mathbb R^{4p}\colon \Phi_p(z)=g(x)\}.
\]

\begin{definition}
\label{def:problem-level-inner-semicontinuity}
Let \(\bar x\in\mathcal F_{\rm NECO}\), and
\(\bar z\in\Phi_p^{-1}(g(\bar x))\). We say that
\(\Phi_p^{-1}\circ g\) is inner semicontinuous at
\((\bar x,\bar z)\) relative to \(\mathcal F_{\rm NECO}\) if, for
every sequence \(\{x^k\}\subseteq\mathcal F_{\rm NECO}\) with
\(x^k\to\bar x\), there exists \(\{z^k\}\subseteq\mathbb R^{4p}\)
such that \(z^k\to\bar z\) and \(g(x^k)=\Phi_p(z^k)\) for all
sufficiently large \(k\).
\end{definition}

The following corollary gives a sufficient condition for
this problem-level continuity in terms of the cone blocks.
\begin{corollary}
\label{cor:problem-level-inner-semicontinuity}
Let \(\bar x\in\mathcal F_{\rm NECO}\) and
\(\bar z\in\Phi_p^{-1}(g(\bar x))\). If
\[
g^{(i)}(\bar x)
\in
\operatorname{int}\mathcal K_{\exp}
\cup
F_{\rm curv}
\cup
\operatorname{ri}F_{-\infty}
\qquad
\forall i=1,\ldots,p,
\]
then \(\Phi_p^{-1}\circ g\) is inner semicontinuous at
\((\bar x,\bar z)\) relative to \(\mathcal F_{\rm NECO}\).
\end{corollary}

\begin{proof}
Let \(\{x^k\}\subseteq\mathcal F_{\rm NECO}\) satisfy
\(x^k\to\bar x\). Then \(g(x^k)\in\mathcal K_{\exp}^p\) and, by
continuity of $g$, \(g(x^k)\to g(\bar x)\). For each \(i=1,\ldots,p\),
apply Proposition~\ref{prop:cone-level-inner-semicontinuity} with
\(y^k=g^{(i)}(x^k)\), \(\bar y=g^{(i)}(\bar x)\), and the prescribed
preimage \(\bar z^{(i)}\). This gives
\(z^{k,(i)}\to\bar z^{(i)}\) and
\(\Phi(z^{k,(i)})=g^{(i)}(x^k)\) for all sufficiently large \(k\).
Since there are finitely many blocks, these relations hold
simultaneously. With
\(z^k\coloneqq(z^{k,(1)},\ldots,z^{k,(p)})\), we obtain
\(z^k\to\bar z\) and \(\Phi_p(z^k)=g(x^k)\), proving the claim.
\end{proof}

\begin{remark}
\label{rem:problem-level-isc-not-necessary}
The blockwise cone-location condition in
Corollary~\ref{cor:problem-level-inner-semicontinuity}
is sufficient but not necessary for problem-level inner
semicontinuity. Indeed, problem-level inner
semicontinuity concerns only sequences of cone points generated by
feasible sequences \(x^k\to\bar x\), rather than all sequences in
\(\mathcal K_{\exp}^p\) approaching \(g(\bar x)\).
This restriction may exclude the sequences responsible for
failure of cone-level inner semicontinuity, allowing problem-level
inner semicontinuity to hold.

For example, in Example~\ref{ex:running-boundary}, take \(\bar x=0\) and
\(\bar z=(0,0,\varepsilon,\bar r)\), where
\(\varepsilon\in\{-1,1\}\) and \(\bar r\geq0\).
Here \(g(x)=(0,1,x)\) fixes the first two cone coordinates;
only the third can vary. Along these points, we can keep
\(r=\bar r\) when choosing nearby preimages. Specifically,
for any feasible sequence \(x^k\to0\), define, for all sufficiently
large \(k\),
\[
z^k\coloneqq
\left(\sqrt{x^k},\sqrt{x^k\bar r},
\varepsilon\sqrt{1-x^k\exp(\bar r)},\bar r\right).
\]
Then \(\Phi(z^k)=g(x^k)\) and \(z^k\to\bar z\), so
\(\Phi^{-1}\circ g\) is inner semicontinuous at \((0,\bar z)\).
However, \(g(0)=(0,1,0)\in F_{\rm ne}\setminus\{0\}\), and
Proposition~\ref{prop:cone-level-inner-semicontinuity} shows that
\(\Phi^{-1}\) is not inner semicontinuous at \((g(0),\bar z)\).
The sequence used in the proof of
Proposition~\ref{prop:cone-level-inner-semicontinuity}
to show failure at \(g(0)\) has positive first coordinates, whereas every point \(g(x)\)
has first coordinate zero. Thus this sequence is excluded
from the problem-level continuity requirement.
\end{remark}

The preceding continuity analysis allows us to compare
local minimizers of the two formulations.
For a local minimizer of RNLP to be strict, its selected
preimage must also be isolated, since changing only the slack variables
leaves the objective unchanged.

\begin{theorem}[Local-minimizer transfer]
\label{thm:local-minimizer-transfer}
Let \(\bar x\in\mathcal F_{\rm NECO}\) and
\(\bar z\in\Phi_p^{-1}(g(\bar x))\). Then the following statements
hold.

\begin{enumerate}[(i)]
\item
If \(\bar x\) is a local minimizer of \eqref{prob:neco}, then
\((\bar x,\bar z)\) is a local minimizer of \eqref{prob:RNLP}.
If, in addition, \(\Phi_p^{-1}\circ g\) is inner semicontinuous at
\((\bar x,\bar z)\) relative to \(\mathcal F_{\rm NECO}\), then the
converse also holds.

\item
The point \((\bar x,\bar z)\) is a strict local minimizer of
\eqref{prob:RNLP} if and only if \(\bar x\) is a strict local
minimizer of \eqref{prob:neco} and
\(g^{(i)}(\bar x)\in F_{\rm curv}\) for every \(i=1,\ldots,p\).\end{enumerate}
\end{theorem}
\begin{proof}
\noindent\emph{{\rm (i)}.}
The forward implication follows from continuity of the
projection \((x,z)\mapsto x\), as in
\cite[Proposition~3.6]{LevinKileelBoumal2025}.
The projection is surjective by \eqref{eq:feasible-set-projection},
and the inner-semicontinuity assumption makes it open at
\((\bar x,\bar z)\).
Since \(\mathcal F_{\rm RNLP}\) is a closed subset of a
finite-dimensional Euclidean space, the reverse implication
follows from \cite[Theorem~A.2]{LevinKileelBoumal2025}.

\smallskip
\noindent\emph{{\rm (ii)}.}
Suppose first that \((\bar x,\bar z)\) is a strict local minimizer of
\eqref{prob:RNLP}. Then \(\bar z\) must be isolated in
\(\Phi_p^{-1}(g(\bar x))\); otherwise, distinct feasible
pairs of RNLP with \(x=\bar x\) and objective value \(f(\bar x)\)
would approach \((\bar x,\bar z)\).
The product structure of the preimage and
Lemma~\ref{lem:preimage-structure-exp-cone} then give
\(g^{(i)}(\bar x)\in F_{\rm curv}\) for every \(i=1,\ldots,p\).

If \(\bar x\) is not a strict local minimizer of
\eqref{prob:neco}, there would exist a sequence
\(\{x^k\}\subseteq\mathcal F_{\rm NECO}\setminus\{\bar x\}\)
such that \(x^k\to\bar x\) and \(f(x^k)\leq f(\bar x)\).
Since every block of \(g(\bar x)\) lies in \(F_{\rm curv}\),
Corollary~\ref{cor:problem-level-inner-semicontinuity} gives
\(z^k\to\bar z\) with \(z^k\in\Phi_p^{-1}(g(x^k))\).
Hence \((x^k,z^k)\ne(\bar x,\bar z)\) is feasible for RNLP,
converges to \((\bar x,\bar z)\), and satisfies
\(f(x^k)\leq f(\bar x)\), which is again a contradiction.

Conversely, suppose that \(\bar x\) is a strict local minimizer of
\eqref{prob:neco} and
\(g^{(i)}(\bar x)\in F_{\rm curv}\) for every \(i\).
By Lemma~\ref{lem:preimage-structure-exp-cone}, \(\bar z\) is isolated in
\(\Phi_p^{-1}(g(\bar x))\).
Choose a neighborhood \(V\) of \(\bar z\) such that
\(V\cap\Phi_p^{-1}(g(\bar x))=\{\bar z\}\), and a neighborhood
\(U\) of \(\bar x\) in which every other feasible point of NECO
has a strictly larger objective value.
Every feasible pair of RNLP
\((x,z)\in(U\times V)\setminus\{(\bar x,\bar z)\}\)
then satisfies \(x\ne\bar x\), and hence \(f(x)>f(\bar x)\).
Therefore, \((\bar x,\bar z)\) is a strict local minimizer of
\eqref{prob:RNLP}.
\end{proof}

The feasible set of RNLP, together with the projection
\((x,z)\mapsto x\), is an instance of the construction of
Levin, Kileel, and Boumal
\cite[Definition~4.1]{LevinKileelBoumal2025}.
Their local-minimizer transfer also applies without requiring
the feasible set to be a smooth manifold.
This framework also clarifies the role of inner
semicontinuity in part~{\rm(i)}.
With the feasible sets and the selected pair
\((\bar x,\bar z)\) fixed, this assumption is necessary
and sufficient for the following property:
for every continuous objective on
\(\mathcal F_{\rm NECO}\), if \((\bar x,\bar z)\)
is a local minimizer of RNLP, then \(\bar x\)
is a local minimizer of NECO
\cite[Theorem~A.2]{LevinKileelBoumal2025}.
For a particular objective, reverse transfer may hold
without this assumption.

The reverse implication in part~{\rm (i)} can fail without
the inner-semicontinuity assumption, as the following
variation of Example~\ref{ex:running-boundary} shows.

\begin{example}
\label{ex:local-minimum-recovery-failure}
Consider Example~\ref{ex:running-boundary} with the
objective replaced by \(f(x)\mathrel{\coloneqq}-x\), keeping \(g(x)=(0,1,x)\).
Take \(\bar x=0\) and \(\bar z=(0,0,1,-1)\).
The constraint \((0,1,x)\in\mathcal K_{\exp}\) is equivalent to \(0\leq x\leq1\),
so every small positive \(x\) improves the objective
over \(\bar x\). Thus \(\bar x\) is not a local minimizer of NECO.
These improving points have no preimages near the selected
\(\bar z\).
The equality constraints for RNLP give \(x=a^2\geq0\) and \(b^2=xr\).
Hence every preimage of a positive \(x\) must have
\(r\geq0\), whereas the selected preimage has \(\bar r=-1\).
Thus improving feasible points of NECO can approach \(\bar x\)
without admitting preimages approaching \(\bar z\):
problem-level inner semicontinuity fails.

This separation also makes \((\bar x,\bar z)\) a local
minimizer of RNLP. Indeed, in a neighborhood where \(r<0\),
the same equalities force every feasible point to have
\(x=0\), so the objective remains constant.
The RNLP therefore has a local minimizer whose projection
is not a local minimizer of NECO.
\end{example}

Theorem~\ref{thm:local-minimizer-transfer}
supports the use of \eqref{prob:RNLP}
to study local minimizers of \eqref{prob:neco}.
In particular, part~{\rm(i)}
can be used to rule out a candidate as a local minimizer of NECO.
A local minimizer of NECO would give a local minimizer of RNLP at every preimage. Thus, if
\((\bar x,\bar z)\) is not a local minimizer of RNLP at even one
preimage \(\bar z\in\Phi_p^{-1}(g(\bar x))\), then
\(\bar x\) cannot be a local minimizer of NECO.
This implication does not require inner semicontinuity.

The two running examples illustrate how RNLP can both
identify local minimizers of NECO and rule out other feasible
points as local minimizers.

\begin{examplecontinuation}{ex:running-nonconvex}{Local minimizers}
\label{cont:running-nonconvex-local-minimizers}
As observed when the example was introduced, every
local minimizer of NECO satisfies \(v=\exp(u)\).
In \eqref{prob:running-nonconvex-rnlp}, substituting
the equality constraints into the objective gives
\[
f(x)
=
\exp(r)\bigl(1-\exp(-b^2)\bigr)+c^2
+\bigl((r-b^2)^2-1\bigr)^2
\geq0.
\]
For each \(s\in\{-1,1\}\), the feasible pair
\[
\bar x=(s,\exp(s)),\qquad
\bar z=(1,0,0,s)
\]
attains zero and is therefore a local minimizer
of~\eqref{prob:running-nonconvex-rnlp}.
Since \(g(\bar x)\in F_{\rm curv}\),
Corollary~\ref{cor:problem-level-inner-semicontinuity}
ensures that \(\Phi^{-1}\circ g\) is inner semicontinuous
at \((\bar x,\bar z)\) relative to \(\mathcal F_{\rm NECO}\).
The reverse implication in
Theorem~\ref{thm:local-minimizer-transfer}{\rm(i)}
then shows that \((\pm1,\exp(\pm1))\) are local minimizers
of \eqref{prob:running-nonconvex-neco}.
\end{examplecontinuation}

\begin{examplecontinuation}{ex:running-boundary}{Local minimizers}
\label{cont:running-boundary-local-minimizers}
The equality \(x=a^2\) in
\eqref{prob:running-boundary-rnlp} gives
\(f(x)=a^4\geq0\).
The feasible pair
\[
\bar x=0,\qquad \bar z=(0,0,1,0)
\]
attains zero and is therefore a local minimizer
of~\eqref{prob:running-boundary-rnlp}.
Remark~\ref{rem:problem-level-isc-not-necessary}
establishes problem-level inner semicontinuity at this pair.
The reverse implication in
Theorem~\ref{thm:local-minimizer-transfer}{\rm(i)}
therefore recovers \(\bar x=0\) as a local minimizer
of \eqref{prob:running-boundary-neco}.

At any \(\bar x\in(0,1]\), choose the preimage
\(\bar z=\bigl(\sqrt{\bar x},0,\sqrt{1-\bar x},0\bigr)\).
Decreasing \(x\) can be accompanied by a small change
in the slack variables: for \(0<t<\sqrt{\bar x}\), take
\[
x(t)\mathrel{\coloneqq}\bar x-t^2,\qquad
z(t)\mathrel{\coloneqq}\bigl(\sqrt{\bar x-t^2},0,
          \sqrt{1-\bar x+t^2},0\bigr).
\]
Then \(\Phi(z(t))=(0,1,\bar x-t^2)=g(x(t))\), and
the pairs approach \((\bar x,\bar z)\) as \(t\downarrow0\).
Since \(f(x(t))=(\bar x-t^2)^2<\bar x^2=f(\bar x)\),
the selected pair \((\bar x,\bar z)\) is not a local minimizer of RNLP.
The forward implication in
Theorem~\ref{thm:local-minimizer-transfer}{\rm(i)}
therefore rules out every \(\bar x\in(0,1]\)
as a local minimizer of NECO.
\end{examplecontinuation}

\section{First-order invariance, KKT transfer, and regularity}
\label{sec:kkt}

We next examine when the KKT conditions of
\eqref{prob:neco} can be recovered from those of
\eqref{prob:RNLP}.
The KKT equations for RNLP with respect to the slack variables
do not automatically
enforce the normal-cone condition~\eqref{eq:kkt_neco_multiplier_blockwise} of NECO.
We characterize the additional multiplier restrictions
needed for this transfer so that KKT multipliers of RNLP can be
validated for NECO. We then compare LICQ for~\eqref{prob:RNLP} with nondegeneracy for \eqref{prob:neco},
which will be used in the second-order analysis.


\subsection{The KKT system for RNLP}
\label{subsec:rnlp-kkt}

The Lagrangian function of \eqref{prob:RNLP} is
\[
\mathcal L(x,z,\mu)
\coloneqq
f(x)+\langle\mu,g(x)-\Phi_p(z)\rangle,
\qquad
\mu\in\mathbb R^{3p}.
\]
A triple \((\bar x,\bar z,\bar\mu)\) is a KKT point of
\eqref{prob:RNLP} if it satisfies
\begin{align}
\nabla f(\bar x)+Dg(\bar x)^\T\bar\mu&=0,
\label{eq:rnlp-kkt-x}\\
D\Phi_p(\bar z)^\T\bar\mu&=0,
\label{eq:rnlp-kkt-z}\\
g(\bar x)-\Phi_p(\bar z)&=0. \notag
\end{align}
Note that equation~\eqref{eq:rnlp-kkt-z} is the KKT equation for RNLP associated with the slack variables.
When \((\bar x,\bar z,\bar\mu)\) is a KKT point of
\eqref{prob:RNLP}, we call \(\bar\mu\) its
\emph{KKT multiplier} at \((\bar x,\bar z)\).

To understand the multiplier restrictions imposed by
\eqref{eq:rnlp-kkt-z}, we examine the derivative of the
slack mapping. Since \(\Phi_p\) acts separately on each
cone block, this calculation can be carried out blockwise.
For the one-block map \(\Phi\), we obtain
\begin{equation}
\label{eq:D_Phi}
D\Phi(a,b,c,r)
=
\begin{pmatrix}
2ar & -2b & 0 & a^2\\
2a\exp(r) & 0 & 2c & a^2\exp(r)\\
2a & 0 & 0 & 0
\end{pmatrix}.
\end{equation}
We next examine whether the multiplier restrictions
in \eqref{eq:rnlp-kkt-z} depend on the selected preimage,
using the image of this derivative.


\subsection{Preimage-independent first-order spaces}
\label{subsec:image-slack-derivative}

We first relate the derivative image of the parameterization to the
tangent lineality space of the exponential cone. Although
\(D\Phi(z)\) itself may depend on the selected preimage, we will show
that its image does not. This lets us describe the directions
generated by the slack mapping and the multipliers satisfying
its KKT equation in terms of the cone point alone.

\begin{lemma}
\label{lem:dphi}
Let \(y\in\mathcal K_{\exp}\) and
\(z\in\Phi^{-1}(y)\). Then one has
\[
\operatorname{Im}D\Phi(z)
\subseteq
\operatorname{lin}T_{\mathcal K_{\exp}}(y).
\]
Moreover, equality holds if and only if
\(y\notin F_{\rm ne}\setminus\{0\}\).
\end{lemma}

\begin{proof}
Write \(z=(a,b,c,r)\). We distinguish cases according to the
location of \(y\).

\smallskip
\noindent\emph{Case 1: \(y_3>0\).}
This case covers
\(y\in\operatorname{int}\mathcal K_{\exp}\cup F_{\rm curv}\).
By Lemma~\ref{lem:preimage-structure-exp-cone},
\(a\neq0\).

If \(y\in\operatorname{int}\mathcal K_{\exp}\), then \(b\) and \(c\)
cannot both vanish. If \(b\neq0\), the columns of
\(D\Phi(z)\) in~\eqref{eq:D_Phi} corresponding to \(b\), \(r\),
and \(a\) are linearly independent. If \(b=0\), then \(c\neq0\),
and the columns corresponding to \(c\), \(r\), and \(a\) are
linearly independent. Hence \(\operatorname{Im}D\Phi(z)
=
\mathbb R^3
=
\operatorname{lin}T_{\mathcal K_{\exp}}(y)\).

If \(y\in F_{\rm curv}\), then \(b=c=0\) and
\(r=y_1/y_3\). The only nonzero columns of \(D\Phi(z)\) are
\begin{equation} \label{eq:curv_D_phi}
    (2ar,2a\exp(r),2a)^\T
\quad\text{and}\quad
(a^2,a^2\exp(r),0)^\T.
\end{equation}
By Proposition~\ref{prop:tangent-normal-lineality-exp-cone},
\[
\operatorname{lin}T_{\mathcal K_{\exp}}(y)
=
\{d\in\mathbb R^3 \colon \nu(y)^\T d=0\},
\]
where
\(\nu(y)=(y_2,-y_3,y_2(1-y_1/y_3))\).
Using \(y_3=a^2\), \(y_2=a^2\exp(r)\), and
\(r=y_1/y_3\), a direct calculation shows that both columns in~\eqref{eq:curv_D_phi} are
orthogonal to \(\nu(y)\). 
Since \(\nu(y)\neq0\), 
\(\operatorname{lin}T_{\mathcal K_{\exp}}(y)=\nu(y)^\perp\)
is two-dimensional, and the two linearly independent columns in~\eqref{eq:curv_D_phi}
span 
\(\operatorname{lin}T_{\mathcal K_{\exp}}(y)\).

\smallskip
\noindent\emph{Case 2: \(y_3=0\).}
This case covers 
\[
\operatorname{ri}F_{-\infty},\qquad
F_{\rm ne}\setminus\{0\},\qquad
F_\infty\setminus\{0\},\qquad
\{0\}.
\]
Lemma~\ref{lem:preimage-structure-exp-cone} gives \(a=0\).
Hence by~\eqref{eq:D_Phi}, we have \(\operatorname{Im}D\Phi(z)
=
\operatorname{span}
\left\{
(b,0,0),\,
(0,c,0)
\right\}\).

First, if \(y\in\operatorname{ri}F_{-\infty}\), then \(b,c\neq0\), and therefore
$
\operatorname{Im}D\Phi(z)
=
\{d\colon d_3=0\}
=
\operatorname{lin}T_{\mathcal K_{\exp}}(y).
$
If \(y\in F_{\rm ne}\setminus\{0\}\), then \(b=0\) and \(c\neq0\),
and so
$
\operatorname{Im}D\Phi(z)
=
\operatorname{span}\{(0,1,0)\}
\subsetneq
\{d \colon d_3=0\}
=
\operatorname{lin}T_{\mathcal K_{\exp}}(y).
$
If \(y\in F_\infty\setminus\{0\}\), then \(b\neq0\) and \(c=0\),
which gives
$
\operatorname{Im}D\Phi(z)
=
\operatorname{span}\{(1,0,0)\}
=
\operatorname{lin}T_{\mathcal K_{\exp}}(y).
$
Finally, if \(y=0\), then \(b=c=0\), and hence
$
\operatorname{Im}D\Phi(z)
=
\{0\}
=
\operatorname{lin}T_{\mathcal K_{\exp}}(0).
$
Thus, \(\operatorname{Im}D\Phi(z)
\subseteq
\operatorname{lin}T_{\mathcal K_{\exp}}(y)\) holds
for every \(y\in\mathcal K_{\exp}\), with equality if and only if
\(y\notin F_{\rm ne}\setminus\{0\}\).
\end{proof}

The next proposition shows the preimage independence
of the derivative image.

\begin{proposition}
\label{prop:dphi-image-independence}
For \(y\in\mathcal K_{\exp}\) and every
\(z,\widehat z\in\Phi^{-1}(y)\), one has \(\operatorname{Im}D\Phi(z)
=
\operatorname{Im}D\Phi(\widehat z)\).
Likewise, for \(y\in\mathcal K_{\exp}^p\) and every
\(z,\widehat z\in\Phi_p^{-1}(y)\), one has \(\operatorname{Im}D\Phi_p(z)
=
\operatorname{Im}D\Phi_p(\widehat z)\).
\end{proposition}

\begin{proof}
The one-block assertion follows from the explicit derivative
images in the proof of Lemma~\ref{lem:dphi}.
The product assertion follows because \(D\Phi_p(z)\) is block diagonal.
\end{proof}

This independence makes the following spaces well-defined:
$
\mathcal R_\Phi(y)
\coloneqq
\operatorname{Im}D\Phi(z),
$ with 
$
z\in\Phi^{-1}(y),
$
as well as its orthogonal complement
$
\mathcal A_\Phi(y)
\coloneqq
\mathcal R_\Phi(y)^\perp
=
\ker D\Phi(z)^\T.
$
Similarly, for \(y\in\mathcal K_{\exp}^p\), define
\[
\mathcal R_{\Phi_p}(y)
\coloneqq
\operatorname{Im}D\Phi_p(z),
\qquad
\mathcal A_{\Phi_p}(y)
\coloneqq
\mathcal R_{\Phi_p}(y)^\perp,
\qquad
z\in\Phi_p^{-1}(y).
\]
Since \(\Phi_p\) acts blockwise, it follows that
\[
\mathcal R_{\Phi_p}(y)
=
\prod_{i=1}^p\mathcal R_\Phi(y^{(i)}),
\qquad
\mathcal A_{\Phi_p}(y)
=
\prod_{i=1}^p\mathcal A_\Phi(y^{(i)}).
\]
In particular, for every \(y\in\mathcal K_{\exp}^p\),
\(z\in\Phi_p^{-1}(y)\), and
\(\mu\in\mathbb R^{3p}\), one has
\begin{equation}
\label{eq:adjoint-space-equivalence}
D\Phi_p(z)^\T\mu=0
\quad\Longleftrightarrow\quad
\mu\in\mathcal A_{\Phi_p}(y).
\end{equation}
Thus, although the matrix \(D\Phi_p(z)\) depends on the selected
preimage, the solution space of~\eqref{eq:rnlp-kkt-z} depends only on
the cone point \(y\).

The relation with tangent lineality extends immediately to the
product cone.

\begin{corollary}
\label{cor:dphip-lineality}
Let \(y\in\mathcal K_{\exp}^p\). Then one has
$
\mathcal R_{\Phi_p}(y)
\subseteq
\operatorname{lin}T_{\mathcal K_{\exp}^p}(y).
$
Moreover, equality holds if and only if \(y^{(i)}\notin F_{\rm ne}\setminus\{0\}\) for all \(i=1,\ldots,p\).
\end{corollary}

\begin{proof}
Both spaces are Cartesian products of their one-block counterparts,
so the conclusion follows from Lemma~\ref{lem:dphi} blockwise.
\end{proof}

The preceding inclusion also shows that a cone normal
satisfies the KKT equation for the slack variables in RNLP.

\begin{corollary}
\label{cor:normal-implies-dphip-orthogonality}
Let \(y\in\mathcal K_{\exp}^p\). Then one has \(N_{\mathcal K_{\exp}^p}(y)
\subseteq
\mathcal A_{\Phi_p}(y)\).
\end{corollary}

\begin{proof}
Every element of \(N_{\mathcal K_{\exp}^p}(y)\) is orthogonal
to \(\operatorname{lin}T_{\mathcal K_{\exp}^p}(y)\), which contains
\(\mathcal R_{\Phi_p}(y)\) by Corollary~\ref{cor:dphip-lineality}.
Hence it belongs to
\(\mathcal R_{\Phi_p}(y)^\perp=\mathcal A_{\Phi_p}(y)\).
\end{proof}


\subsection{KKT transfer and normal compatibility}
\label{subsec:neco-to-rnlp-kkt}
\label{subsec:normal-compatible}

The two formulations~\eqref{prob:neco} and \eqref{prob:RNLP} share the KKT equation with respect to \(x\), but impose different conditions on the multiplier. Let us now establish their relations.

\begin{theorem}[{KKT transfer from NECO to RNLP}]
\label{thm:neco-kkt-implies-rnlp-kkt}
Let \((\bar x,\bar\lambda)\) be a KKT point of
\eqref{prob:neco}, and
\(\bar z\in\Phi_p^{-1}(g(\bar x))\). Then
\((\bar x,\bar z,\bar\lambda)\) is a KKT point of
\eqref{prob:RNLP}.
\end{theorem}

\begin{proof}
Since \(\bar z\in\Phi_p^{-1}(g(\bar x))\), we have \(g(\bar x)=\Phi_p(\bar z)\).
Moreover, the KKT conditions of \eqref{prob:neco} give
\[
\nabla f(\bar x)+Dg(\bar x)^\T\bar\lambda=0 \quad \text{ and } \quad \bar\lambda
\in
N_{\mathcal K_{\exp}^p}(g(\bar x)).
\]
By Corollary~\ref{cor:normal-implies-dphip-orthogonality}, we have \(\bar\lambda
\in
\mathcal A_{\Phi_p}(g(\bar x))\).
Since \(\bar z\in\Phi_p^{-1}(g(\bar x))\),
\eqref{eq:adjoint-space-equivalence} yields
$
D\Phi_p(\bar z)^\T\bar\lambda=0.
$
Thus all the KKT conditions of \eqref{prob:RNLP} hold.
\end{proof}

The next result follows from
Proposition~\ref{prop:dphi-image-independence} and shows that, if the KKT conditions for RNLP hold at one
preimage, they hold at every preimage of the same cone
point with the same \(x\) and multiplier.

\begin{corollary}
\label{cor:preimage-independent-rnlp-kkt}
If
\((\bar x,\bar z,\mu)\) is a KKT point of
\eqref{prob:RNLP}, then
\((\bar x,\widehat z,\mu)\) is a KKT point of
\eqref{prob:RNLP} for every
\(\widehat z\in\Phi_p^{-1}(g(\bar x))\).
\end{corollary}

\begin{proof}
Feasibility holds for every such \(\widehat z\), and 
\eqref{eq:rnlp-kkt-x}  does not involve the preimage. For \eqref{eq:rnlp-kkt-z},
\[
D\Phi_p(\bar z)^\T\mu=0
\quad\Longleftrightarrow\quad
\mu\in\mathcal A_{\Phi_p}(g(\bar x))
\quad\Longleftrightarrow\quad
D\Phi_p(\widehat z)^\T\mu=0,
\]
where both equivalences follow from
\eqref{eq:adjoint-space-equivalence}.
\end{proof}

To recover KKT conditions for NECO from a KKT point of RNLP,
we must check an additional restriction on the multiplier.
For a KKT point \((\bar x,\bar z,\bar\mu)\) of RNLP,
\eqref{eq:rnlp-kkt-z} is equivalent to
\(\bar\mu\in\mathcal A_{\Phi_p}(g(\bar x))\).
On the other hand, the KKT conditions for NECO require
\(\bar\mu\in N_{\mathcal K_{\exp}^p}(g(\bar x))\).
By Corollary~\ref{cor:normal-implies-dphip-orthogonality},
\[
N_{\mathcal K_{\exp}^p}(g(\bar x))
\subseteq
\mathcal A_{\Phi_p}(g(\bar x)),
\]
but this inclusion may be strict.
Thus, a KKT multiplier of RNLP need not satisfy the
normal-cone condition for NECO.

A similar first-order mismatch occurs in squared-slack
reformulations of problems with inequality
constraints~\cite{DingWright2025}
or semidefinite
constraints~\cite[Section~3.1]{DingWright2025NSDP}.
The reformulated KKT conditions need not enforce
the required sign restrictions on scalar multipliers
or positive semidefiniteness of matrix multipliers.
For the exponential-cone reformulation,
the required multiplier restrictions depend on the
location of the cone point.
The following definition specifies these restrictions.

\begin{definition}
\label{def:normal-compatible-one-block}
Let \(y\in\mathcal K_{\exp}\) and \(\mu\in\mathbb R^3\). We say that
\(\mu\) is \emph{normal-compatible} with \(y\) if the corresponding
condition in Table~\ref{tab:normal-compatibility} holds.
\end{definition}

\begin{table}[htbp]
\caption{Normal-compatibility conditions.}
\label{tab:normal-compatibility}
\centering
\small
\begin{tabular}{@{}ll@{}}
\toprule
Location of \(y\) & normal-compatibility \\
\midrule
\(\operatorname{int}\mathcal K_{\exp}\)
& \(\mu=0\)\\
\(F_{\rm curv}\)
& \(\mu_2\leq0\)\\
\(\operatorname{ri}F_{-\infty} \cup F_{\rm ne}\setminus\{0\}\)
& \(\mu_1=\mu_2=0,\ \mu_3\leq0\)\\
\(F_\infty\setminus\{0\}\)
& \(\mu_2\leq0,\ \mu_3\leq0\)\\
\(0\)
& \(\mu\in\mathcal K_{\exp}^{\circ}\)\\
\bottomrule
\end{tabular}
\end{table}

For the product cone, normal-compatibility is imposed blockwise.

\begin{definition}
\label{def:normal-compatible-product}
Let \(y\in\mathcal K_{\exp}^p\) and
\(\mu\in\mathbb R^{3p}\). We say that
\(\mu\) is \emph{blockwise normal-compatible} with \(y\) if
\(\mu^{(i)}\) is normal-compatible with \(y^{(i)}\) for every
\(i=1,\ldots,p\).
\end{definition}

Within the solution space \(\mathcal A_\Phi(y)\) of
\(D\Phi(z)^\T\mu=0\), normal compatibility is equivalent
to membership in the normal cone, as the next proposition shows.

\begin{proposition}
\label{prop:normal-characterization-one-block}
Let \(y\in\mathcal K_{\exp}\) and \(\mu\in\mathcal A_\Phi(y)\).
Then the following statements are equivalent:
\begin{enumerate}[(i)]
\item One has \(\mu\in N_{\mathcal K_{\exp}}(y)\).

\item The multiplier \(\mu\) is normal-compatible with \(y\).
\end{enumerate}
\end{proposition}

\begin{proof}
If \(y\in F_{\rm curv}\), Lemma~\ref{lem:dphi} and
Proposition~\ref{prop:tangent-normal-lineality-exp-cone} give
\(\mathcal A_\Phi(y)=\operatorname{span}\{\nu(y)\}\).
Thus \(\mu=\alpha\nu(y)\) for some \(\alpha\in\mathbb R\).
By Table~\ref{tab:normal-compatibility}, normal compatibility
is equivalent to \(\mu_2\leq0\).
Since \(\mu_2=-\alpha y_3\) and \(y_3>0\), this is
equivalent to \(\alpha\geq0\), and hence to
\(\mu\in\mathbb R_+\nu(y)=N_{\mathcal K_{\exp}}(y)\).

If \(y\in F_\infty\setminus\{0\}\), the same results give
\(\mathcal A_\Phi(y)=\{\mu\in\mathbb R^3\colon\mu_1=0\}\).
Hence \(\mu_1=0\) by assumption.
Normal compatibility requires \(\mu_2\leq0\) and
\(\mu_3\leq0\), which, together with \(\mu_1=0\),
characterize \(N_{\mathcal K_{\exp}}(y)\).

For every remaining cone subset, the normal-compatibility
conditions in Table~\ref{tab:normal-compatibility} coincide
with the normal-cone formulas in
Table~\ref{tab:first-order-geometry}.
The equivalence therefore follows directly.
\end{proof}

The characterization extends blockwise to the product cone.

\begin{corollary}
\label{cor:normal-characterization-product}
Let \(y\in\mathcal K_{\exp}^p\) and
\(\mu\in\mathcal A_{\Phi_p}(y)\). Then the following statements are
equivalent:
\begin{enumerate}[(i)]
\item One has
\(\mu\in N_{\mathcal K_{\exp}^p}(y)\).

\item The multiplier \(\mu\) is blockwise normal-compatible with \(y\).
\end{enumerate}
\end{corollary}

\begin{proof}
Applying Proposition~\ref{prop:normal-characterization-one-block} to
each block gives the equivalence.
\end{proof}

The preceding characterization gives the exact condition under which a KKT point of RNLP also satisfies the KKT conditions for NECO.

\begin{theorem}[{Characterization of KKT points of NECO}]
\label{thm:rnlp-neco-kkt-characterization}
Suppose that
\((\bar x,\bar z,\bar\mu)\) is a KKT point of
\eqref{prob:RNLP}.
Then the following statements are equivalent:
\begin{enumerate}[(i)]
\item The pair \((\bar x,\bar\mu)\) is a KKT point of
\eqref{prob:neco}.

\item The multiplier \(\bar\mu\) is blockwise normal-compatible with
\(g(\bar x)\).
\end{enumerate}
\end{theorem}

\begin{proof}
Feasibility for RNLP gives
\(\bar z\in\Phi_p^{-1}(g(\bar x))\) and
\(\bar x\in\mathcal F_{\rm NECO}\). Moreover,
\eqref{eq:rnlp-kkt-x} is precisely
\eqref{eq:kkt_neco_stationary} with
\(\bar\lambda=\bar\mu\). By
\eqref{eq:adjoint-space-equivalence}, \eqref{eq:rnlp-kkt-z} gives
\(\bar\mu\in\mathcal A_{\Phi_p}(g(\bar x))\).

The remaining KKT condition for NECO is
\(\bar\mu\in N_{\mathcal K_{\exp}^p}(g(\bar x))\),
which is equivalent to blockwise normal compatibility by
Corollary~\ref{cor:normal-characterization-product}.
\end{proof}

Theorem~\ref{thm:rnlp-neco-kkt-characterization}, together with Theorem~\ref{thm:neco-kkt-implies-rnlp-kkt}, can guide methods that identify KKT points of NECO through RNLP. Theorem~\ref{thm:neco-kkt-implies-rnlp-kkt} shows every KKT point of NECO yields a KKT point of RNLP
at every preimage with the same multiplier.
KKT points of RNLP therefore provide candidates for
KKT points of NECO.
Given a KKT point \((\bar x,\bar z,\bar\mu)\) of RNLP,
Theorem~\ref{thm:rnlp-neco-kkt-characterization}
gives an exact criterion for accepting the
candidate \((\bar x,\bar\mu)\):
\((\bar x,\bar\mu)\) is a KKT point of NECO
if and only if \(\bar\mu\) is
blockwise normal-compatible with \(g(\bar x)\).
If this condition fails, the candidate
\((\bar x,\bar\mu)\) must be rejected.

We return to the two running examples to illustrate
these results.
In the first, normal compatibility follows automatically
from the KKT equations for RNLP, and hence
every KKT point of RNLP yields a KKT point of NECO
with the same multiplier.
In the second, at \(x=0\), normal compatibility
identifies a single KKT point of NECO from the family
of KKT points of RNLP; at \(x=1\), it rejects
the corresponding candidate KKT point for NECO.

\begin{examplecontinuation}{ex:running-nonconvex}{KKT conditions}
\label{cont:running-nonconvex-kkt}
We use the KKT conditions of
\eqref{prob:running-nonconvex-rnlp}
to obtain candidate KKT points for~\eqref{prob:running-nonconvex-neco}.
The KKT equation for \(v\) in RNLP gives \(\mu_2=-1\).
If \(g(u,v)\) were an interior cone point,
Proposition~\ref{prop:tangent-normal-lineality-exp-cone}
and Lemma~\ref{lem:dphi} would give
\(\operatorname{Im}D\Phi(z)=\mathbb R^3\), so
\eqref{eq:rnlp-kkt-z} would force \(\mu=0\), a contradiction.
Since \(g_3(u,v)=1\), it follows that \(g(u,v)\in F_{\rm curv}\).

Lemma~\ref{lem:preimage-structure-exp-cone} now gives
\(a=\pm1\), \(b=c=0\), \(r=u\), and \(v=\exp(u)\).
The \(r\)-component of \eqref{eq:rnlp-kkt-z}
then gives \(\mu_1=\exp(u)\).
The KKT equation for \(u\) reduces to \(4u(u^2-1)=0\),
while the \(a\)-component of \eqref{eq:rnlp-kkt-z}
gives \(\mu_3=\exp(u)(1-u)\).
Thus, the KKT points \((\bar x,\bar z,\mu)\) of RNLP
are given by
\[
\bar x=(s,\exp(s)),
\qquad
\bar z=(\pm1,0,0,s),
\qquad
\mu=(\exp(s),-1,\exp(s)(1-s)),
\]
where \(s\in\{-1,0,1\}\).

Each candidate multiplier has \(\mu_2=-1\), so it
satisfies the normal-compatibility requirement
\(\mu_2\leq0\) on \(F_{\rm curv}\).
By Theorem~\ref{thm:rnlp-neco-kkt-characterization},
for each \(s\), the pair \((\bar x,\mu)\) is a
KKT point of \eqref{prob:running-nonconvex-neco}.
Theorem~\ref{thm:neco-kkt-implies-rnlp-kkt}
ensures that there are no other KKT points of NECO.
The values \(s=\pm1\) correspond to the minimizers identified in
\hyperref[cont:running-nonconvex-local-minimizers]
{the local-minimizer analysis of this example};
\(s=0\) gives an additional KKT point of NECO.
\end{examplecontinuation}

\begin{examplecontinuation}{ex:running-boundary}{KKT conditions}
\label{cont:running-boundary-kkt}
For \eqref{prob:running-boundary-rnlp}, the
KKT condition with respect to \(x\) is
\(2x+\mu_3=0\).
For \(0<x<1\), the cone point \(g(x)\) is interior.
As in the preceding example,
\eqref{eq:rnlp-kkt-z} forces \(\mu=0\),
contradicting \(2x+\mu_3=0\).
Thus, the KKT conditions for RNLP require
\(\bar x=0\) or \(\bar x=1\).

At \(\bar x=0\), we have \(\mu_3=0\).
By Lemma~\ref{lem:preimage-structure-exp-cone},
every preimage has the form
\(\bar z=(0,0,\varepsilon,r)\), where
\(\varepsilon\in\{-1,1\}\) and \(r\in\mathbb R\).
At each such preimage, \eqref{eq:rnlp-kkt-z}
reduces to \(2\varepsilon\mu_2=0\), giving \(\mu_2=0\)
while leaving \(\mu_1\) unrestricted.
Thus the KKT multipliers of RNLP are
$\mu=(m,0,0)$ with $m\in\mathbb R$, at every preimage.
Normal compatibility accepts exactly the candidate
with \(m=0\).

At \(\bar x=1\), the cone point \(g(1)=(0,1,1)\)
lies on \(F_{\rm curv}\).
Lemma~\ref{lem:preimage-structure-exp-cone} gives
\(\bar z=(\pm1,0,0,0)\).
The KKT equation for \(x\) gives \(\mu_3=-2\), while
the \(r\)- and \(a\)-components of \eqref{eq:rnlp-kkt-z}
give \(\mu_1+\mu_2=0\) and \(\mu_2+\mu_3=0\).
Hence \(\mu=(-2,2,-2)\).
Since \(\mu_2=2>0\), this multiplier
fails the normal-compatibility requirement on
\(F_{\rm curv}\).
Thus, the candidate
\((\bar x,\mu)=(1,(-2,2,-2))\)
is not a KKT point of NECO.

By Theorems~\ref{thm:neco-kkt-implies-rnlp-kkt}
and~\ref{thm:rnlp-neco-kkt-characterization},
the only KKT point of
\eqref{prob:running-boundary-neco} is
\((\bar x,\mu)=(0,0)\).
\end{examplecontinuation}


\subsection{Constraint qualifications}
\label{sec:cq}

We now compare the linear independence constraint qualification
(LICQ) for \eqref{prob:RNLP} with nondegeneracy for
\eqref{prob:neco}. This determines when a rank condition
for the equality reformulation can establish regularity of the
original problem.
The comparison depends on the relation between the
preimage-independent derivative image
\(\mathcal R_{\Phi_p}(g(\bar x))\) and the tangent lineality space
\(\operatorname{lin}T_{\mathcal K_{\exp}^p}(g(\bar x))\).

\begin{definition}
\label{def:rnlp-licq}
For \((\bar x,\bar z)\in\mathcal F_{\rm RNLP}\), LICQ holds for
\eqref{prob:RNLP} at \((\bar x,\bar z)\) if
\begin{equation*}
\operatorname{Im}Dg(\bar x)
+
\operatorname{Im}D\Phi_p(\bar z)
=
\mathbb R^{3p}.
\end{equation*}
\end{definition}

By Proposition~\ref{prop:dphi-image-independence} and the definition of \(\mathcal R_{\Phi_p}(y)\), the LICQ condition is equivalently written as
\begin{equation}
\label{eq:rnlp-licq-preimage-independent}
\operatorname{Im}Dg(\bar x)
+
\mathcal R_{\Phi_p}(g(\bar x))
=
\mathbb R^{3p}.
\end{equation}
Therefore, LICQ for RNLP is independent of the selected
\(\bar z\in\Phi_p^{-1}(g(\bar x))\).

Corollary~\ref{cor:dphip-lineality} shows that
\(\mathcal R_{\Phi_p}(g(\bar x))
\subseteq
\operatorname{lin}T_{\mathcal K_{\exp}^p}(g(\bar x))\), with
equality whenever no block of \(g(\bar x)\) lies in
\(F_{\rm ne}\setminus\{0\}\). This yields the following comparison of LICQ and nondegeneracy.

\begin{theorem}[{LICQ for RNLP and nondegeneracy for NECO}]
\label{thm:licq-nondegeneracy}
Let \((\bar x,\bar z)\in\mathcal F_{\rm RNLP}\). Then the following
statements hold:
\begin{enumerate}[(i)]
\item
If LICQ holds for \eqref{prob:RNLP} at
\((\bar x,\bar z)\), then \eqref{prob:neco} is nondegenerate at
\(\bar x\).

\item
Suppose that
\(g^{(i)}(\bar x)\notin F_{\rm ne}\setminus\{0\}\) for every
\(i=1,\ldots,p\). Then LICQ holds for \eqref{prob:RNLP} at
\((\bar x,\bar z)\) if and only if \eqref{prob:neco} is
nondegenerate at \(\bar x\).
\end{enumerate}
\end{theorem}

\begin{proof}
Both assertions follow directly from
Corollary~\ref{cor:dphip-lineality} by comparing \eqref{eq:neco-nondegeneracy} and
\eqref{eq:rnlp-licq-preimage-independent}.
\end{proof}

Theorem~\ref{thm:licq-nondegeneracy} makes LICQ for RNLP a
sufficient test for nondegeneracy for NECO.
We first check LICQ for RNLP at the points identified in
the KKT analysis. When it holds, part~{\rm(i)} gives
nondegeneracy for NECO. When it fails, we examine
nondegeneracy for NECO separately.

\begin{examplecontinuation}{ex:running-nonconvex}{Constraint qualifications}
\label{cont:running-nonconvex-cq}
For each \(s\in\{-1,0,1\}\), consider the point
\(\bar x=(s,\exp(s))\).
Then
\(\operatorname{Im}Dg(\bar x)=\{v\in\mathbb R^3\colon v_3=0\}\).
The
\hyperref[cont:running-nonconvex-kkt]
{KKT analysis of Example~\ref*{ex:running-nonconvex}}
gives the preimages
\(\bar z=(a,0,0,s)\), where \(a=\pm1\).
By \eqref{eq:D_Phi}, the first column of \(D\Phi(\bar z)\)
is \(2a(s,\exp(s),1)^\T\).
Its third coordinate is nonzero, and so
$
\operatorname{Im}Dg(\bar x)
+\operatorname{Im}D\Phi(\bar z)
=\mathbb R^3.
$
Thus LICQ holds for \eqref{prob:running-nonconvex-rnlp}
at every corresponding preimage.
Theorem~\ref{thm:licq-nondegeneracy}{\rm(i)}
then shows that \eqref{prob:running-nonconvex-neco}
is nondegenerate at each of the three points.
\end{examplecontinuation}

\begin{examplecontinuation}{ex:running-boundary}{Constraint qualifications}
\label{cont:running-boundary-cq}
At \(\bar x=1\), the corresponding preimages are
\(\bar z=(a,0,0,0)\), where \(a=\pm1\).
The column of \(Dg(1)\) and the first and fourth columns
of \(D\Phi(\bar z)\) are, respectively,
$(0,0,1)^\T$, $(0,2a,2a)^\T$, and $(1,1,0)^\T$.
These vectors are linearly independent because \(a\neq0\).
Hence LICQ holds for \eqref{prob:running-boundary-rnlp}
at both preimages.
Theorem~\ref{thm:licq-nondegeneracy}{\rm(i)}
therefore gives nondegeneracy
of~\eqref{prob:running-boundary-neco} at \(\bar x=1\).

At \(\bar x=0\), we again begin with LICQ for RNLP.
For the preimages described in
\hyperref[cont:running-boundary-kkt]
{the KKT analysis of Example~\ref*{ex:running-boundary}},
\eqref{eq:D_Phi} gives
\(\operatorname{Im}D\Phi(\bar z)
=\operatorname{span}\{(0,1,0)\}\).
Together with
\(\operatorname{Im}Dg(\bar x)
=\operatorname{span}\{(0,0,1)\}\),
this yields
\[
\operatorname{Im}Dg(\bar x)
+\operatorname{Im}D\Phi(\bar z)
=\operatorname{span}\{(0,0,1),(0,1,0)\}
\neq\mathbb R^3.
\]
Thus LICQ for \eqref{prob:running-boundary-rnlp} fails at
\((\bar x,\bar z)\) for every preimage \(\bar z\).

To check nondegeneracy for NECO separately, note that
\(g(\bar x)=(0,1,0)\in F_{\rm ne}\setminus\{0\}\).
Proposition~\ref{prop:tangent-normal-lineality-exp-cone}
gives
\(\operatorname{lin}T_{\mathcal K_{\exp}}(g(\bar x))
=\{v\in\mathbb R^3\colon v_3=0\}\).
Consequently,
$
\operatorname{Im}Dg(\bar x)
+\operatorname{lin}T_{\mathcal K_{\exp}}(g(\bar x))
=\mathbb R^3,
$
and thus \eqref{prob:running-boundary-neco}
is nondegenerate at \(\bar x\).
\end{examplecontinuation}

The analysis of constraint qualifications in
\hyperref[cont:running-boundary-cq]
{Example~\ref*{ex:running-boundary}}
shows that nondegeneracy for NECO need not imply LICQ for RNLP
and that this failure persists at every preimage.
In this example, the failure arises from the strict
inclusion in Lemma~\ref{lem:dphi} for the proposed
mapping \(\Phi\).
One may therefore ask whether another \(C^2\)
slack mapping could restore LICQ for RNLP.
The following proposition shows that
at every nonzero point of the nonexposed ray \(F_{\rm ne}\),
\(\Phi\) already attains the largest derivative image
permitted by the cone geometry.
Consequently, replacing \(\Phi\) with another such
\(C^2\) slack mapping cannot make nondegeneracy for NECO
sufficient for LICQ for RNLP in general.

\begin{proposition}[An intrinsic restriction on smooth slack mappings]
\label{prop:smooth-map-nonexposed-ray}
Let \(U\subseteq\mathbb R^m\) be open
and \(\Psi \colon U\to\mathbb R^3\) be twice continuously
differentiable with \(\Psi(U)\subseteq\mathcal K_{\exp}\).
If \(\bar z\in U\) satisfies \(\bar y\coloneqq\Psi(\bar z)=(0,\beta,0) \in F_{\rm ne}\) where \(\beta>0\),
then
\[
\operatorname{Im}D\Psi(\bar z)
\subseteq
\operatorname{span}F_{\rm ne}
\subsetneq
\operatorname{lin}T_{\mathcal K_{\exp}}(\bar y).
\]

\end{proposition}

\begin{proof}
Fix \(h\in\mathbb R^m\) and define
\(\gamma(t)\coloneqq\Psi(\bar z+th)\)
for \(t\in \R\) such that \(\bar z+th \in U\).
Since \(\Psi(U)\subseteq\mathcal K_{\exp}\),
we have \(\gamma(t)\in\mathcal K_{\exp}\) and hence
\(\gamma_3(t)\geq0\).
Together with \(\gamma_3(0)=0\), this shows that
\(t=0\) is a local minimizer of \(\gamma_3\).
Therefore, \(\gamma_3'(0)=0\).
Since \(\gamma_3\) is twice continuously differentiable
and \(\gamma_3(0)=\gamma_3'(0)=0\), Taylor expansion at \(0\) gives
\(\gamma_3(t)=O(t^2)\).
Moreover, continuity of \(\gamma_2\) and
\(\gamma_2(0)=\beta>0\) ensure that \(\gamma_2\) remains
positive and bounded near \(0\).
Thus, there exist constants
\(C,M,\delta>0\) such that
\[
0\leq\gamma_3(t)\leq Ct^2,
\qquad
0<\gamma_2(t)\leq M
\qquad (|t|<\delta).
\]
Shrink \(\delta\), if necessary, so that
\(C\delta^2<M/e\).
If \(\gamma_3(t)>0\), then \(\gamma(t) \in \mathcal K_{\exp}\) gives
\[
\gamma_1(t)
\leq
\gamma_3(t)
\log\!\left(\frac{\gamma_2(t)}{\gamma_3(t)}\right)
\leq
\gamma_3(t)
\log\!\left(\frac{M}{\gamma_3(t)}\right).
\]
Since \(s\mapsto s\log(M/s)\) is increasing on
\((0,M/e)\), it follows that
\[
\gamma_1(t)
\leq
Ct^2\log\!\left(\frac{M}{Ct^2}\right)
=o(|t|)
\qquad (t\to0).
\]
If \(\gamma_3(t)=0\), then \(\gamma(t) \in \mathcal K_{\exp}\) gives \(\gamma_1(t)\leq0\).
For \(0<|t|<\delta\), we also have \(Ct^2<M/e\), so
\[
\gamma_1(t)\leq0
<
Ct^2\log\!\left(\frac{M}{Ct^2}\right).
\]
Hence, the same upper bound holds in this case.
Using \(\gamma_1(0)=0\) and considering both signs of \(t\)
yields \(\gamma_1'(0)=0\).
Consequently,
$
D\Psi(\bar z)h
=\gamma'(0)
=(0,\gamma_2'(0),0).
$
Since \(h\) was arbitrary and
\(\operatorname{span}F_{\rm ne}
=\operatorname{span}\{(0,1,0)\}\),
this proves the first inclusion.
The strict inclusion follows immediately from
Proposition~\ref{prop:tangent-normal-lineality-exp-cone}.
\end{proof}

By Lemma~\ref{lem:dphi}, the proposed mapping \(\Phi\)
attains the largest derivative image permitted by this
restriction. The failure of LICQ for RNLP in the example
therefore reflects a restriction on every such smooth slack
representation. Applying the proposition to
\hyperref[cont:running-boundary-cq]
{Example~\ref*{ex:running-boundary}},
we obtain, for any such \(\Psi\) and any
\(\bar z\) with \(\Psi(\bar z)=g(0)\),
\[
\operatorname{Im}Dg(0)+\operatorname{Im}D\Psi(\bar z)
\subseteq
\operatorname{span}\{(0,1,0),(0,0,1)\}
\subsetneq\mathbb R^3.
\]
Thus, changing the \(C^2\) slack mapping cannot restore
LICQ for RNLP in this example.


\section{Second-order multiplier recovery}
\label{sec:soc}

The preceding first-order optimality comparison shows that the KKT
equation for RNLP~\eqref{eq:rnlp-kkt-z} need not imply the normal-cone
condition~\eqref{eq:kkt_neco_multiplier_blockwise}. For scalar squared-slack and square
matrix-factor reformulations, Ding and Wright show that
second-order necessary conditions recover the missing
multiplier restrictions
\cite{DingWright2025, DingWright2025NSDP}. We now ask when
second-order information from RNLP recovers a KKT multiplier of NECO.
After introducing the second-order necessary conditions for NECO and RNLP, we give an equivalent characterization of the missing
normal-cone condition and determine the minimum number of
preimages needed for its recovery.

\subsection{Second-order necessary conditions}
\label{subsec:direct-sigma-soc}
\label{subsec:rnlp-soc}

We first recall the second-order tangent set and the associated sigma term used
in the direct second-order analysis of \eqref{prob:neco}. Let
\(K\subseteq\mathbb R^m\) be closed, \(y\in K\), and
\(v\in T_K(y)\). The second-order tangent set to \(K\) at \(y\) in
the direction \(v\) is defined by
\[
T_K^2(y,v)
\coloneqq
\left\{
\eta\in\mathbb R^m
 \colon 
\begin{array}{l}
\exists\,t_k\downarrow0,\ \exists\,\eta^k\to\eta
\text{ such that}\\
y+t_kv+\frac12t_k^2\eta^k\in K
\quad \forall k
\end{array}
\right\}.
\]
For \(\lambda\in\mathbb R^m\) and 
\(T\subseteq\mathbb R^m\), the support function of \(T\) at
\(\lambda\) is defined in the extended-real sense by
\begin{equation*}
\sigma(\lambda,T)
\coloneqq
\sup_{\eta\in T}\langle\lambda,\eta\rangle.
\end{equation*}
We use the convention
\(\sigma(\lambda,\emptyset)=-\infty\).
For later use, we show a consequence of the
definition of the second-order tangent set.

\begin{lemma}
\label{lem:parameterization-second-order-tangent}
Let \(U\subseteq\mathbb R^m\) be open, 
\(K\subseteq\mathbb R^q\) be closed, and 
\(\Psi\colon U\to\mathbb R^q\) be twice continuously differentiable
with \(\Psi(U)\subseteq K\). Then, for every \(z\in U\) and
every \(w\in\mathbb R^m\), one has
\begin{equation}
\label{eq:parameterization-second-order-tangent}
D^2\Psi(z)[w,w]
\in T_K^2\bigl(\Psi(z),D\Psi(z)w\bigr).
\end{equation}
In particular, \(T_K^2(\Psi(z),v)\neq\emptyset\) for every
\(v\in\operatorname{Im}D\Psi(z)\).
\end{lemma}

\begin{proof}
Since \(U\) is open, \(0\in T_U^2(z,w)\).
The inclusion follows from the forward inclusion established
in the proof of
\cite[Proposition~13.13]{RockafellarWets1998}.
For any \(v\in\operatorname{Im}D\Psi(z)\), choosing \(w\)
such that \(D\Psi(z)w=v\) proves the nonemptiness assertion.
\end{proof}

For a feasible point \(\bar x\) of \eqref{prob:neco},
the critical cone is defined by
\[
C(\bar x)
\coloneqq
\left\{
d\in\mathbb R^n\colon
Dg(\bar x)d
\in T_{\mathcal K_{\exp}^p}\bigl(g(\bar x)\bigr),
\quad
\langle\nabla f(\bar x),d\rangle=0
\right\}.
\]
For \(\bar x\in\mathcal F_{\rm NECO}\),
\(\lambda\in\mathbb R^{3p}\), and \(d\in\mathbb R^n\) with
\(Dg(\bar x)d\in T_{\mathcal K_{\exp}^p}(g(\bar x))\), define
\begin{equation}
\label{eq:direct-sigma-expression}
\Theta_{\bar x,\lambda}(d)
\coloneqq
d^\T\nabla_{xx}^2L(\bar x,\lambda)d
-
\sigma\!\left(
\lambda,
T_{\mathcal K_{\exp}^p}^2
\bigl(g(\bar x),Dg(\bar x)d\bigr)
\right).
\end{equation}
The following definition states the second-order
necessary condition for NECO.

\begin{definition}
\label{def:direct-sonc-neco}
Suppose that \(\bar x\) satisfies the KKT conditions for
\eqref{prob:neco}. We say that the
\emph{direct sigma-term second-order necessary condition for
\eqref{prob:neco}}
(\emph{NECO SONC}) holds at \(\bar x\) if
\[
\sup_{\lambda\in\Lambda(\bar x)}
\left\{
d^\T\nabla_{xx}^2L(\bar x,\lambda)d
-
\sigma\bigl(\lambda,\mathcal T(d)\bigr)
\right\}
\geq0
\]
for every \(d\in C(\bar x)\) and every nonempty convex set
\[
\mathcal T(d)
\subseteq
T_{\mathcal K_{\exp}^p}^2
\bigl(g(\bar x),Dg(\bar x)d\bigr).
\]
\end{definition}

By the general sigma-term second-order necessary condition
\cite[Theorem~3.45]{BonnansShapiro2000}, applied
directly to the closed convex set
\(\mathcal K_{\exp}^p\), every local minimizer of
\eqref{prob:neco} at which Robinson's constraint qualification
holds satisfies the NECO SONC.
See also \cite[Theorem~5.7]{Xu2026Thesis} for the
single exponential-cone case.

Definition~\ref{def:direct-sonc-neco} involves the supremum
over the entire multiplier set \(\Lambda(\bar x)\).
The direction \(d\) and the set \(\mathcal T(d)\)
are fixed before taking this supremum.
The definition does not require a single multiplier
to satisfy the inequality for every such choice.
In the subsequent comparison with \eqref{prob:RNLP},
we fix a multiplier \(\lambda\in\Lambda(\bar x)\)
and examine \(\Theta_{\bar x,\lambda}(d)\)
as the direction \(d\) varies.

We next formulate the second-order necessary condition for RNLP.
For a feasible point \((\bar x,\bar z)\) of \eqref{prob:RNLP}, define the
linearized space
\[
\mathcal S_{\rm R}(\bar x,\bar z)
\coloneqq
\left\{
(d,w)\in\mathbb R^n\times\mathbb R^{4p}
 \colon 
Dg(\bar x)d-D\Phi_p(\bar z)w=0
\right\}.
\]
At a KKT point \((\bar x,\bar z,\bar\mu)\) of
\eqref{prob:RNLP},
the KKT conditions~\eqref{eq:rnlp-kkt-x}
and~\eqref{eq:rnlp-kkt-z}
give
\begin{equation} \label{eq:rnlp-critical-orthogonality}
    \langle\nabla f(\bar x),d\rangle = -\langle\bar\mu,Dg(\bar x)d\rangle =  -\langle\bar\mu,D\Phi_p(\bar z)w\rangle 
    = 0 \quad \forall (d,w)\in\mathcal S_{\rm R}(\bar x,\bar z).
\end{equation} 
For \(\bar\mu\in\mathbb R^{3p}\), the Hessian of the Lagrangian for RNLP gives the quadratic form
\begin{equation}
\label{eq:Q-RNLP}
Q_{\bar x,\bar z,\bar\mu}(d,w)
\coloneqq
d^\T\nabla_{xx}^2L(\bar x,\bar\mu)d
-
\sum_{i=1}^p
\left\langle
\bar\mu^{(i)},
D^2\Phi(\bar z^{(i)})[w^{(i)},w^{(i)}]
\right\rangle.
\end{equation}

Using~\eqref{eq:rnlp-critical-orthogonality}, we formulate
the RNLP SONC on \(\mathcal S_{\rm R}(\bar x,\bar z)\) below.

\begin{definition}
\label{def:rnlp-sonc}
Let \((\bar x,\bar z,\bar\mu)\) be a KKT point of
\eqref{prob:RNLP}. We say that the
\emph{second-order necessary condition for \eqref{prob:RNLP}}
(\emph{RNLP SONC}) holds at
\((\bar x,\bar z,\bar\mu)\) if
\[
Q_{\bar x,\bar z,\bar\mu}(d,w)\geq0
\qquad
 \forall 
(d,w)\in\mathcal S_{\rm R}(\bar x,\bar z).
\]
\end{definition}

The RNLP SONC is the standard second-order necessary
condition for equality-constrained nonlinear programming.
If \((\bar x,\bar z)\) is a local minimizer of \eqref{prob:RNLP}
and LICQ holds at \((\bar x,\bar z)\), then the RNLP SONC holds for
its KKT multiplier; see, for example,
\cite{bertsekas2016,NocedalWright2006}.

\subsection{Second-order recovery of KKT conditions for NECO}
\label{subsec:second-order-recovery-normality}

Subsection~\ref{subsec:normal-compatible} identified the additional
multiplier restrictions needed to recover the
KKT conditions for NECO
from the KKT conditions for RNLP.
We now characterize these restrictions through
second-order inequalities involving only slack directions.
A related mechanism occurs for squared slacks: second-order necessary
conditions recover multiplier signs
not enforced by the reformulated KKT conditions%
~\cite[Theorem~3.3]{DingWright2025}.
Under nondegeneracy assumptions, it was shown that inequality
multipliers with the wrong sign produce directions of negative
curvature for the reformulated Lagrangian~\cite[Theorem~3.8]{JongenStein2003}. This explains why
second-order information can detect multiplier restrictions
that the first-order equations do not enforce.

Suppose that the RNLP SONC holds at
\((\bar x,\bar z,\bar\mu)\).
By the definition of \(\mathcal S_{\rm R}(\bar x,\bar z)\),
one has
\((0,w)\in\mathcal S_{\rm R}(\bar x,\bar z)\)
if and only if \(w\in\ker D\Phi_p(\bar z)\).
Hence,
\[
Q_{\bar x,\bar z,\bar\mu}(0,w)
=
-\sum_{i=1}^p
\left\langle
\bar\mu^{(i)},
D^2\Phi\bigl(\bar z^{(i)}\bigr)
\bigl[w^{(i)},w^{(i)}\bigr]
\right\rangle
\geq0
\]
for every \(w\in\ker D\Phi_p(\bar z)\).

The RNLP SONC gives this inequality only at the selected
preimage \(\bar z\).
As explained after \eqref{eq:Phi-decomposition}, for a single
block, varying \(r\) in \(z=(0,b,c,r)\), with \(b,c\) fixed,
leaves \(\Phi(z)\) and \(D\Phi(z)\) unchanged but can change
the quadratic expression in this inequality.
This motivates the condition introduced below, which requires
the inequality to hold
at every preimage of \(g(\bar x)\), with the same multiplier
\(\bar\mu\).
We first formulate the condition for one exponential-cone block.

\begin{definition}
\label{def:all-preimage-slack-condition-one-block}
Let \(y\in\mathcal K_{\exp}\) and
\(\mu\in\mathbb R^3\). We say that the
\emph{all-preimage slack-only second-order condition} (\emph{AP-SOC}) holds at
\((y,\mu)\) if
\begin{equation}
\label{eq:ap-soc-one-block}
-\left\langle
\mu,D^2\Phi(z)[w,w]
\right\rangle
\geq0
\qquad
\forall z\in\Phi^{-1}(y),\quad
\forall w\in\ker D\Phi(z).
\end{equation}
\end{definition}

To evaluate the AP-SOC for~\eqref{prob:RNLP},
we use the explicit second derivative of \(\Phi\).
For \(z=(a,b,c,r)\) and \(w=(w_a,w_b,w_c,w_r)\),
direct differentiation gives
\begin{equation}
\label{eq:D2Phi-full}
\begin{aligned}
D^2\Phi(z)[w,w]
=
\bigl(&
2rw_a^2+4aw_aw_r-2w_b^2,\\
&
2\exp(r)w_a^2+4a\exp(r)w_aw_r
+a^2\exp(r)w_r^2+2w_c^2,
2w_a^2
\bigr).
\end{aligned}
\end{equation}

When \(y_3=0\), the second derivative of \(\Phi\) has the
following additional property, which will be used in tangent-cone
recovery and the comparison of second-order conditions.

\begin{lemma}
\label{lem:second-derivative-flat-preimages}
If \(y \in \mathcal K_{\exp} \setminus \left(\operatorname{int}\mathcal K_{\exp} \cup F_{\rm curv} \right) \), then for every
\(z=(a,b,c,r)\in\Phi^{-1}(y)\) and every
\(w=(w_a,w_b,w_c,w_r)\in\mathbb R^4\), one has
\begin{equation}\label{eq:D2Phi-flat-preimage}
    D^2\Phi(z)[w,w] 
=\Phi(\sqrt{2}w_a,\sqrt{2}w_b,\sqrt{2}w_c,r)
\in\mathcal K_{\exp}.
\end{equation}
\end{lemma}

\begin{proof}
The assumption on \(y\) gives \(y_3=0\), so
\eqref{eq:general-preimage-relations} yields \(a=0\).
Substituting \(a=0\) into \eqref{eq:D2Phi-full} 
gives 
\[
D^2\Phi(z)[w,w]
=2w_a^2(r,\exp(r),1)+(-2w_b^2,2w_c^2,0).
\]
It then follows from \eqref{eq:Phi-decomposition} that
\eqref{eq:D2Phi-flat-preimage} holds.
\end{proof}

Recall from Proposition~\ref{prop:normal-characterization-one-block}
that, for \(\mu\in\mathcal A_\Phi(y)\),
normal compatibility with \(y\) is equivalent to
\(\mu\in N_{\mathcal K_{\exp}}(y)\).
To show that the AP-SOC characterizes normal compatibility,
we use the following lemma.
Following the derivative representation at a fixed
preimage in
\cite[Proposition~3.26(a)]{LevinKileelBoumal2025},
we consider the corresponding sets over all preimages
and take the closure of their union.
For the matrix parametrization \(R\mapsto RR^\T\),
\cite[Proposition~3.32]{LevinKileelBoumal2025}
already establishes a tangent-cone representation
at each fixed preimage.
The following lemma establishes the representation
for \(\Phi\) by considering all preimages and their limits.

\begin{lemma}[Tangent-cone recovery from all preimages]
\label{lem:tangent-cone-all-preimages}
For every \(y\in\mathcal K_{\exp}\), one has
\begin{equation}
\label{eq:tangent-cone-all-preimages}
N_{\mathcal K_{\exp}}(y)^\circ
=
\operatorname{cl}\Bigl(
\mathcal A_\Phi(y)^\perp
+\bigcup_{z\in\Phi^{-1}(y)}
\left\{D^2\Phi(z)[w,w] \colon w\in\ker D\Phi(z)\right\}
\Bigr).
\end{equation}
\end{lemma}

\begin{proof}
Since \(N_{\mathcal K_{\exp}}(y)^\circ
=T_{\mathcal K_{\exp}}(y)\) and
\(\mathcal A_\Phi(y)^\perp=\mathcal R_\Phi(y)\),
we prove the two inclusions in
\eqref{eq:tangent-cone-all-preimages}.

We first show that the right-hand side is contained in
\(T_{\mathcal K_{\exp}}(y)\).
For \(z\in\Phi^{-1}(y)\) and \(w\in\ker D\Phi(z)\),
Lemma~\ref{lem:parameterization-second-order-tangent},
applied to \(\Phi\), gives
\[
D^2\Phi(z)[w,w]
\in T^2_{\mathcal K_{\exp}}(y,0)
=T_{\mathcal K_{\exp}}(y).
\]
Moreover, Lemma~\ref{lem:dphi} gives
\(\mathcal R_\Phi(y)\subseteq
\operatorname{lin}T_{\mathcal K_{\exp}}(y)\).
Because the tangent cone is a closed convex cone,
adding vectors in \(\mathcal R_\Phi(y)\) and taking
the closure preserve this inclusion.
This proves the first inclusion.

For the reverse inclusion, we distinguish three cases.

\smallskip
\noindent\emph{Case 1:
\(y\in\operatorname{int}\mathcal K_{\exp}\).}
Lemma~\ref{lem:dphi} gives
\(\mathcal R_\Phi(y)=\mathbb R^3
=T_{\mathcal K_{\exp}}(y)\).
Taking \(w=0\) at any preimage therefore proves
the reverse inclusion.

\smallskip
\noindent\emph{Case 2: \(y\in F_{\rm curv}\).}
By Lemma~\ref{lem:preimage-structure-exp-cone},
every preimage has the form \(z=(a,0,0,r)\),
where \(a\neq0\).
Equations~\eqref{eq:D_Phi} and~\eqref{eq:D2Phi-full}
give
\[
\begin{aligned}
\ker D\Phi(z)
&=\{(0,w_b,w_c,0)\colon w_b,w_c\in\mathbb R\},\\
D^2\Phi(z)[w,w]
&=(-2w_b^2,2w_c^2,0)
\qquad (w\in\ker D\Phi(z)).
\end{aligned}
\]
Lemma~\ref{lem:dphi} and
Proposition~\ref{prop:tangent-normal-lineality-exp-cone}
give \(\mathcal R_\Phi(y)=\nu(y)^\perp\), where
\(\nu(y)=(y_2,-y_3,y_2(1-y_1/y_3))\).

Let \(v\in T_{\mathcal K_{\exp}}(y)\).
By Proposition~\ref{prop:tangent-normal-lineality-exp-cone},
one has \(\langle\nu(y),v\rangle\leq0\).
Since \(y_3>0\), choose \(w_c\geq0\) satisfying
\(2y_3w_c^2=-\langle\nu(y),v\rangle\).
For \(w=(0,0,w_c,0)\), the preceding formulas give
\(w\in\ker D\Phi(z)\) and
\(D^2\Phi(z)[w,w]=(0,2w_c^2,0)\).
Consequently,
\[
\left\langle
\nu(y),v-D^2\Phi(z)[w,w]
\right\rangle
=
\langle\nu(y),v\rangle+2y_3w_c^2
=0.
\]
Thus \(v-D^2\Phi(z)[w,w]\in\mathcal R_\Phi(y)\),
which proves the reverse inclusion.

\smallskip
\noindent\emph{Case 3: \(y_3=0\).}
We treat the four remaining subsets together.
Write \(y=(-b^2,c^2,0)\) and
\(z_r\mathrel{\coloneqq}(0,b,c,r)\in\Phi^{-1}(y)\) for \(r\in\mathbb R\).
Equation~\eqref{eq:D_Phi} gives
\[
\begin{aligned}
\mathcal R_\Phi(y)
&=\operatorname{span}\{(b,0,0),(0,c,0)\},\\
\ker D\Phi(z_r)
&=\{w\in\mathbb R^4\colon bw_b=cw_c=0\}.
\end{aligned}
\]
Since \(y\in\mathcal R_\Phi(y)\) and
\(\mathcal R_\Phi(y)\) is a subspace,
we have
\(\operatorname{span}\{y\}\subseteq\mathcal R_\Phi(y)\).
Thus, \eqref{eq:tangent-cone-convex-cone} gives
\begin{equation}
\label{eq:flat-tangent-cone-inclusion}
T_{\mathcal K_{\exp}}(y)
=
\operatorname{cl}\bigl(
\operatorname{span}\{y\}+\mathcal K_{\exp}
\bigr)
\subseteq
\operatorname{cl}\bigl(
\mathcal R_\Phi(y)+\mathcal K_{\exp}
\bigr).
\end{equation}

We next prove
\begin{equation}
\label{eq:flat-preimage-recovery-inclusion}
\mathcal R_\Phi(y)+\mathcal K_{\exp}
\subseteq
\mathcal R_\Phi(y)
+\bigcup_{z\in\Phi^{-1}(y)}
\left\{
D^2\Phi(z)[w,w] \colon w\in\ker D\Phi(z)
\right\}.
\end{equation}

Let \(u\in\mathcal K_{\exp}\).
By surjectivity of \(\Phi\), write \(u=\Phi(A,B,C,r)\).
At \(z_r=(0,b,c,r)\), with the same \(r\) as in the
representation of \(u\), choose
\(w_a=A/\sqrt{2}\) and \(w_r=0\).
Set \(w_b=B/\sqrt{2}\) if \(b=0\), and \(w_b=0\) if
\(b\neq0\). Similarly, set \(w_c=C/\sqrt{2}\) if \(c=0\),
and \(w_c=0\) if \(c\neq0\).
Then \(w\in\ker D\Phi(z_r)\).
Equation~\eqref{eq:D2Phi-flat-preimage} gives
\(D^2\Phi(z_r)[w,w]=\Phi(A,\sqrt{2}w_b,\sqrt{2}w_c,r)\).
Together with \eqref{eq:Phi-decomposition}, this yields
\[
u-D^2\Phi(z_r)[w,w]
=(-B^2+2w_b^2,\ C^2-2w_c^2,\ 0).
\]
The choices of \(w_b\) and \(w_c\) make the first coordinate
zero when \(b=0\) and the second coordinate zero when \(c=0\).
Using the formula for \(\mathcal R_\Phi(y)\) above, we obtain
\[
u-D^2\Phi(z_r)[w,w]\in\mathcal R_\Phi(y).
\]

Now take any \(h\in\mathcal R_\Phi(y)\).
Using the same \(z_r\) and \(w\), we can write
\[
h+u
=
h+\bigl(u-D^2\Phi(z_r)[w,w]\bigr)
+
D^2\Phi(z_r)[w,w].
\]
Note that \(h+\bigl(u-D^2\Phi(z_r)[w,w]\bigr) \in \mathcal R_\Phi(y)\)
because \(\mathcal R_\Phi(y)\) is a subspace closed under addition.
Since \(h\in\mathcal R_\Phi(y)\) and
\(u\in\mathcal K_{\exp}\) were arbitrary, this proves~\eqref{eq:flat-preimage-recovery-inclusion}.
Taking closures in
\eqref{eq:flat-preimage-recovery-inclusion}
and using~\eqref{eq:flat-tangent-cone-inclusion} yields
\[
\begin{aligned}
T_{\mathcal K_{\exp}}(y)
&\subseteq
\operatorname{cl}\bigl(
\mathcal R_\Phi(y)+\mathcal K_{\exp}
\bigr)\\
&\subseteq
\operatorname{cl}\Bigl(
\mathcal R_\Phi(y)
+\bigcup_{z\in\Phi^{-1}(y)}
\left\{
D^2\Phi(z)[w,w] \colon w\in\ker D\Phi(z)
\right\}
\Bigr).
\end{aligned}
\]
This proves the reverse inclusion.
Together with the first inclusion, this proves
\eqref{eq:tangent-cone-all-preimages}.
\end{proof}

The preceding lemma gives the following characterization
of normal compatibility through the AP-SOC.

\begin{proposition}
\label{prop:all-preimage-normal-characterization-one-block}
Let \(y\in\mathcal K_{\exp}\) and
\(\mu\in\mathcal A_\Phi(y)\).
Then the following statements are equivalent:
\begin{enumerate}[(i)]
\item The AP-SOC holds at
\((y,\mu)\).

\item The multiplier \(\mu\) is normal-compatible with \(y\).
\end{enumerate}
\end{proposition}

\begin{proof}
Taking polars in \eqref{eq:tangent-cone-all-preimages}
shows that, for \(\mu\in\mathcal A_\Phi(y)\), the AP-SOC
holds at \((y,\mu)\) if and only if
\(\mu\in N_{\mathcal K_{\exp}}(y)\).
The equivalence therefore follows from
Proposition~\ref{prop:normal-characterization-one-block}.
\end{proof}

By
Lemma~\ref{lem:preimage-structure-exp-cone},
\(\Phi^{-1}(y)\) is infinite for every
\(y\in\mathcal K_{\exp}\setminus F_{\rm curv}\), and
Definition~\ref{def:all-preimage-slack-condition-one-block}
requires the inequality~\eqref{eq:ap-soc-one-block}
to hold at every \(z\in\Phi^{-1}(y)\).
We therefore ask how many selected preimages suffice for
their slack-only inequalities to recover normal-compatibility
for every \(\mu\in\mathcal A_\Phi(y)\).
Conversely, normal-compatibility implies these
inequalities at every preimage, and hence at any selected
preimages, by
Proposition~\ref{prop:all-preimage-normal-characterization-one-block}.

The following one-block definition formalizes the minimum number
of selected preimages needed for this recovery.

\begin{definition}
\label{def:preimage-test-complexity}
For \(y\in\mathcal K_{\exp}\) and \(m\in\mathbb N_0\), let
\(\mathcal P_m(y)\) mean that there exist
\(z^1,\ldots,z^m\in\Phi^{-1}(y)\) such that, for every
\(\mu\in\mathcal A_\Phi(y)\), the inequalities
\begin{equation}
\label{eq:selected-preimage-slack-inequalities}
-\langle\mu,D^2\Phi(z^j)[w,w]\rangle\geq0
\qquad
\forall w\in\ker D\Phi(z^j),\quad j=1,\ldots,m
\end{equation}
imply that \(\mu\) is normal-compatible with \(y\).
Define
\[
\kappa_\Phi(y)
\coloneqq
\inf\{m\in\mathbb N_0\colon\mathcal P_m(y)\}.
\]
Here, \(\mathbb N_0\mathrel{\coloneqq}\{0,1,\ldots\}\) and 
the infimum of the empty set is \(+\infty\).
When \(m=0\), no preimage tests are imposed:
\(\mathcal P_0(y)\) holds precisely when every
\(\mu\in\mathcal A_\Phi(y)\) is normal-compatible with \(y\).
\end{definition}

The selected preimages may depend on \(y\), but the same
selection must work for every \(\mu\in\mathcal A_\Phi(y)\).
By Proposition~\ref{prop:normal-characterization-one-block},
normal-compatibility is equivalent to
membership in \(N_{\mathcal K_{\exp}}(y)\)
within \(\mathcal A_\Phi(y)\).
Thus, \(\kappa_\Phi(y)\) gives the minimum number of
selected preimages needed to recover
this normal-cone condition.
When \(\kappa_\Phi(y)=+\infty\),
for every finite selection, there exists
a multiplier in \(\mathcal A_\Phi(y)\) that satisfies all
the selected inequalities but violates the normal-cone condition.

The following theorem determines \(\kappa_\Phi(y)\)
on each of the six subsets of the exponential cone.

\begin{theorem}[{Minimum number of preimage tests}]
\label{thm:preimage-test-complexity}
For every \(y\in\mathcal K_{\exp}\), one has 
\[
\kappa_\Phi(y)
=
\begin{cases}
0,
&
y\in\operatorname{int}\mathcal K_{\exp},
\\
1,
&
y\in F_{\rm curv}\cup\operatorname{ri}F_{-\infty},
\\
+\infty,
&
y\in
\bigl(F_{\rm ne}\setminus\{0\}\bigr)
\cup
\bigl(F_\infty\setminus\{0\}\bigr)
\cup
\{0\}.
\end{cases}
\]
Moreover, when \(\kappa_\Phi(y)=1\), the slack-only inequality at
any preimage of \(y\) is sufficient.
\end{theorem}

\begin{proof}
We consider the six subsets of the exponential cone.

\smallskip
\noindent\emph{Case 1:
\(y\in\operatorname{int}\mathcal K_{\exp}\).}
By Lemma~\ref{lem:dphi},
\(\mathcal A_\Phi(y)=\{0\}\). Hence every multiplier in
\(\mathcal A_\Phi(y)\) is normal-compatible without a test, and
\(\kappa_\Phi(y)=0\).

\smallskip
\noindent\emph{Case 2: \(y\in F_{\rm curv}\).}
Set \(r\mathrel{\coloneqq}y_1/y_3\). Every preimage has the form
\(z=(a,0,0,r)\), where \(a^2=y_3>0\), and
\(\ker D\Phi(z)=\{(0,w_b,w_c,0)\colon w_b,w_c\in\mathbb R\}\).
For \(\mu\in\mathcal A_\Phi(y)\),
\eqref{eq:D_Phi} gives
\(\mu_1=-\exp(r)\mu_2\).
Equation~\eqref{eq:D2Phi-full} therefore reduces the test to
\begin{equation}
\label{eq:ap-soc-curved-one-preimage}
-2\mu_2\bigl(\exp(r)w_b^2+w_c^2\bigr)\geq0
\qquad
\forall w_b,w_c\in\mathbb R.
\end{equation}
Thus, the test at any
preimage holds if and only if \(\mu_2\leq0\), which is
normal-compatibility on \(F_{\rm curv}\). Hence
\(\mathcal P_1(y)\) holds.

Zero tests do not suffice because
\(\mu=(-\exp(r),1,\exp(r)(r-1))\) belongs to
\(\mathcal A_\Phi(y)\) but is not normal-compatible with \(y\).
Consequently, \(\kappa_\Phi(y)=1\).

\smallskip
\noindent\emph{Case 3:
\(y\in\operatorname{ri}F_{-\infty}\).}
Every preimage is of the form \(z=(0,b,c,r)\),
where \(b,c\neq0\), and
\(\ker D\Phi(z)=\{(w_a,0,0,w_r) \colon w_a,w_r\in\mathbb R\}\).
Exploiting \(\mu\in\mathcal A_\Phi(y)\) gives
\(\mu_1=\mu_2=0\), so \eqref{eq:D2Phi-full} reduces the test to
\begin{equation}
\label{eq:ap-soc-flat-one-preimage}
-2\mu_3w_a^2\geq0
\qquad
\forall w_a\in\mathbb R.
\end{equation}
Thus, the test at any
preimage holds if and only if \(\mu_3\leq0\), which is
normal-compatibility on \(\operatorname{ri}F_{-\infty}\).
Hence \(\mathcal P_1(y)\) holds.

Zero tests fail because
\((0,0,1)\in\mathcal A_\Phi(y)\) is not normal-compatible with
\(y\). Therefore, \(\kappa_\Phi(y)=1\).

For each of the remaining cases, zero tests fail: one may take
\(\mu=(1,0,0)\) in Case~4 and \(\mu=(0,0,1)\) in Cases~5 and~6.
In each case, the specified multiplier belongs to
\(\mathcal A_\Phi(y)\) but is not normal-compatible with \(y\).
It remains to show that no positive finite number of tests is
sufficient. 
Fix an arbitrary integer \(m\geq1\) and arbitrary preimages
\(z^1,\ldots,z^m\in\Phi^{-1}(y)\), and write
\(z^j=(a_j,b_j,c_j,r_j)\).
For these finitely many selected preimages, we construct
\(\mu\in\mathcal A_\Phi(y)\) satisfying
\eqref{eq:selected-preimage-slack-inequalities}
but not normal-compatible with \(y\).
Thus, these inequalities do not imply normal-compatibility.

\smallskip
\noindent\emph{Case 4:
\(y\in F_{\rm ne}\setminus\{0\}\).}
Write \(y=(0,y_2,0)\). Then
\(\mathcal A_\Phi(y)=\{\mu\in\mathbb R^3 \colon \mu_2=0\}\).
Set \(R\mathrel{\coloneqq}\max_{1\leq j\leq m}r_j\) and
\(\mu\mathrel{\coloneqq}(1,0,-R)\). Every \(w\in\ker D\Phi(z^j)\) satisfies
\(w_c=0\), and \eqref{eq:D2Phi-full} gives
\[
-\left\langle
\mu,D^2\Phi(z^j)[w,w]
\right\rangle
=
2w_b^2+2(R-r_j)w_a^2
\geq0.
\]
Thus \(\mu\) satisfies
\eqref{eq:selected-preimage-slack-inequalities}
but is not normal-compatible with \(y\), because
\(\mu_1=1\neq0\).
Since both \(m\geq1\) and the selected preimages were arbitrary,
\(\mathcal P_m(y)\) fails for every \(m\geq1\).
Together with the failure of \(\mathcal P_0(y)\), this gives
\(\kappa_\Phi(y)=+\infty\).

\smallskip
\noindent\emph{Case 5:
\(y\in F_\infty\setminus\{0\}\).}
Write \(y=(y_1,0,0)\).
By Lemma~\ref{lem:preimage-structure-exp-cone}
and \eqref{eq:D_Phi},
\(\mathcal A_\Phi(y)=\{\mu\in\mathbb R^3\colon\mu_1=0\}\),
and every \(w\in\ker D\Phi(z^j)\) satisfies \(w_b=0\).
Set \(m_0\mathrel{\coloneqq}\min_{1\leq j\leq m}\exp(r_j)>0\) and
\(\mu\mathrel{\coloneqq}(0,-1,m_0/2)\in\mathcal A_\Phi(y)\).
Hence \eqref{eq:D2Phi-full} gives
\[
-\langle\mu,D^2\Phi(z^j)[w,w]\rangle
=
2\left(\exp(r_j)-\frac{m_0}{2}\right)w_a^2
+2w_c^2
\geq0
\]
for every \(j=1,\ldots,m\) and every
\(w\in\ker D\Phi(z^j)\), since \(\exp(r_j)\geq m_0\).
Thus \eqref{eq:selected-preimage-slack-inequalities} holds,
although \(\mu_3=m_0/2>0\) violates normal-compatibility.
As in Case~4, this construction rules out every finite
selection of preimages, giving \(\kappa_\Phi(y)=+\infty\).

\smallskip
\noindent\emph{Case 6: \(y=0\).}
Here \(\mathcal A_\Phi(0)=\mathbb R^3\) and
\(\ker D\Phi(z^j)=\mathbb R^4\).
Take \(R\mathrel{\coloneqq}\max_{1\leq j\leq m}r_j\) and
\(\mu\mathrel{\coloneqq}(1,0,-R)\), as in Case~4.
For every \(j\) and every \(w\in\mathbb R^4\),
\eqref{eq:D2Phi-full} gives
\[
-\langle\mu,D^2\Phi(z^j)[w,w]\rangle
=2w_b^2+2(R-r_j)w_a^2\geq0.
\]
Thus \eqref{eq:selected-preimage-slack-inequalities} holds.
However, for \(r>R\), one has
\(\langle\mu,\Phi(1,0,0,r)\rangle=r-R>0\).
Hence \(\mu\notin\mathcal K_{\exp}^{\circ}\), so
\(\mu\) is not normal-compatible with the origin.
The same argument as in Case~4 therefore gives
\(\kappa_\Phi(0)=+\infty\).
\end{proof}

\begin{remark}[Countable families of preimage tests]
\label{rem:countable-preimage-tests}
The value \(\kappa_\Phi(y)=+\infty\) in
Theorem~\ref{thm:preimage-test-complexity} does not mean that
every preimage must be tested.
For each of the three subsets where \(\kappa_\Phi(y)=+\infty\),
we specify below a countably infinite selection of preimages,
independent of \(\mu\).
For each selection, suppose that a fixed
\(\mu\in\mathcal A_\Phi(y)\) satisfies
\begin{equation} \label{eq:soc_without_ap}
    -\langle\mu,D^2\Phi(z)[w,w]\rangle\geq0 \quad \forall w\in\ker D\Phi(z)
\end{equation}
for every selected preimage \(z\).
We show that these inequalities imply normal compatibility
with \(y\).

For \(y\in F_{\rm ne}\setminus\{0\}\), select
\(z^j\mathrel{\coloneqq}(0,0,\sqrt{y_2},j) \in\Phi^{-1}(y)\),
\(j=1,2,\ldots\).
Since \(\mu\in\mathcal A_\Phi(y)\), one has \(\mu_2=0\).
Applying \eqref{eq:soc_without_ap} at \(z=z^j\)
to \(w=(0,1,0,0)\) and \(w=(1,0,0,0)\),
and using \eqref{eq:D2Phi-full}, gives
\[
\mu_1\geq0,
\qquad
j\mu_1+\mu_3\leq0
\qquad
\forall j=1,2,\ldots.
\]
Since these integer values are unbounded above,
the selected inequalities force \(\mu_1=0\) and then
\(\mu_3\leq0\).
These restrictions, together with \(\mu_2=0\),
give normal-compatibility.

For \(y\in F_\infty\setminus\{0\}\), select
\(z^j\mathrel{\coloneqq}(0,\sqrt{-y_1},0,-j)\), \(j=1,2,\ldots\).
Since \(\mu\in\mathcal A_\Phi(y)\), one has \(\mu_1=0\).
Applying \eqref{eq:soc_without_ap} at each \(z^j\)
to \(w=(0,0,1,0)\) and \(w=(1,0,0,0)\),
and using \eqref{eq:D2Phi-full}, gives
\[\mu_2\leq0, \quad \exp(-j)\mu_2+\mu_3\leq0 \quad \forall  j=1,2,\ldots.\]
Letting \(j\to\infty\) yields \(\mu_3\leq0\), so
\(\mu\) is normal-compatible with \(y\).

At \(y=0\), select \(z_q\mathrel{\coloneqq}(0,0,0,q)\) for \(q\in\mathbb Q\), where \(\mathbb Q\) denotes the set of rational numbers.
Using \eqref{eq:D2Phi-full},
\eqref{eq:soc_without_ap} at all these selected preimages gives \[
\mu_1\geq0,
\qquad
\mu_2\leq0,
\qquad
q\mu_1+\exp(q)\mu_2+\mu_3\leq0
\quad
\forall q\in\mathbb Q.
\]
Since \(r\mapsto r\mu_1+\exp(r)\mu_2+\mu_3\) is continuous
and \(\mathbb Q\) is dense in \(\mathbb R\), the last inequality above
implies \(r\mu_1+\exp(r)\mu_2+\mu_3\leq0\)
for every \(r\in\mathbb R\).
The decomposition~\eqref{eq:Phi-decomposition} now gives
\(\langle\mu,\Phi(a,b,c,r)\rangle\leq0\) for every
\((a,b,c,r)\in\mathbb R^4\).
Surjectivity of \(\Phi\) therefore yields
\(\mu\in\mathcal K_{\exp}^{\circ}\), which is
normal-compatibility at the origin.
\end{remark}

The classification in
Theorem~\ref{thm:preimage-test-complexity}
concerns the proposed parameterization \(\Phi\).
Alternative \(C^2\) surjections onto \(\mathcal K_{\exp}\)
may permit recovery of the normal-cone condition from
finitely many selected preimages at some points where
\(\kappa_\Phi(y)=+\infty\).
We leave a systematic study of such alternative parameterizations to future work.

The result in
Theorem~\ref{thm:preimage-test-complexity}
shows both agreement with and differences from
standard squared-variable models.
For scalar squared slacks, second-order necessary
conditions at one selected preimage recover the
missing multiplier sign~\cite[Theorem~3.3]{DingWright2025}.
For the unrestricted square factorization
\(C(x)=FF^\T\) in nonlinear semidefinite programming,
second-order necessary conditions at any one square
factor recover positive semidefiniteness of the
multiplier~\cite[Theorem~3.1]{DingWright2025NSDP}.
For the exponential-cone reformulation, one preimage
likewise suffices on the curved boundary and the relative
interior of the flat face. On the three remaining boundary
subsets, however, no finite selection suffices to recover
the normal-cone condition for every
\(\mu\in\mathcal A_\Phi(y)\).

We next extend the AP-SOC and its characterization
of normal-compatibility to the product cone.

\begin{definition}
\label{def:all-preimage-slack-condition-product}
Let \(y\in\mathcal K_{\exp}^p\) and
\(\mu\in\mathbb R^{3p}\). We say that the
\emph{blockwise AP-SOC} holds at
\((y,\mu)\) if the AP-SOC holds at \(\bigl(y^{(i)},\mu^{(i)}\bigr)\)
for every \(i=1,\ldots,p\).
\end{definition}

The blockwise condition is equivalent to
\begin{equation}
\label{eq:ap-soc-product}
-\sum_{i=1}^p
\left\langle
\mu^{(i)},
D^2\Phi\bigl(z^{(i)}\bigr)
\bigl[w^{(i)},w^{(i)}\bigr]
\right\rangle
\geq0
\end{equation}
for every \(z\in\Phi_p^{-1}(y)\) and every
\(w\in\ker D\Phi_p(z)\).
Indeed, since \(D\Phi_p(z)\) is block diagonal,
\(w\in\ker D\Phi_p(z)\) implies
\(w^{(i)}\in\ker D\Phi(z^{(i)})\) for every \(i\).
Thus the blockwise conditions imply
\eqref{eq:ap-soc-product} by summation.
Conversely, fix a block \(i\) and arbitrary
\(z^{(i)}\in\Phi^{-1}(y^{(i)})\) and
\(w^{(i)}\in\ker D\Phi(z^{(i)})\).
Choose preimages for the remaining blocks to obtain
\(z\in\Phi_p^{-1}(y)\), and set \(w^{(j)}=0\) for every
\(j\neq i\).
Block diagonality then gives \(w\in\ker D\Phi_p(z)\).
Applying \eqref{eq:ap-soc-product}, all terms except the
\(i\)-th vanish, giving the one-block inequality.
Since \(i\), \(z^{(i)}\), and \(w^{(i)}\) were arbitrary,
the blockwise condition follows.

Applying the one-block characterization to each block gives
the following equivalence for the multiplier restrictions for NECO.

\begin{corollary}
\label{cor:all-preimage-normal-characterization-product}
Let \(y\in\mathcal K_{\exp}^p\) and
\(\mu\in\mathcal A_{\Phi_p}(y)\).
Then the following statements are equivalent:
\begin{enumerate}[(i)]
\item The blockwise AP-SOC
holds at \((y,\mu)\).
\item The multiplier \(\mu\) is blockwise normal-compatible with
\(y\).
\end{enumerate}
\end{corollary}

\begin{proof}
By the product structure,
\(\mu\in\mathcal A_{\Phi_p}(y)\) is equivalent to
\(\mu^{(i)}\in\mathcal A_\Phi(y^{(i)})\) for every
\(i=1,\ldots,p\). Both assertions are also blockwise. The result
therefore follows by applying Proposition~
\ref{prop:all-preimage-normal-characterization-one-block} to every
block.
\end{proof}

Together with Theorem~\ref{thm:rnlp-neco-kkt-characterization},
Corollary~\ref{cor:all-preimage-normal-characterization-product}
shows that a KKT point of RNLP gives a KKT point of NECO with the same
multiplier if and only if the blockwise AP-SOC holds.
When every cone block lies in a subset where at most one
preimage is needed, the RNLP SONC itself supplies this condition.

\begin{corollary}
\label{cor:rnlp-sonc-one-preimage-recovery}
If the RNLP SONC holds at
\((\bar x,\bar z,\mu)\), and
\[
g^{(i)}(\bar x)\in
\operatorname{int}\mathcal K_{\exp}
\cup F_{\rm curv}
\cup\operatorname{ri}F_{-\infty}
\qquad
\forall i=1,\ldots,p,
\]
then the blockwise AP-SOC holds at
\((g(\bar x),\mu)\), and the pair \((\bar x,\mu)\) is a KKT point
of \eqref{prob:neco}.
\end{corollary}

\begin{proof}
By Definition~\ref{def:rnlp-sonc},
\((\bar x,\bar z,\mu)\) is a KKT point of
\eqref{prob:RNLP}. Feasibility for RNLP and
\eqref{eq:rnlp-kkt-z}, together with
\eqref{eq:adjoint-space-equivalence}, give
\(\mu\in\mathcal A_{\Phi_p}(g(\bar x))\).
Hence \(\mu^{(i)}\in\mathcal A_\Phi(g^{(i)}(\bar x))\)
for every \(i=1,\ldots,p\).

Fix \(i\in\{1,\ldots,p\}\) and
\(w^{(i)}\in\ker D\Phi(\bar z^{(i)})\).
Let \(w\in\mathbb R^{4p}\) have \(i\)-th block \(w^{(i)}\)
and all other blocks zero. Then
\((0,w)\in\mathcal S_{\rm R}(\bar x,\bar z)\), so the
RNLP SONC gives
\[
-\left\langle
\mu^{(i)},
D^2\Phi(\bar z^{(i)})[w^{(i)},w^{(i)}]
\right\rangle
\geq0.
\]
Thus, at the fixed preimage \(\bar z^{(i)}\), this inequality
holds for every \(w^{(i)}\in\ker D\Phi(\bar z^{(i)})\).

If \(g^{(i)}(\bar x)\in\operatorname{int}\mathcal K_{\exp}\),
then \(\mathcal A_\Phi(g^{(i)}(\bar x))=\{0\}\), so
\(\mu^{(i)}=0\) is normal-compatible with \(g^{(i)}(\bar x)\).
If \(g^{(i)}(\bar x)\in F_{\rm curv} \cup \operatorname{ri}F_{-\infty}\), the calculations in the proof of
Theorem~\ref{thm:preimage-test-complexity} show that the inequality
at any one preimage, together with
\(\mu^{(i)}\in\mathcal A_\Phi(g^{(i)}(\bar x))\),
implies normal compatibility.
Hence \(\mu^{(i)}\) is normal-compatible with
\(g^{(i)}(\bar x)\) in every case.

Since \(i\) was arbitrary, \(\mu\) is blockwise normal-compatible
with \(g(\bar x)\).
Corollary~\ref{cor:all-preimage-normal-characterization-product}
therefore gives the blockwise AP-SOC, while
Theorem~\ref{thm:rnlp-neco-kkt-characterization}
shows that \((\bar x,\mu)\) is a KKT point of
\eqref{prob:neco}.
\end{proof}

We return to the two running examples to illustrate
the AP-SOC characterization in
Corollary~\ref{cor:all-preimage-normal-characterization-product}
and its role in identifying KKT points of NECO through RNLP.
In the first, the KKT conditions for RNLP already imply
normal compatibility, so the AP-SOC follows automatically.
In the second, the AP-SOC recovers multiplier restrictions
missing from the KKT conditions for RNLP.
The inequality at one selected preimage, however,
may be insufficient.

\begin{examplecontinuation}{ex:running-nonconvex}{Multiplier recovery}
\label{cont:running-nonconvex-recovery}
For \(s\in\{-1,0,1\}\), consider a KKT point
\((\bar x,\bar z,\mu)\)
of~\eqref{prob:running-nonconvex-rnlp} identified in
\hyperref[cont:running-nonconvex-kkt]{the KKT analysis}, with
\(\bar x=(s,\exp(s))\) and \(\bar z=(a,0,0,s)\), where
\(a=\pm1\).
Equation~\eqref{eq:D_Phi} shows that every
\(w\in\ker D\Phi(\bar z)\) has the form
\(w=(0,w_b,w_c,0)\).
Since \(\mu_1=\exp(s)\) and \(\mu_2=-1\),
equation~\eqref{eq:D2Phi-full} gives
\[
-\langle\mu,D^2\Phi(\bar z)[w,w]\rangle
=2(\exp(s)w_b^2+w_c^2)\geq0.
\]
This inequality holds at both preimages for each \(s\), so
the AP-SOC holds at \((g(\bar x),\mu)\) for all three
values of \(s\).
It imposes no new multiplier restriction in this example,
consistent with the earlier conclusion that every
KKT point of RNLP yields a KKT point
of~\eqref{prob:running-nonconvex-neco}.
\end{examplecontinuation}

\begin{examplecontinuation}{ex:running-boundary}{Multiplier recovery}
\label{cont:running-boundary-recovery}
The \hyperref[cont:running-boundary-kkt]{KKT analysis}
shows that KKT points
of~\eqref{prob:running-boundary-rnlp} can occur only at
\(\bar x=0\) or \(\bar x=1\).
We examine the AP-SOC at these two points to determine
which corresponding pairs \((\bar x,\mu)\)
satisfy the KKT conditions
of~\eqref{prob:running-boundary-neco}.

At \(\bar x=0\), the multiplier family of RNLP is
\(\mu=(m,0,0)\), \(m\in\mathbb R\), while the KKT
conditions for NECO require \(m=0\).
We now recover this restriction through the AP-SOC.
Fix a preimage \(\bar z=(0,0,\varepsilon,r)\), where
\(\varepsilon\in\{-1,1\}\) and \(r\in\mathbb R\).
It follows from~\eqref{eq:D_Phi} that
\(\ker D\Phi(\bar z)=\{w\in\mathbb R^4\colon w_c=0\}\).
Using~\eqref{eq:D2Phi-full}, we can express the AP-SOC
at \((g(0),\mu)\) as
\[
2m(w_b^2-rw_a^2)\geq0
\qquad
\forall r\in\mathbb R,\quad
\forall (w_a,w_b)\in\mathbb R^2,
\]
with the same multiplier of RNLP \(\mu=(m,0,0)\)
for every \(r\).
Taking \(w_a=1\) and \(w_b=0\) gives
\(-2mr\geq0\) for every \(r\in\mathbb R\).
This forces \(m=0\); conversely, \(m=0\) makes
the displayed expression zero.
Thus, the AP-SOC holds if and only if \(m=0\),
recovering the KKT point
\((\bar x,\mu)=(0,0)\) of NECO identified earlier.

Here, \(g(0)=(0,1,0)\in F_{\rm ne}\setminus\{0\}\).
Theorem~\ref{thm:preimage-test-complexity}
shows that no finite selection of preimages works
for every \(\mu\in\mathcal A_\Phi(g(0))\).
In this example, the KKT conditions for RNLP further
restrict the multipliers to \(\mu=(m,0,0)\).
Even within this family, the inequality at a selected
preimage may be insufficient.
To see this, consider
\(\bar z=(0,0,\varepsilon,0)\), for which \(r=0\).
The inequality reduces to \(2mw_b^2\geq0\)
and requires only \(m\geq0\).
In particular, the multiplier of RNLP
\(\mu=(1,0,0)\), for which \(m=1\),
satisfies the inequality for every
\((w_a,w_b)\in\mathbb R^2\) at this preimage.
Keep this multiplier fixed and consider another
preimage \(\widehat z=(0,0,\varepsilon,1)\), for which \(r=1\).
The direction \(w=(1,0,0,0)\) belongs to
\(\ker D\Phi(\widehat z)\) and gives
\[
-\langle\mu,D^2\Phi(\widehat z)[w,w]\rangle=-2<0.
\]
Thus, this multiplier satisfies the inequality at
\(\bar z\) but fails it at \(\widehat z\), so the AP-SOC
does not hold.
The inequality at the selected preimage \(\bar z\)
therefore does not ensure that the corresponding
candidate \((\bar x,\mu)\) is a KKT point of NECO.

We next consider the KKT points of RNLP with \(\bar x=1\).
The \hyperref[cont:running-boundary-kkt]{KKT analysis}
gives \(\mu=(-2,2,-2)\) and
\(\bar z=(\pm1,0,0,0)\).
Since \(g(1)=(0,1,1)\in F_{\rm curv}\),
Theorem~\ref{thm:preimage-test-complexity}
shows that one selected preimage suffices for recovery.
Choose \(\bar z=(1,0,0,0)\).
For \(w=(0,0,1,0)\), equations~\eqref{eq:D_Phi}
and~\eqref{eq:D2Phi-full} give
\[
D\Phi(\bar z)w=0,
\qquad
-\langle\mu,D^2\Phi(\bar z)[w,w]\rangle=-4<0.
\]
Thus, the inequality fails at this preimage,
so the AP-SOC does not hold at \((g(1),\mu)\).
The corresponding candidate
\((\bar x,\mu)=(1,(-2,2,-2))\)
is therefore not a KKT point of NECO.

Together with the analysis at \(\bar x=0\),
this shows that the AP-SOC identifies exactly
the KKT point \((\bar x,\mu)=(0,0)\) of NECO
from the candidates obtained through RNLP.
\end{examplecontinuation}

\section[Curvature consistency and optimality conditions]
{Curvature consistency and optimality conditions}
\label{sec:curvature-optimality-comparisons}

With a KKT multiplier of NECO fixed, we now examine
the relation between the quadratic form for RNLP and the
direct second-order expression for NECO.
We establish an exact curvature identity and
use it to compare the second-order necessary and sufficient conditions for NECO and RNLP.
These comparisons determine when the second-order
optimality conditions for RNLP can rule out a point as a
local minimizer of NECO or certify a strict local minimizer.

\subsection{Direction sets and reduced curvature}

To compare the two formulations, we first identify which
directions in the original variables can be paired with
slack directions.
Let \(\bar x\in\mathcal F_{\rm NECO}\) and
\(\bar z\in\Phi_p^{-1}(g(\bar x))\).
For each direction \(d\in\mathbb R^n\) in \(x\),
define the set of slack directions by
\[
\mathcal W(\bar x,\bar z;d)
\coloneqq
\left\{
w\in\mathbb R^{4p} \colon
Dg(\bar x)d-D\Phi_p(\bar z)w=0
\right\}.
\]
Thus, \(w\in\mathcal W(\bar x,\bar z;d)\) if and only if
\((d,w)\in\mathcal S_{\rm R}(\bar x,\bar z)\).

We next identify the directions \(d\) for which
\(\mathcal W(\bar x,\bar z;d)\) is nonempty.
The image of \(D\Phi_p(\bar z)\) is
independent of the selected preimage; specifically,
\[
\operatorname{Im}D\Phi_p(\bar z)
=\mathcal R_{\Phi_p}(g(\bar x)).
\]
We therefore collect these directions in the set
\[
\mathcal S_{\Phi_p}(\bar x)
\coloneqq
\left\{
d\in\mathbb R^n\colon
Dg(\bar x)d\in
\mathcal R_{\Phi_p}\bigl(g(\bar x)\bigr)
\right\}.
\]
The following proposition states the equivalence of
these descriptions.

\begin{proposition}
\label{prop:projected-direction-preimage-independence}
For every \(\bar x\in\mathcal F_{\rm NECO}\), every
\(\bar z\in\Phi_p^{-1}(g(\bar x))\), and every
\(d\in\mathbb R^n\), the following statements are equivalent:
\begin{enumerate}[(i)]
\item One has \(d\in\mathcal S_{\Phi_p}(\bar x)\).
\item The set \(\mathcal W(\bar x,\bar z;d)\) is nonempty.
\item There exists \(w\in\mathbb R^{4p}\) such that
\((d,w)\in\mathcal S_{\rm R}(\bar x,\bar z)\).
\end{enumerate}

\end{proposition}

\begin{proof}
Statements~{\rm (ii)} and~{\rm (iii)} are equivalent by definition.
Moreover,~{\rm (ii)} holds exactly when
\(Dg(\bar x)d\in\operatorname{Im}D\Phi_p(\bar z)\).  By
Proposition~\ref{prop:dphi-image-independence}, the latter image is
\(\mathcal R_{\Phi_p}(g(\bar x))\), proving~{\rm (i)}.
\end{proof}

At a KKT point \((\bar x,\bar\lambda)\) of NECO, the following result shows that
every direction \(d\in\mathcal S_{\Phi_p}(\bar x)\)
belongs to the critical cone \(C(\bar x)\) for NECO.

\begin{lemma}
\label{lem:S-contained-critical-cone}
If \((\bar x,\bar\lambda)\) is a KKT point of
\eqref{prob:neco}, then \(\mathcal S_{\Phi_p}(\bar x)\subseteq C(\bar x)\).
\end{lemma}

\begin{proof}
Let \(d\in\mathcal S_{\Phi_p}(\bar x)\). By
Corollary~\ref{cor:dphip-lineality}, 
one has \(Dg(\bar x)d \in
\operatorname{lin}T_{\mathcal K_{\exp}^p}(g(\bar x)) \subseteq T_{\mathcal K_{\exp}^p}(g(\bar x))\).
By \eqref{eq:kkt_neco_multiplier_blockwise} and
Corollary~\ref{cor:normal-implies-dphip-orthogonality},
\(\bar\lambda\in
\mathcal A_{\Phi_p}(g(\bar x))
=\mathcal R_{\Phi_p}(g(\bar x))^\perp\).
Since \(d\in\mathcal S_{\Phi_p}(\bar x)\), one has
\(Dg(\bar x)d\in\mathcal R_{\Phi_p}(g(\bar x))\).
Hence \(\langle\bar\lambda,Dg(\bar x)d\rangle=0\).
The KKT equation~\eqref{eq:kkt_neco_stationary} therefore gives
\(\langle\nabla f(\bar x),d\rangle
=-\langle\bar\lambda,Dg(\bar x)d\rangle=0\).
Thus \(d\in C(\bar x)\).
\end{proof}

For a fixed multiplier, the direct second-order
expression for NECO depends on \(d\), while the quadratic form for RNLP
also depends on the slack direction \(w\).
To compare them, we take the infimum over the corresponding
slack directions.
For \(\lambda\in\mathbb R^{3p}\) and
\(d\in\mathcal S_{\Phi_p}(\bar x)\), define
\begin{equation*}
q_{\bar x,\bar z,\lambda}(d)
\coloneqq
\inf_{w\in\mathcal W(\bar x,\bar z;d)}
Q_{\bar x,\bar z,\lambda}(d,w).
\end{equation*}
To compare the term involving \(D^2\Phi\)
with the sigma term in the second-order conditions for NECO,
we also define
\begin{equation*}
\Gamma_{\Phi_{p}}(\bar x,\bar z,\lambda;d)
\coloneqq
\sup_{w\in\mathcal W(\bar x,\bar z;d)}
\sum_{i=1}^p
\left\langle
\lambda^{(i)},
D^2\Phi(\bar z^{(i)})
[w^{(i)},w^{(i)}]
\right\rangle.
\end{equation*}
The supremum and infimum are understood in the extended-real sense.
The Hessian term in~\eqref{eq:Q-RNLP} is independent
of \(w\), and the sum involving \(D^2\Phi\) enters with
a minus sign. Hence,
\begin{equation}
\label{eq:q-and-reformulation-curvature}
q_{\bar x,\bar z,\lambda}(d)
=
d^\T\nabla^2_{xx}L(\bar x,\lambda)d
-
\Gamma_{\Phi_{p}}(\bar x,\bar z,\lambda;d).
\end{equation}

The following lemma expresses the RNLP SONC
in terms of \(q\), so that the condition involves
only directions in the original variable \(x\).

\begin{lemma}
\label{lem:rnlp-sonc-q-equivalence}
Let \((\bar x,\bar z,\mu)\) be a KKT point of
\eqref{prob:RNLP}. Then the RNLP SONC holds at
\((\bar x,\bar z,\mu)\) if and only if
\[
q_{\bar x,\bar z,\mu}(d)\geq0
\qquad
\forall d\in\mathcal S_{\Phi_p}(\bar x).
\]
\end{lemma}

\begin{proof}
By Proposition~
\ref{prop:projected-direction-preimage-independence},
\(d\in\mathcal S_{\Phi_p}(\bar x)\) if and only if
\(\mathcal W(\bar x,\bar z;d)\neq\emptyset\). For every such \(d\),
the definition of \(q_{\bar x,\bar z,\mu}(d)\) shows that
\(q_{\bar x,\bar z,\mu}(d)\geq0\) if and only if
\(Q_{\bar x,\bar z,\mu}(d,w)\geq0\) for every
\(w\in\mathcal W(\bar x,\bar z;d)\).
Since the KKT conditions hold by assumption, this is exactly
the RNLP SONC in Definition~\ref{def:rnlp-sonc}.
\end{proof}

\subsection{Exact curvature consistency}

To compare the second-order term induced by \(\Phi\) with the
sigma term, we need explicit formulas for
\(T^2_{\mathcal K_{\exp}}(y,v)\)
when \(v\in\mathcal R_\Phi(y)\).
These are exactly the directions of the form
\(v=D\Phi(z)w\), where \(w\) is a slack direction.
The following lemma specializes~\cite[Theorem~5.6]{Xu2026Thesis}
to these directions, using Lemma~\ref{lem:dphi}
to describe \(\mathcal R_\Phi(y)\).

\begin{lemma}[{\cite[Theorem~5.6]{Xu2026Thesis}}]
\label{lem:second-order-tangent-derivative-image}
Let \(y\in\mathcal K_{\exp}\) and
\(v\in\mathcal R_\Phi(y)\). Then the following statements hold:
\begin{enumerate}[(i)]
\item If \(y\in\operatorname{int}\mathcal K_{\exp}\), then
\(T^2_{\mathcal K_{\exp}}(y,v)=\mathbb R^3\).

\item If \(y\in F_{\rm curv}\), recall from
Proposition~\ref{prop:tangent-normal-lineality-exp-cone} that
\(\nu(y)=\bigl(y_2,-y_3,y_2(1-y_1/y_3)\bigr)\).
Then
\begin{equation}
\label{eq:second-order-tangent-curved}
T^2_{\mathcal K_{\exp}}(y,v)
=
\left\{
\eta\in\mathbb R^3 \colon
\langle\nu(y),\eta\rangle
\leq
-\frac{y_2}{y_3}
\left(v_1-\frac{y_1}{y_3}v_3\right)^2
\right\}.
\end{equation}

\item If \(y\in\operatorname{ri}F_{-\infty}\), then
\(T^2_{\mathcal K_{\exp}}(y,v)
=\{\eta\in\mathbb R^3\colon \eta_3\geq0\}\).

\item If \(y\in F_{\rm ne}\setminus\{0\}\), then
\(T^2_{\mathcal K_{\exp}}(y,v)
=\{\eta\in\mathbb R^3\colon\eta_3\geq0\}\).

\item If \(y\in F_\infty\setminus\{0\}\), then
\(T^2_{\mathcal K_{\exp}}(y,v)
=\{\eta\in\mathbb R^3\colon\eta_2\geq0,\ \eta_3\geq0\}\).

\item If \(y=0\), then \(v=0\) and
\(T^2_{\mathcal K_{\exp}}(0,0)=\mathcal K_{\exp}\).
\end{enumerate}
\end{lemma}

The next proposition establishes the curvature identity for
a single exponential-cone block.
For a multiplier in the normal cone,
the sigma term equals the supremum of the second-order
term induced by \(\Phi\), over slack directions satisfying
\(D\Phi(z)w=v\).
This identity holds at every preimage of the same cone point.

\begin{proposition}
\label{prop:block-curvature-consistency}
Let \(y\in\mathcal K_{\exp}\),
\(z\in\Phi^{-1}(y)\), and
\(\lambda\in N_{\mathcal K_{\exp}}(y)\). Then, for every
\(v\in\operatorname{Im}D\Phi(z)\), one has
\begin{equation}
\label{eq:block-curvature-consistency}
\sup_{\substack{w\in\mathbb R^4\\D\Phi(z)w=v}}
\left\langle
\lambda,D^2\Phi(z)[w,w]
\right\rangle
=
\sigma\!\left(
\lambda,T^2_{\mathcal K_{\exp}}(y,v)
\right).
\end{equation}
\end{proposition}

\begin{proof}
Write \(z=(a,b,c,r)\) and fix
\(v\in\operatorname{Im}D\Phi(z)=\mathcal R_\Phi(y)\).
Thus the set of \(w\) satisfying \(D\Phi(z)w=v\) is nonempty.
We use the preimage descriptions in
Lemma~\ref{lem:preimage-structure-exp-cone}
and the derivative formulas
\eqref{eq:D_Phi} and \eqref{eq:D2Phi-full}.

Suppose that \(y\notin F_{\rm curv}\).
Since \(\lambda\in N_{\mathcal K_{\exp}}(y)\),
the normal-cone formulas in
Proposition~\ref{prop:tangent-normal-lineality-exp-cone}
and Lemma~\ref{lem:second-order-tangent-derivative-image}
give \(\langle\lambda,\eta\rangle\leq0\) for every
\(\eta\in T^2_{\mathcal K_{\exp}}(y,v)\).
Thus, by the definition of the support function,
\[
\sigma\!\left(\lambda,T^2_{\mathcal K_{\exp}}(y,v)\right)
=
\sup_{\eta\in T^2_{\mathcal K_{\exp}}(y,v)}
\langle\lambda,\eta\rangle
\leq0.
\]
Moreover, Lemma~\ref{lem:second-order-tangent-derivative-image} gives
\(0\in T^2_{\mathcal K_{\exp}}(y,v)\), so
\(\sigma(\lambda,T^2_{\mathcal K_{\exp}}(y,v))
\geq\langle\lambda,0\rangle=0\).
Consequently,
\(\sigma(\lambda,T^2_{\mathcal K_{\exp}}(y,v))=0\)
whenever \(y\notin F_{\rm curv}\).
For such cases, it remains to show that the supremum
on the left-hand side of
\eqref{eq:block-curvature-consistency} is zero.
For \(y\in F_{\rm curv}\), we compute both sides and show
that they agree.
We consider three cases, treating the four subsets
with \(y_3=0\) together.

\smallskip
\noindent\emph{Case 1:
\(y\in\operatorname{int}\mathcal K_{\exp}\).}
Here \(\lambda=0\), so the supremum is also zero.

\smallskip
\noindent\emph{Case 2: \(y\in F_{\rm curv}\).}
By Lemma~\ref{lem:preimage-structure-exp-cone},
\(b=c=0\), \(a^2=y_3>0\), and \(r=y_1/y_3\).
By Proposition~\ref{prop:tangent-normal-lineality-exp-cone},
\(\lambda=\alpha \nu(y)\) for some \(\alpha\geq0\).
By \eqref{eq:second-order-tangent-curved}, every
\(\eta\in T^2_{\mathcal K_{\exp}}(y,v)\) satisfies
\[
\alpha \langle \nu(y), \eta \rangle 
\leq-\alpha \frac{y_2 }{y_3}(v_1-rv_3)^2.
\]
This bound is attained on the boundary of the halfspace.
Hence, by the definition of the support function and
\(y_2/y_3=\exp(r)\), one has
\[
\sigma(\lambda,T^2_{\mathcal K_{\exp}}(y,v))
=-\alpha\exp(r)(v_1-rv_3)^2.
\]

Let \(w=(w_a,w_b,w_c,w_r)\) satisfy \(D\Phi(z)w=v\).
The first and third components give
\(v_1=2arw_a+a^2w_r\) and \(v_3=2aw_a\).
Hence \(w_a=v_3/(2a)\) and
\(w_r=(v_1-rv_3)/a^2\), whereas \(w_b\) and \(w_c\)
are free. Substituting into \eqref{eq:D2Phi-full} and using
\(y_2=a^2\exp(r)\), we obtain
\[
\left\langle
\lambda,D^2\Phi(z)[w,w]
\right\rangle
=
-\alpha a^4\exp(r)w_r^2
-2\alpha a^2\exp(r)w_b^2
-2\alpha a^2w_c^2.
\]
Taking \(w_b=w_c=0\) attains the supremum, which is exactly \(\sigma(\lambda,T^2_{\mathcal K_{\exp}}(y,v))\).

\smallskip
\noindent\emph{Case 3: \(y_3=0\).}
Here \(z=(0,b,c,r)\).
Lemma~\ref{lem:second-derivative-flat-preimages}
gives \(D^2\Phi(z)[w,w]\in\mathcal K_{\exp}\) for every
\(w\in\mathbb R^4\).
Since \(\lambda\in N_{\mathcal K_{\exp}}(y)
\subseteq\mathcal K_{\exp}^{\circ}\),
one has \(\langle\lambda,D^2\Phi(z)[w,w]\rangle\leq0\).

We show that zero is attained under \(D\Phi(z)w=v\).
By \eqref{eq:D_Phi},
\(\mathcal R_\Phi(y)=
\operatorname{span}\{(b,0,0),(0,c,0)\}\).
Thus \(v_3=0\), with \(v_1=0\) if \(b=0\) and
\(v_2=0\) if \(c=0\).
Choose \(w_a=w_r=0\),
\(w_b=-v_1/(2b)\) if \(b\neq0\), and
\(w_c=v_2/(2c)\) if \(c\neq0\), setting \(w_b=0\)
when \(b=0\) and \(w_c=0\) when \(c=0\).
Then \(D\Phi(z)w=v\), and
\eqref{eq:D2Phi-flat-preimage} gives
\[
D^2\Phi(z)[w,w]
=(-2w_b^2,2w_c^2,0)
\in\mathcal R_\Phi(y).
\]
Corollary~\ref{cor:normal-implies-dphip-orthogonality}
with \(p=1\) gives
\(\lambda\in\mathcal R_\Phi(y)^\perp\), so
\(\langle\lambda,D^2\Phi(z)[w,w]\rangle=0\).
The supremum is therefore zero, equal to the sigma term.
This proves \eqref{eq:block-curvature-consistency}.
\end{proof}

We now extend the identity to the product cone.
For a multiplier satisfying the normal-cone condition,
the following theorem shows that \(q\) equals the direct
second-order expression for NECO at every preimage.
This identity underlies the subsequent comparison of the
second-order conditions for NECO and RNLP.

\begin{theorem}[Exact curvature consistency]
\label{thm:product-curvature-consistency}
Let \(\bar x\in\mathcal F_{\rm NECO}\),
\(\bar z\in\Phi_p^{-1}(g(\bar x))\), and
\(\lambda\in
N_{\mathcal K_{\exp}^p}(g(\bar x))\).
Then, one has
\begin{equation}
\label{eq:product-curvature-consistency}
\Gamma_{\Phi_{p}}(\bar x,\bar z,\lambda;d)
=
\sigma\!\left(
\lambda,
T^2_{\mathcal K_{\exp}^p}
\bigl(g(\bar x),Dg(\bar x)d\bigr)
\right)
\qquad
\forall d\in\mathcal S_{\Phi_p}(\bar x).
\end{equation}
Consequently, one has
\begin{equation}
\label{eq:q-theta-consistency}
q_{\bar x,\bar z,\lambda}(d)
=
\Theta_{\bar x,\lambda}(d)
\qquad
\forall d\in\mathcal S_{\Phi_p}(\bar x).
\end{equation}
\end{theorem}

\begin{proof}
Fix \(d\in\mathcal S_{\Phi_p}(\bar x)\), and set
\(y\mathrel{\coloneqq}g(\bar x)\) and \(v\mathrel{\coloneqq}Dg(\bar x)d\). Write
\(T\mathrel{\coloneqq}T^2_{\mathcal K_{\exp}^p}(y,v)\) and
\(T_i\mathrel{\coloneqq}T^2_{\mathcal K_{\exp}}(y^{(i)},v^{(i)})\)
for \(i=1,\ldots,p\).

We first express the left-hand side of
\eqref{eq:product-curvature-consistency}
as a sum of one-block sigma terms. By Proposition~
\ref{prop:projected-direction-preimage-independence},
\(\mathcal W(\bar x,\bar z;d)\) is nonempty.
Its defining equation consists of the independent block equations
\(D\Phi(\bar z^{(i)})w^{(i)}=v^{(i)}\).
Thus \(v^{(i)}\in\operatorname{Im}D\Phi(\bar z^{(i)})\)
for every \(i\). The product structure of the normal cone also
gives \(\lambda^{(i)}\in N_{\mathcal K_{\exp}}(y^{(i)})\).
Therefore, the block structure of \(D^2\Phi_p\) and
Proposition~\ref{prop:block-curvature-consistency} yield
\[
\begin{aligned}
\Gamma_{\Phi_{p}}(\bar x,\bar z,\lambda;d)
&=
\sum_{i=1}^p
\sup_{D\Phi(\bar z^{(i)})w^{(i)}=v^{(i)}}
\left\langle
\lambda^{(i)},
D^2\Phi(\bar z^{(i)})[w^{(i)},w^{(i)}]
\right\rangle\\
&=
\sum_{i=1}^p\sigma(\lambda^{(i)},T_i).
\end{aligned}
\]

We next show that the right-hand side of
\eqref{eq:product-curvature-consistency},
namely \(\sigma(\lambda,T)\), is bounded above by the  sum.
To this end, we first establish
\(T\subseteq\prod_{i=1}^pT_i\).
For \(\eta\in T\), there exist \(t_k\downarrow0\) and
\(\eta^k\to\eta\) such that
\(y+t_kv+\frac12t_k^2\eta^k\in\mathcal K_{\exp}^p\).
Taking each block in this relation gives
\(\eta^{(i)}\in T_i\), proving the inclusion.
By the definition of the support function,
\[
\sigma(\lambda,T)
\leq
\sigma\!\left(\lambda,\prod_{i=1}^pT_i\right)
=
\sum_{i=1}^p\sigma(\lambda^{(i)},T_i)
=
\Gamma_{\Phi_{p}}(\bar x,\bar z,\lambda;d).
\]

It remains to prove the opposite inequality,
\(\Gamma_{\Phi_{p}}(\bar x,\bar z,\lambda;d)
\leq\sigma(\lambda,T)\).
Let \(w\in\mathcal W(\bar x,\bar z;d)\).
Since \(\Phi_p(\bar z)=y\) and \(D\Phi_p(\bar z)w=v\),
Lemma~\ref{lem:parameterization-second-order-tangent},
applied to \(\Phi_p\), gives
\(D^2\Phi_p(\bar z)[w,w]\in T\).
In particular, \(T\) is nonempty, and the definition of the
support function yields
\[\langle\lambda,D^2\Phi_p(\bar z)[w,w]\rangle
\leq\sigma(\lambda,T).\]
Taking the supremum of the left-hand side over
\(w\in\mathcal W(\bar x,\bar z;d)\) gives
\(\Gamma_{\Phi_{p}}(\bar x,\bar z,\lambda;d)
\leq\sigma(\lambda,T)\).
The two inequalities establish
\eqref{eq:product-curvature-consistency}.

Finally, \eqref{eq:q-theta-consistency} follows from
\eqref{eq:direct-sigma-expression},
\eqref{eq:q-and-reformulation-curvature}, and
\eqref{eq:product-curvature-consistency}.
\end{proof}

Related curvature comparisons appear in the squared-slack
NSOCP and NSDP analyses
\cite[equation~(24)]{FukudaFukushima2016NSOCP},
\cite[equation~(13)]{LourencoFukudaFukushima2017NSDP}.
Those calculations are carried out within proofs of
second-order sufficient-condition transfer.
Specifically, the proofs of
\cite[Proposition~3.3]{FukudaFukushima2016NSOCP}
and \cite[Proposition~10]{LourencoFukudaFukushima2017NSDP}
assume the second-order sufficient condition for the
squared-slack reformulation, which also implies strict
complementarity.
Although the exponential cone lacks the symmetric-cone
square structure, Theorem~\ref{thm:product-curvature-consistency}
establishes its curvature identity under the normal-cone
condition on \(\lambda\), without nondegeneracy, strict
complementarity, or a second-order optimality assumption.
Thus, curvature consistency in this setting
can be established independently of the additional
assumptions used in those optimality-condition proofs.

For unrestricted square matrix factors, the lemmas in
\cite[Lemmas~2.1--2.2]{DingWright2025NSDP}
also imply a corresponding semidefinite curvature identity
without additional regularity or second-order optimality
assumptions. Their construction uses the square factorization
of a positive semidefinite matrix.
Our theorem establishes this correspondence for the
exponential cone, without using a
symmetric-cone square representation.

For fixed \(\bar x\), \(\lambda\), and
\(d\in\mathcal S_{\Phi_p}(\bar x)\), Theorem~\ref{thm:product-curvature-consistency} also shows
that \(\Gamma_{\Phi_p}\) and \(q\) are independent of the
selected preimage, even though the corresponding slack
directions may differ.
Under the theorem's assumptions, one may therefore choose
a convenient preimage when evaluating these quantities.
The equality
\(q_{\bar x,\bar z,\lambda}(d)=\Theta_{\bar x,\lambda}(d)\)
allows the second-order expression for NECO on these directions
to be obtained from the quadratic form for RNLP.
Note that this identity does not imply that the RNLP expression
is cheaper to evaluate. The explicit formulas in
Lemma~\ref{lem:second-order-tangent-derivative-image}
also permit direct evaluation of the sigma term on these
directions. The identity establishes that the slack
derivatives represent the same conic curvature; a comparison
of numerical methods requires further analysis and experiment.
We use this representation to compare necessary conditions
in the next subsection.

\subsection{Comparison of the NECO and RNLP SONCs}
\label{subsec:comparison-sonc}

We now compare the NECO SONC and the RNLP SONC.
Two differences need to be addressed in this comparison.
First, the NECO SONC involves the entire multiplier set
\(\Lambda(\bar x)\), whereas the RNLP SONC is formulated for a fixed
multiplier. Second, the NECO SONC is imposed on the critical cone
\(C(\bar x)\), whereas, after eliminating the direction \(w\), the
RNLP SONC is imposed on the preimage-independent set
\(\mathcal S_{\Phi_p}(\bar x)\).

We first isolate the exact correspondence for one fixed multiplier.

\begin{proposition}
\label{prop:fixed-multiplier-sonc-equivalence}
Let \((\bar x,\lambda)\) be a KKT point of
\eqref{prob:neco} and
\(\bar z\in\Phi_p^{-1}\bigl(g(\bar x)\bigr)\).
Then the RNLP SONC holds at
\((\bar x,\bar z,\lambda)\) if and only if
\begin{equation}
\label{eq:fixed-multiplier-sigma-inequality}
\Theta_{\bar x,\lambda}(d)\geq0
\qquad
\forall d\in\mathcal S_{\Phi_p}(\bar x).
\end{equation}
\end{proposition}

\begin{proof}
By Theorem~\ref{thm:neco-kkt-implies-rnlp-kkt},
\((\bar x,\bar z,\lambda)\) is a KKT point of
\eqref{prob:RNLP}. By Lemma~\ref{lem:rnlp-sonc-q-equivalence}, the
RNLP SONC is equivalent to
\(q_{\bar x,\bar z,\lambda}(d)\geq0\) for every
\(d\in\mathcal S_{\Phi_p}(\bar x)\). The result now follows from
Theorem~\ref{thm:product-curvature-consistency}.
\end{proof}

For a fixed KKT multiplier \(\lambda\) of NECO,
\eqref{eq:fixed-multiplier-sigma-inequality} does not involve the
preimage. Hence the RNLP SONC holds at one preimage if and only if
it holds at every preimage of \(g(\bar x)\).
The fixed-multiplier inequality
\eqref{eq:fixed-multiplier-sigma-inequality} differs from
the NECO SONC in
Definition~\ref{def:direct-sonc-neco}, which involves the supremum
over \(\Lambda(\bar x)\) and all directions in \(C(\bar x)\).

To pass from the fixed-multiplier inequality on
\(\mathcal S_{\Phi_p}(\bar x)\) to the NECO SONC on
\(C(\bar x)\), we first give a sufficient condition under which
these direction sets coincide.

\begin{lemma}
\label{lem:S-equals-critical-cone}
Suppose that strict complementarity holds at
\((\bar x,\bar\lambda)\) and that
\(g^{(i)}(\bar x)\notin F_{\rm ne}\setminus\{0\}\) for every
\(i=1,\ldots,p\). Then
\(\mathcal S_{\Phi_p}(\bar x)=C(\bar x)\).
\end{lemma}

\begin{proof}
Set \(y\coloneqq g(\bar x)\). By the KKT conditions for NECO,
\(d\in C(\bar x)\) if and only if
\(Dg(\bar x)d\in
T_{\mathcal K_{\exp}^p}(y)\cap\{\bar\lambda\}^{\perp}\).
By strict complementarity,
\(\bar\lambda\in
\operatorname{ri}N_{\mathcal K_{\exp}^p}(y)
=\operatorname{ri}\bigl(T_{\mathcal K_{\exp}^p}(y)^\circ\bigr)\).
Hence Lemma~\ref{lem:ri-polar-lineality} gives
\(T_{\mathcal K_{\exp}^p}(y)\cap\{\bar\lambda\}^{\perp}
=\operatorname{lin}T_{\mathcal K_{\exp}^p}(y)\).
Moreover, since no block of \(y\) belongs to
\(F_{\rm ne}\setminus\{0\}\),
Corollary~\ref{cor:dphip-lineality} gives
\(\operatorname{lin}T_{\mathcal K_{\exp}^p}(y)
=\mathcal R_{\Phi_p}(y)\).
Consequently, \(d\in C(\bar x)\) if and only if
\(Dg(\bar x)d\in\mathcal R_{\Phi_p}(y)\), which is
equivalent to \(d\in\mathcal S_{\Phi_p}(\bar x)\).
\end{proof}

We can now state the main comparison between the NECO and RNLP SONCs.

\begin{theorem}[{Transfer of second-order necessary conditions}]
\label{thm:sonc-comparison}
The following statements hold:
\begin{enumerate}[(i)]
\item If the NECO SONC holds at \(\bar x\) and
\eqref{prob:neco} is nondegenerate at \(\bar x\), then
\(\Lambda(\bar x)=\{\bar\lambda\}\) for some
\(\bar\lambda\), and the RNLP SONC holds at
\((\bar x,\bar z,\bar\lambda)\) for every
\(\bar z\in\Phi_p^{-1}(g(\bar x))\).

\item If the RNLP SONC holds at
\((\bar x,\bar z,\bar\lambda)\),
\((\bar x,\bar\lambda)\) is a KKT point of NECO satisfying
strict complementarity, and
\(g^{(i)}(\bar x)\notin F_{\rm ne}\setminus\{0\}\) for every
\(i=1,\ldots,p\), then the NECO SONC holds at~\(\bar x\).
\end{enumerate}
\end{theorem}

\begin{proof}
We first prove~(i).
The NECO SONC includes the KKT conditions for NECO.
By nondegeneracy for NECO,
\(\Lambda(\bar x)=\{\bar\lambda\}\).

Fix any \(\bar z\in\Phi_p^{-1}(g(\bar x))\).
By Theorem~\ref{thm:neco-kkt-implies-rnlp-kkt},
\((\bar x,\bar z,\bar\lambda)\) is a KKT point of
\eqref{prob:RNLP}.
Let \((d,w)\in\mathcal S_{\rm R}(\bar x,\bar z)\).
By Proposition~\ref{prop:projected-direction-preimage-independence},
one has \(d\in\mathcal S_{\Phi_p}(\bar x)\). By
Lemma~\ref{lem:S-contained-critical-cone}, we have
\(d\in C(\bar x)\).
Since \(\Phi_p(\bar z)=g(\bar x)\) and
\(D\Phi_p(\bar z)w=Dg(\bar x)d\),
Lemma~\ref{lem:parameterization-second-order-tangent},
applied to \(\Phi_p\), shows that
\(D^2\Phi_p(\bar z)[w,w] \in T^2_{\mathcal K_{\exp}^p}
\bigl(g(\bar x),Dg(\bar x)d\bigr)\).

Applying the NECO SONC
with \(\mathcal T(d)=\{D^2\Phi_p(\bar z)[w,w]\}\), and using
\(\Lambda(\bar x)=\{\bar\lambda\}\), we obtain
\[
Q_{\bar x,\bar z,\bar\lambda}(d,w)=
d^\T\nabla_{xx}^2L(\bar x,\bar\lambda)d
-
\langle\bar\lambda,D^2\Phi_p(\bar z)[w,w]\rangle
\geq0.
\]
Since \((d,w)\) was arbitrary, the RNLP SONC holds at
\((\bar x,\bar z,\bar\lambda)\).

We next prove~(ii). By
Lemma~\ref{lem:S-equals-critical-cone}, the assumptions give
\(\mathcal S_{\Phi_p}(\bar x)=C(\bar x)\). The RNLP SONC and
Proposition~\ref{prop:fixed-multiplier-sonc-equivalence} therefore
give \(\Theta_{\bar x,\bar\lambda}(d)\geq0\) for every
\(d\in C(\bar x)\).

Let \(d\in C(\bar x)\), and let \(\mathcal T(d)\) be any
nonempty convex subset of
\(T^2_{\mathcal K_{\exp}^p}
\bigl(g(\bar x),Dg(\bar x)d\bigr)\). By the monotonicity of the
support function, we have
\[
d^\T\nabla_{xx}^2L(\bar x,\bar\lambda)d
-
\sigma\bigl(\bar\lambda,\mathcal T(d)\bigr)
\geq
\Theta_{\bar x,\bar\lambda}(d)
\geq0.
\]
Since \(\bar\lambda\in\Lambda(\bar x)\), the supremum over
\(\Lambda(\bar x)\) is no smaller than the expression at
\(\bar\lambda\). Hence, the NECO SONC holds at \(\bar x\).
\end{proof}

The following corollary applies part~{\rm(i)} of
Theorem~\ref{thm:sonc-comparison} to local minimizers of NECO.
Under nondegeneracy for NECO, it gives the RNLP SONC as a
necessary condition, even when LICQ for RNLP fails.

\begin{corollary}
\label{cor:local-min-nondegeneracy-rnlp-sonc}
Suppose that \(\bar x\) is a nondegenerate local minimizer of
\eqref{prob:neco}.
Then \(\Lambda(\bar x)=\{\bar\lambda\}\) for some
\(\bar\lambda\), and, for every
\(\bar z\in\Phi_p^{-1}(g(\bar x))\),
the RNLP SONC holds at
\((\bar x,\bar z,\bar\lambda)\).
\end{corollary}

\begin{proof}
Nondegeneracy implies Robinson's constraint
qualification~\eqref{eq:neco-robinson-cq}.
Since \(\bar x\) is a local minimizer, the necessary-condition
result recalled after Definition~\ref{def:direct-sonc-neco}
gives the NECO SONC.
Part~(i) of Theorem~\ref{thm:sonc-comparison} gives the conclusion.
\end{proof}

Corollary~\ref{cor:local-min-nondegeneracy-rnlp-sonc}
provides a second-order necessary condition for
a point to be a local minimizer of NECO through RNLP.
The condition uses the Hessian of the Lagrangian for RNLP
on directions satisfying the linearized equality constraints,
without requiring explicit evaluation of the sigma term as well as calculating the second-order tangent sets.
Suppose that \((\bar x,\bar\lambda)\) is a KKT point of NECO
and NECO is nondegenerate at \(\bar x\).
If the RNLP SONC fails
at \((\bar x,\bar z,\bar\lambda)\) for even one simple
\(\bar z\in\Phi_p^{-1}(g(\bar x))\), then
\(\bar x\) cannot be a local minimizer of NECO.
Meanwhile, the required nondegeneracy for NECO is weaker than
LICQ for RNLP by Theorem~\ref{thm:licq-nondegeneracy}.
The \hyperref[cont:running-boundary-cq]
{constraint-qualification analysis of
Example~\ref*{ex:running-boundary}}
shows that NECO can be nondegenerate while LICQ for RNLP
fails on the nonexposed ray.

\subsection{Comparison of the NECO and RNLP SOSCs}
\label{subsec:neco-rnlp-sosc}

We now introduce and compare the second-order
sufficient conditions for NECO and RNLP.

\begin{definition}
\label{def:direct-sosc-neco}
Suppose that \(\bar x\) satisfies the KKT conditions for
\eqref{prob:neco}. We say that the
\emph{direct sigma-term second-order sufficient condition for
\eqref{prob:neco}}
(\emph{NECO SOSC}) holds at \(\bar x\) if
\[
\sup_{\lambda\in\Lambda(\bar x)}
\Theta_{\bar x,\lambda}(d)>0
\qquad
\forall d\in C(\bar x)\setminus\{0\}.
\]
\end{definition}

Under Robinson's constraint qualification~\eqref{eq:neco-robinson-cq}
and \(C^2\)-cone reducibility of
\(\mathcal K_{\exp}^p\) at \(g(\bar x)\), the NECO SOSC is
necessary and sufficient for local quadratic growth at \(\bar x\);
see \cite[Theorem~3.137] {BonnansShapiro2000} for the general conic-constraint problem and 
\cite[Theorem~5.8]{Xu2026Thesis} for the single exponential-cone case.
These two regularity assumptions concern the quadratic-growth
characterization; they are not imposed as additional assumptions in
Theorem~\ref{thm:sosc-comparison-neco-rnlp} below.

\begin{definition}
\label{def:rnlp-sosc}
Let \((\bar x,\bar z,\bar\mu)\) be a KKT point of
\eqref{prob:RNLP}. We say that the
\emph{second-order sufficient condition for \eqref{prob:RNLP}}
(\emph{RNLP SOSC}) holds at
\((\bar x,\bar z,\bar\mu)\) if
\[
Q_{\bar x,\bar z,\bar\mu}(d,w)>0
\qquad
\forall (d,w)\in
\mathcal S_{\rm R}(\bar x,\bar z)\setminus\{(0,0)\}.
\]
\end{definition}

The RNLP SOSC is the standard second-order sufficient condition
for equality-constrained nonlinear programming, which guarantees
that \((\bar x,\bar z)\) is a strict local minimizer of
\eqref{prob:RNLP}, without requiring LICQ;
see, for example, \cite{bertsekas2016,NocedalWright2006}.

The RNLP SOSC also requires strict positivity for
directions \((0,w)\) with \(w\ne0\) and
\(D\Phi_p(\bar z)w=0\), which leave \(x\) unchanged.
The following theorem shows that this requirement forces
every cone block to lie on the curved boundary and the
multiplier to satisfy strict complementarity.
It then characterizes the remaining positivity condition
in terms of the direct expression
\(\Theta_{\bar x,\bar\mu}\).

\begin{theorem}[{Characterization of RNLP SOSC}]
\label{thm:sosc-comparison-neco-rnlp}
Let \(\bar x\in\mathbb R^n\),
\(\bar z\in\Phi_p^{-1}(g(\bar x))\), and
\(\bar\mu\in\mathbb R^{3p}\). Then the following statements are
equivalent:
\begin{enumerate}[(i)]
\item The RNLP SOSC holds at
\((\bar x,\bar z,\bar\mu)\).

\item One has
\(g^{(i)}(\bar x)\in F_{\rm curv}\) for every
\(i=1,\ldots,p\),
\((\bar x,\bar\mu)\) is a KKT point of NECO satisfying
strict complementarity, and
\[
\Theta_{\bar x,\bar\mu}(d)>0
\qquad
\forall d\in C(\bar x)\setminus\{0\}.
\]
\end{enumerate}
\end{theorem}

\begin{proof}
\noindent\emph{(i) \(\Rightarrow\) (ii).}
Suppose first that the RNLP SOSC holds.
By Definition~\ref{def:rnlp-sosc},
\((\bar x,\bar z,\bar\mu)\) is a KKT point of RNLP.
Fix \(i\in\{1,\ldots,p\}\) and write
\(\bar z^{(i)}=(\bar a_i,\bar b_i,\bar c_i,\bar r_i)\).
We first show that every block \(g^{(i)}(\bar x)\) belongs to
\(F_{\rm curv}\). By the decomposition of \(\mathcal K_{\exp}\),
it suffices to exclude \(g^{(i)}_3(\bar x)=0\), i.e., \( g^{(i)}(\bar x) \in
\operatorname{ri}F_{-\infty}
\cup\bigl(F_{\rm ne}\setminus\{0\}\bigr)
\cup\bigl(F_\infty\setminus\{0\}\bigr)
\cup\{0\}\), and
\(\operatorname{int}\mathcal K_{\exp}\).

\smallskip
\noindent\emph{Case 1: \(g^{(i)}_3(\bar x)=0\).}
Lemma~\ref{lem:preimage-structure-exp-cone} gives \(\bar a_i=0\).
Choose \(w\) with \(w^{(i)}=(0,0,0,1)\) and all other blocks zero.
Then \eqref{eq:D_Phi} and \eqref{eq:D2Phi-full} give
\((0,w)\in\mathcal S_{\rm R}(\bar x,\bar z)\setminus\{(0,0)\}\)
and \(Q_{\bar x,\bar z,\bar\mu}(0,w)=0\), contradicting the RNLP SOSC.

\smallskip
\noindent\emph{Case 2:
\(g^{(i)}(\bar x)\in\operatorname{int}\mathcal K_{\exp}\).}
Lemma~\ref{lem:preimage-structure-exp-cone} and
\eqref{eq:D_Phi} give
\(\operatorname{Im}D\Phi(\bar z^{(i)})=\mathbb R^3\).
Hence,
\eqref{eq:rnlp-kkt-z} gives \(\bar\mu^{(i)}=0\).
On the other hand, one has
\(\dim\ker D\Phi(\bar z^{(i)})=1\).
Choose
\(w^{(i)}\in\ker D\Phi(\bar z^{(i)})\setminus\{0\}\)
and set all other blocks of \(w\) to zero.
Again, \((0,w)\) is a nonzero element of
\(\mathcal S_{\rm R}(\bar x,\bar z)\) with
\(Q_{\bar x,\bar z,\bar\mu}(0,w)=0\), a contradiction.

By feasibility for RNLP and the two cases above, every block lies on
\(F_{\rm curv}\).

\smallskip
We next show that strict complementarity for
\eqref{prob:neco} holds at \((\bar x,\bar\mu)\).
For each curved block, Lemma~\ref{lem:preimage-structure-exp-cone}
gives
\(\bar z^{(i)}=(\bar a_i,0,0,\bar r_i)\), where
\(\bar a_i^2=g_3^{(i)}(\bar x)>0\) and
\(\bar r_i=g_1^{(i)}(\bar x)/g_3^{(i)}(\bar x)\).
Equation~\eqref{eq:rnlp-kkt-z} gives
\(\bar\mu^{(i)}_1=-\exp(\bar r_i)\bar\mu^{(i)}_2\) and
\(\bar\mu^{(i)}_3=\exp(\bar r_i)(\bar r_i-1)
\bar\mu^{(i)}_2\).
Choose \(w\) with \(w^{(i)}=(0,0,1,0)\) and all other blocks
zero. Since \((0,w)\) is a nonzero element of
\(\mathcal S_{\rm R}(\bar x,\bar z)\), the RNLP SOSC and
\eqref{eq:D2Phi-full} give
\(Q_{\bar x,\bar z,\bar\mu}(0,w)=-2\bar\mu^{(i)}_2>0\).
Thus \(\bar\mu^{(i)}_2<0\), and
Proposition~\ref{prop:tangent-normal-lineality-exp-cone} gives
\(\bar\mu^{(i)}\in\operatorname{ri}
N_{\mathcal K_{\exp}}(g^{(i)}(\bar x))\) for every \(i\).
Together with feasibility for RNLP and \eqref{eq:rnlp-kkt-x},
this gives the KKT conditions for NECO and strict complementarity
at \((\bar x,\bar\mu)\).

\smallskip
Finally, we show that
\(\Theta_{\bar x,\bar\mu}(d)>0\) for every
\(d\in C(\bar x)\setminus\{0\}\).
Lemma~\ref{lem:S-equals-critical-cone} gives
\(\mathcal S_{\Phi_p}(\bar x)=C(\bar x)\).
Since \(Q_{\bar x,\bar z,\bar\mu}\) is a quadratic form,
compactness of the unit sphere in
\(\mathcal S_{\rm R}(\bar x,\bar z)\) and
the RNLP SOSC imply that there exists \(\alpha>0\) such that
\[
Q_{\bar x,\bar z,\bar\mu}(d,w)
\geq\alpha\bigl(\|d\|^2+\|w\|^2\bigr)
\geq\alpha\|d\|^2
\qquad
\forall(d,w)\in\mathcal S_{\rm R}(\bar x,\bar z).
\]
For \(d\in C(\bar x)\setminus\{0\}\),
Proposition~\ref{prop:projected-direction-preimage-independence}
ensures \(\mathcal W(\bar x,\bar z;d)\neq\emptyset\).
Taking the infimum over \(w\in\mathcal W(\bar x,\bar z;d)\)
and applying Theorem~\ref{thm:product-curvature-consistency}
yields
\[
\Theta_{\bar x,\bar\mu}(d)
=q_{\bar x,\bar z,\bar\mu}(d)
\geq\alpha\|d\|^2 > 0
\qquad
\forall d\in C(\bar x)\setminus\{0\}.
\]

\smallskip
\noindent\emph{(ii) \(\Rightarrow\) (i).}
Conversely, assume~{\rm (ii)}.
The KKT conditions for NECO and
Theorem~\ref{thm:neco-kkt-implies-rnlp-kkt} give the KKT conditions for RNLP at \((\bar x,\bar z,\bar\mu)\).
Lemma~\ref{lem:S-equals-critical-cone} again gives
\(\mathcal S_{\Phi_p}(\bar x)=C(\bar x)\).
Let \((d,w)\in\mathcal S_{\rm R}(\bar x,\bar z)\setminus\{(0,0)\}\).
If \(d\ne0\), then \(d\in C(\bar x)\setminus\{0\}\), and
Theorem~\ref{thm:product-curvature-consistency} gives
\[
Q_{\bar x,\bar z,\bar\mu}(d,w)
\geq q_{\bar x,\bar z,\bar\mu}(d)
=\Theta_{\bar x,\bar\mu}(d)>0.
\]
It remains to verify the RNLP SOSC for directions with
\(d=0\) and \(w\neq0\).
Since \((0,w)\in\mathcal S_{\rm R}(\bar x,\bar z)\),
we have \(D\Phi_p(\bar z)w=0\).
The block-diagonal structure of \(D\Phi_p(\bar z)\) therefore
gives \(w^{(i)}\in\ker D\Phi(\bar z^{(i)})\) for every \(i\).
Since \(g^{(i)}(\bar x)\in F_{\rm curv}\), write
\(\bar z^{(i)}=(\bar a_i,0,0,\bar r_i)\) with \(\bar a_i^2>0\).
Moreover, exploiting \(w^{(i)}\in\ker D\Phi(\bar z^{(i)})\) and~\eqref{eq:D_Phi} gives
\(w^{(i)}=(0,w_b^{(i)},w_c^{(i)},0)\).
Strict complementarity gives \(\bar\mu_2^{(i)}<0\), so
\eqref{eq:D2Phi-full} yields
\[
Q_{\bar x,\bar z,\bar\mu}(0,w)
=-2\sum_i\bar\mu^{(i)}_2
\bigl(\exp(\bar r_i)(w_b^{(i)})^2+(w_c^{(i)})^2\bigr)>0,
\]
which follows from \(w\neq0\), and hence
one of these \(b\)- or \(c\)-components is nonzero.
Thus the RNLP SOSC holds.
\end{proof}

Theorem~\ref{thm:sosc-comparison-neco-rnlp} also shows that the RNLP
SOSC is independent of the selected preimage, because
part~{\rm(ii)} depends only on \(\bar x\) and \(\bar\mu\).
This is analogous to the preimage independence of the
RNLP SONC for a fixed KKT multiplier of NECO established in
Proposition~\ref{prop:fixed-multiplier-sonc-equivalence}.

Note that part~{\rm(ii)} of
Theorem~\ref{thm:sosc-comparison-neco-rnlp}
requires the same prescribed multiplier \(\bar\mu\)
to satisfy
\(\Theta_{\bar x,\bar\mu}(d)>0\)
for every \(d\in C(\bar x)\setminus\{0\}\).
In the NECO SOSC, on the other hand, the supremum over \(\Lambda(\bar x)\)
is taken separately for each direction, so a multiplier
giving a positive value may depend on \(d\).
Thus, the inequality in part~{\rm(ii)} implies the one in the NECO SOSC, but the NECO
SOSC alone does not ensure positivity for the prescribed
\(\bar\mu\).

Under nondegeneracy, the KKT multiplier of NECO is unique, so
the positivity requirement in part~{\rm(ii)}
becomes equivalent to the NECO SOSC.
The curved-boundary and strict-complementarity requirements
remain additional conditions.
The theorem therefore yields the following
exact comparison between the NECO and RNLP SOSCs.

\begin{corollary}
\label{cor:direct-rnlp-sosc}
Let \((\bar x,\bar\lambda)\) be a KKT point of
\eqref{prob:neco},
\(\bar z\in\Phi_p^{-1}(g(\bar x))\), and suppose that
\eqref{prob:neco} is nondegenerate at \(\bar x\). Then the RNLP SOSC
holds at \((\bar x,\bar z,\bar\lambda)\) if and only if the
following statements hold:
\begin{enumerate}[(i)]
\item The NECO SOSC holds at \(\bar x\).

\item One has
\(g^{(i)}(\bar x)\in F_{\rm curv}\) for every
\(i=1,\ldots,p\).

\item Strict complementarity for \eqref{prob:neco} holds at
\((\bar x,\bar\lambda)\).
\end{enumerate}
\end{corollary}

\begin{proof}
By nondegeneracy,
\(\Lambda(\bar x)=\{\bar\lambda\}\). Consequently, the NECO SOSC is
equivalent to
\(\Theta_{\bar x,\bar\lambda}(d)>0\) for every
\(d\in C(\bar x)\setminus\{0\}\). The conclusion now follows from
Theorem~\ref{thm:sosc-comparison-neco-rnlp}.
\end{proof}

Thus the RNLP SOSC is more restrictive than the NECO SOSC: it can
hold only when every block lies on \(F_{\rm curv}\) and strict
complementarity holds.
For a cone point in
\(\mathcal K_{\exp}\setminus F_{\rm curv}\), every preimage
has other preimages arbitrarily close to it.
Thus, if even one cone block lies outside \(F_{\rm curv}\),
we can keep \(\bar x\) fixed and change only the corresponding
slack variables while preserving feasibility and the objective
value \(f(\bar x)\).
The pair \((\bar x,\bar z)\) therefore cannot be a strict
local minimizer of RNLP, so the RNLP SOSC cannot hold.
This differs from scalar squared slacks, whose preimages
at a positive value are isolated.
Similar slack variation occurs in unrestricted matrix
factorizations, since orthogonal transformations preserve
\(FF^\T\) \cite{DingWright2025NSDP}. Nondegeneracy of NECO does not ensure LICQ for RNLP,
and a strict local minimizer of NECO need not yield
a point satisfying the RNLP SOSC.
Consequently, convergence results for NLP methods
that assume these conditions cannot be applied
automatically to RNLP.

Whenever the RNLP SOSC holds, it certifies a strict local
minimizer of NECO.
Indeed, the standard NLP sufficiency theorem shows that
\((\bar x,\bar z)\) is a strict local minimizer of RNLP
without requiring LICQ
\cite[Theorem~12.6]{NocedalWright2006}.
Theorem~\ref{thm:local-minimizer-transfer}{\rm(ii)}
then shows that \(\bar x\) is a strict local minimizer
of NECO, without an additional nondegeneracy assumption.

\subsection{Solutions of the running examples}
\label{subsec:running-example-solutions}
We now use the second-order conditions for RNLP to revisit
the two running examples. The first illustrates how calculations for RNLP recover the second-order expression for NECO,
rule out a candidate, and certify local minimizers of NECO.
The second shows that RNLP SONC remains applicable without
LICQ for RNLP, while RNLP SOSC can fail at a strict local minimizer of NECO.

\begin{examplecontinuation}{ex:running-nonconvex}{Second-order conditions}
\label{cont:running-nonconvex-second-order}
The \hyperref[cont:running-nonconvex-kkt]{KKT analysis}
gives the KKT points \((\bar x,\bar z,\mu)\)
of~\eqref{prob:running-nonconvex-rnlp} with
\[
\bar x=(s,\exp(s)),\qquad
\bar z=(a,0,0,s),\qquad
\mu=(\exp(s),-1,\exp(s)(1-s)),
\]
where \(s\in\{-1,0,1\}\) and \(a=\pm1\).
To apply the transfer results, recall that
\((\bar x,\mu)\) is also a KKT point
of~\eqref{prob:running-nonconvex-neco} at each candidate.
The \hyperref[cont:running-nonconvex-cq]
{constraint qualification analysis}
established nondegeneracy for NECO at all three points,
so \(\Lambda(\bar x)=\{\mu\}\).

By~\eqref{eq:D_Phi}, the equation
\(Dg(\bar x)d=D\Phi(\bar z)w\)
shows that the directions in
\(\mathcal S_{\rm R}(\bar x,\bar z)\) are exactly
\[
d=(t,\exp(s)t),\qquad
w=(0,w_b,w_c,t),\qquad
t,w_b,w_c\in\mathbb R.
\]
Direct substitution into~\eqref{eq:Q-RNLP} yields
\begin{equation}
\label{eq:running-nonconvex-Q}
Q_{\bar x,\bar z,\mu}(d,w)
=(12s^2-4)t^2+2\exp(s)w_b^2+2w_c^2.
\end{equation}

Identity~\eqref{eq:q-theta-consistency} allows us
to obtain the second-order expression for NECO by minimizing
this quadratic form over the slack directions.
Since \(\mu_2=-1<0\), strict complementarity holds.
As \(g(\bar x)\in F_{\rm curv}\),
Lemma~\ref{lem:S-equals-critical-cone} gives
\[
C(\bar x)=\mathcal S_{\Phi_p}(\bar x)
=\{(t,\exp(s)t)\colon t\in\mathbb R\}.
\]
For fixed \(d\), the minimum of
\(Q_{\bar x,\bar z,\mu}(d,w)\)
over the slack directions is attained at \(w_b=w_c=0\).
Theorem~\ref{thm:product-curvature-consistency}
therefore gives the expression for NECO directly:
\[
q_{\bar x,\bar z,\mu}(d)=
\Theta_{\bar x,\mu}(d)
=(12s^2-4)t^2.
\]
Since
\(d^\T\nabla_{xx}^2L(\bar x,\mu)d
=(12s^2-4-\exp(s))t^2\),
identity~\eqref{eq:q-and-reformulation-curvature} also gives
\[
\Gamma_{\Phi_p}(\bar x,\bar z,\mu;d)=-\exp(s)t^2.
\]
By the same theorem, this is the sigma term for NECO.
Thus, neither the sigma term nor the second-order tangent set
needs to be calculated separately.

At \(s=0\), choosing \(t\ne0\) and \(w_b=w_c=0\)
in~\eqref{eq:running-nonconvex-Q} gives
\(Q_{\bar x,\bar z,\mu}(d,w)=-4t^2<0\),
so the RNLP SONC fails at both preimages.
By Corollary~\ref{cor:local-min-nondegeneracy-rnlp-sonc},
this failure rules out \(\bar x=(0,1)\)
as a local minimizer
of~\eqref{prob:running-nonconvex-neco}.

At \(s=\pm1\), \eqref{eq:running-nonconvex-Q} is positive on every
nonzero \((d,w)\in\mathcal S_{\rm R}(\bar x,\bar z)\).
The RNLP SOSC therefore holds,
and each corresponding pair \((\bar x,\bar z)\)
is a strict local minimizer
of \eqref{prob:running-nonconvex-rnlp}.
Theorem~\ref{thm:local-minimizer-transfer}{\rm(ii)}
then shows that the corresponding points
\(\bar x=(\pm1,\exp(\pm1))\) are strict local minimizers
of \eqref{prob:running-nonconvex-neco}.
This sufficiency argument does not require
LICQ for RNLP or nondegeneracy for NECO.
Corollary~\ref{cor:direct-rnlp-sosc}
also gives the NECO SOSC.

Finally, we determine global optimality directly
from~\eqref{prob:running-nonconvex-neco}.
Feasibility gives
\(v\geq\exp(u)\), so
\(f(u,v)=(v-\exp(u))+(u^2-1)^2\geq0\), with equality exactly
at \((-1,\exp(-1))\) and \((1,\exp(1))\).
These are therefore the two global minimizers of NECO.
\end{examplecontinuation}

\begin{examplecontinuation}{ex:running-boundary}{Second-order conditions}
\label{cont:running-boundary-second-order}
The \hyperref[cont:running-boundary-kkt]{KKT analysis}
and the
\hyperref[cont:running-boundary-recovery]{AP-SOC analysis}
leave \((\bar x,\mu)=(0,0)\) as the only KKT point
of~\eqref{prob:running-boundary-neco};
the candidate at \(\bar x=1\) has already been rejected.
We therefore focus on \(\bar x=0\).

The
\hyperref[cont:running-boundary-local-minimizers]
{local-minimizer analysis} established that this point
is a local minimizer of NECO. The
\hyperref[cont:running-boundary-cq]
{analysis related to constraint qualifications} showed that NECO is
nondegenerate, although LICQ
for~\eqref{prob:running-boundary-rnlp} fails at every preimage.
Corollary~\ref{cor:local-min-nondegeneracy-rnlp-sonc}
then guarantees the RNLP SONC at every preimage with
the unique KKT multiplier \(\mu=0\) of NECO.
We verify this condition directly.

Every preimage has the form
\(\bar z=(0,0,\varepsilon,r)\), where
\(\varepsilon\in\{-1,1\}\) and \(r\in\mathbb R\).
The equation \(Dg(0)d=D\Phi(\bar z)w\) forces
\(d=0\) and \(w_c=0\), so
\(Q_{0,\bar z,0}(0,w)=0\) throughout
\(\mathcal S_{\rm R}(0,\bar z)\).
The RNLP SONC therefore holds at every preimage.

Consequently,
\(\mathcal S_{\Phi_p}(0)=\{0\}\), whereas
\(C(0)=\mathbb R_+\).
Identity~\eqref{eq:q-theta-consistency} in
Theorem~\ref{thm:product-curvature-consistency}
equates \(q_{0,\bar z,0}(d)\) with
\(\Theta_{0,0}(d)\) for
\(d\in\mathcal S_{\Phi_p}(0)\).
Here, it therefore gives only
\(q_{0,\bar z,0}(0)=\Theta_{0,0}(0)=0\).
To check the NECO SOSC, we must also examine
the critical directions \(d>0\) directly.

The RNLP SOSC cannot hold here:
\(g(0)=(0,1,0)\in F_{\rm ne}\setminus\{0\}\),
whereas Theorem~\ref{thm:sosc-comparison-neco-rnlp}
requires the cone point to lie on \(F_{\rm curv}\).
Thus, the RNLP SOSC cannot certify this minimizer of NECO.

We now examine the second-order condition directly
for~\eqref{prob:running-boundary-neco}.
To evaluate \(\Theta_{0,0}(d)\) for
\(d\in C(0)=\mathbb R_+\), we first calculate its sigma term.
Since the multiplier is zero, it suffices to show that
\(T^2_{\mathcal K_{\exp}}(g(0),Dg(0)d)\) is nonempty.
We establish this by showing that it contains the zero vector.
For a fixed critical direction \(d\geq0\), consider
moving from \(x=0\) to \(x=td\).
Since \(g(x)=(0,1,x)\), the corresponding path is
\[
g(td)=g(0)+tDg(0)d=(0,1,td).
\]
If \(d>0\) and \(t>0\), this point belongs to the
exponential cone if and only if
\(td\exp(0/(td))=td\leq1\), which holds for all
sufficiently small \(t>0\).
If \(d=0\), the path is constant at \(g(0)\).
Thus,
\[
g(0)+tDg(0)d+\tfrac12t^2\,0
=(0,1,td)\in\mathcal K_{\exp}
\]
for all sufficiently small \(t\geq0\).
By the definition of the second-order tangent set,
the zero second-order correction is admissible, so \(0\in T^2_{\mathcal K_{\exp}}
\bigl(g(0),Dg(0)d\bigr)\).

With \(\lambda=0\), the sigma term is therefore zero
because this set is nonempty and
\(\langle0,\eta\rangle=0\) for every element \(\eta\)
of the set. Hence,
\[
\sigma\!\left(
0,T^2_{\mathcal K_{\exp}}(g(0),Dg(0)d)
\right)=0
\qquad\forall d\geq0.
\]
Moreover, \(L(x,0)=f(x)=x^2\), so
\(\nabla_{xx}^2L(0,0)=2\). Consequently,
\[
\Theta_{0,0}(d)
=d^2\nabla_{xx}^2L(0,0)-0
=2d^2.
\]
This value is positive for every nonzero
critical direction \(d\in C(0)=\mathbb R_+\).
Thus the NECO SOSC holds, even though the RNLP SOSC fails.

Finally, the original feasible set is \([0,1]\), and
\(f(x)=x^2\geq0\), with equality only at \(x=0\).
Hence \(\bar x=0\) is the unique global minimizer
and a strict local minimizer
of \eqref{prob:running-boundary-neco}.
The corresponding pairs are also global minimizers
of~\eqref{prob:running-boundary-rnlp}.
However, varying \(r\) preserves feasibility and the objective
value, so none of these pairs is a strict local
minimizer of RNLP.
\end{examplecontinuation}


\section{Discussion and conclusion}
\label{sec:interpretation}
\label{sec:conclusion}

In this work, we show how a smooth equality reformulation supports
the analysis of NECO through classical nonlinear programming
conditions.
The resulting correspondences provide explicit criteria
for using calculations for RNLP to identify KKT points of NECO
and to rule out a point as a local minimizer of NECO
or establish that it is a strict local minimizer of NECO.

Normal compatibility characterizes when a KKT
point of RNLP yields a KKT point of NECO with the same multiplier.
Its equivalent AP-SOC characterization recovers the required
multiplier restrictions through second-order inequalities
in slack directions.
The minimum number of preimages needed for uniform recovery
is \(0\), \(1\), or \(+\infty\), depending on the point of the exponential cone.
For the proposed parameterization, this classification
identifies the cone subsets on which finitely many preimages
suffice for uniform multiplier recovery.

The exact curvature identity connects the quadratic
form for RNLP with the sigma term for NECO at every preimage and
underlies the second-order comparisons.
Under nondegeneracy for NECO, RNLP SONC with the unique KKT multiplier of NECO is necessary for a local minimizer of NECO
even when LICQ for RNLP fails.
This provides a necessary condition expressed through
the Hessian of the reformulated Lagrangian, without
separately calculating the sigma term.
RNLP SOSC certifies that a point is a strict local
minimizer of NECO, and its characterization identifies
the geometric and multiplier restrictions under which
this certificate is available.

The proposed reformulation also introduces difficulties.
A KKT point of RNLP need not yield a KKT point of NECO
with the same multiplier, and LICQ for RNLP can fail
even when NECO is nondegenerate.
Moreover, whenever a cone block lies outside the curved
boundary, the corresponding slack variables can vary
arbitrarily close to a given preimage without changing
feasibility or the objective value.
Such feasible pairs cannot be strict local minimizers
of RNLP, and the standard SOSC for RNLP cannot hold.
Thus, convergence guarantees for NLP methods that require
LICQ or SOSC cannot be applied automatically to this
reformulation.

Future work will pursue two directions.
First, we will investigate how alternative smooth
parameterizations affect multiplier recovery and the
isolation of preimages, and which difficulties of the
proposed reformulation can be reduced.
This investigation must distinguish properties specific
to \(\Phi\) from intrinsic restrictions.
In particular, the derivative restriction on
\(F_{\rm ne}\setminus\{0\}\) in
Proposition~\ref{prop:smooth-map-nonexposed-ray}
cannot be removed by changing the \(C^2\) slack mapping
within the class considered there.
Second, we will extend the analysis framework to other
nonsymmetric cones like the perspective cones, investigating smooth equality
reformulations and their properties concerning KKT recovery,
curvature consistency, and second-order optimality-condition
transfer.

\section*{Acknowledgments}

This work was supported by the JST SPRING (JPMJSP2110), 
Japan Society for the Promotion of Science,
Grant-in-Aid for Scientific Research (C) (JP25K15002), 
Grant-in-Aid for Scientific Research (B) (JP25K03082).

\paragraph{Use of AI tools.}
During the preparation of this work, ChatGPT (OpenAI, GPT-6.1 Sol) was used to polish the language, assist in generating TikZ code for a figure, and explore ideas for numerical examples. The authors reviewed and edited all the content, taking responsibility for the accuracy of the final text and examples.

\bibliographystyle{plain}
\bibliography{references}

\end{document}